\documentclass[11pt]{amsart}
\usepackage{amscd} 
\usepackage[dvips]{graphicx}
\usepackage{amssymb}

\numberwithin{equation}{section} 

\newcommand{\bFFaa}{\mathbf{F_2^{13}}\text{\bf -a}}
\newcommand{\bFFab}{\mathbf{F_2^{13}}\text{\bf -b}}
\newcommand{\bFFba}{\mathbf{F_2^{23}}\text{\bf -a}}
\newcommand{\bFFbb}{\mathbf{F_2^{23}}\text{\bf -b}}
\newcommand{\bFaa}{\mathbf{F_1^{12}}\text{\bf -a}}
\newcommand{\bFab}{\mathbf{F_1^{12}}\text{\bf -b}}
\newcommand{\bFabone}{\mathbf{F_1^{12}}\text{\bf -b1}}
\newcommand{\bFabtwo}{\mathbf{F_1^{12}}\text{\bf -b2}}
\newcommand{\bFba}{\mathbf{F_1^{23}}\text{\bf -a}}
\newcommand{\bFbb}{\mathbf{F_1^{23}}\text{\bf -b}}
\newcommand{\bFbbone}{\mathbf{F_1^{23}}\text{\bf -b1}}
\newcommand{\bFbbtwo}{\mathbf{F_1^{23}}\text{\bf -b2}}

\newcommand{\bH}{\mathbf{H}}
\newcommand{\bp}{\mathbf{p}}
\newcommand{\bu}{\mathbf{u}}
\newcommand{\bv}{\mathbf{v}}
\newcommand{\bV}{\mathbf{V}}

\newcommand{\bP}{\mathbf{P}}
\newcommand{\bQ}{\mathbf{Q}}
\newcommand{\bR}{\mathbf{R}}
\newcommand{\bEtwoR}{\mathbf{E}\text{\bf -}\mathbf{2R}}
\newcommand{\bEthreeR}{\mathbf{E}\text{\bf -}\mathbf{3R}}
\newcommand{\bEoneC}{\mathbf{E}\text{\bf -}\mathbf{1C}}
\newcommand{\bEtwoC}{\mathbf{E}\text{\bf -}\mathbf{2C}}

\newcommand{\bC}{\mathbb{C}}
\newcommand{\bN}{\mathbb{N}}

\newcommand{\bZ}{\mathbb{Z}}

\newcommand{\id}{{\mathrm{id}}}
\newcommand{\sgn}{{\mathrm{sgn}}}
\newcommand{\Tab}{{\mathrm{Tab}}}
\newcommand{\Rep}{{\mathrm{Rep}}}
\newcommand{\Image}{{\mathrm{Im}}}
\newcommand{\cA}{\mathcal{A}}

\newcommand{\cY}{\mathcal{Y}}

\newcommand{\cW}{\mathcal{W}}
\newcommand{\cZ}{\mathcal{Z}}

\newcommand{\fg}{\mathfrak{g}}

\newcommand{\fP}{\mathfrak{P}}
\newcommand{\fS}{\mathfrak{S}}
\newcommand{\fsl}{\mathfrak{sl}}

\newcommand{\hg}{\hat{\fg}}
\newcommand{\chig}{\chi}

\newcommand{\Uqg}{U_q(\fg )}

\newcommand{\Uqhg}{U_q(\hg )}

\newcommand{\lprod}{\overset{\leftarrow}{\prod }}
\newcommand{\rprod}{\overset{\rightarrow}{\prod }}

\usepackage{amsthm}
\newtheorem{thm}{Theorem}[section]
\newtheorem{lem}[thm]{Lemma}
\newtheorem{prop}[thm]{Proposition}

\newtheorem{conj}[thm]{Conjecture}

\theoremstyle{definition}
\newtheorem{defn}[thm]{Definition}

\theoremstyle{remark}
\newtheorem{rem}[thm]{Remark}

\begin{document}
\title[Paths, tableaux, and $q$-characters]{
Paths, tableaux, and $q$-characters of \\
quantum affine algebras: the $C_n$ case}

\author[W.~Nakai]{}
%
\author[T.~Nakanishi]{}

\maketitle

\begin{center}
\vspace{-8pt}
{\sc Wakako Nakai\footnote{e-mail: {\tt m99013c@math.nagoya-u.ac.jp}}
and Tomoki Nakanishi\footnote{e-mail: {\tt nakanisi@math.nagoya-u.ac.jp}}}\\
{\scriptsize {\it Graduate School of Mathematics, 
Nagoya University,
Nagoya 464-8602, Japan}}
\end{center}
\vspace{13pt}

{\small 
{\sc Abstract}.
For the quantum affine algebra $U_q(\hat{\mathfrak{g}})$
with $\mathfrak{g}$ of classical type,
let  $\chi_{\lambda/\mu,a}$ be the Jacobi-Trudi type
determinant  for the generating series
of  the (supposed) $q$-characters of the  fundamental representations.
We conjecture that $\chi_{\lambda/\mu,a}$ is the $q$-character
of a certain finite dimensional representation of
$U_q(\hat{\mathfrak{g}})$.
We study the tableaux description of  $\chi_{\lambda/\mu,a}$
using the path method due to Gessel-Viennot.
It immediately reproduces the
tableau rule by Bazhanov-Reshetikhin for $A_n$ and by 
Kuniba-Ohta-Suzuki
for $B_n$.
For $C_n$, we derive the explicit tableau rule for 
skew diagrams $\lambda/\mu$ of
three rows and of two columns.
}
\vspace{11pt}


\section{Introduction}
Let $\fg $ be a simple Lie algebra over $\bC$ and $\hg$ be the 
corresponding
non-twisted affine Lie algebra. Let $\Uqhg $ be the quantum affine 
algebra,
namely,
the quantized universal enveloping algebra of $\hg$ \cite{D, Jim1}.
The {\it$q$-character\/} of $\Uqhg$,
introduced in \cite{FR1}, is an injective ring homomorphism
$$\chi_q:\Rep \, (\Uqhg )\to \bZ[Y_{i,a}^{\pm 1}]_{i=1, \dots ,n; a\in 
\bC^{\times}},$$
where $\Rep \, (\Uqhg )$ is the Grothendieck ring of the category of
the finite dimensional representations of $\Uqhg$.
Like the usual character for $\fg$,
$\chi_q(V)$ contains essential data of each representation $V$.
Also, it is a powerful tool to investigate the ring structure
of $\Rep \, (\Uqhg )$.
Unfortunately, not much is known
about the explicit formula of $\chi_q(V)$ so far.

The $q$-character is designed to be a ``universalization'' 
of the family of the transfer matrices of the solvable vertex models 
\cite{Bax} associated to various $R$-matrices \cite{Baz, Jim2, KRS, R}. 
The tableaux descriptions of
the spectra of the transfer matrices of a vertex model
associated to $\Uqhg$ were studied
in \cite{BR, KOS, KS} for $\fg$ of classical type.
Then, one can interpret their results in the context of
the $q$-character in the following way:
Let $\chig_{\lambda/\mu,a}$ be the Jacobi-Trudi determinant
(\ref{eq:det})
for the generating series of  the (supposed) $q$-characters of the
fundamental representations of $\Uqhg$,
where $\lambda/\mu$ is a skew diagram and $a\in \bC$.
For $A_n$ and $B_n$, $\chig_{\lambda/\mu,a}$
is conjectured to be the $q$-character
of the finite dimensional irreducible
representation of $\Uqhg$ associated to $\lambda/\mu$ and $a$.
The determinant $\chig_{\lambda/\mu,a}$
allows the description by the semistandard tableaux
  of shape $\lambda/\mu $
for $A_n$ \cite{BR},
and by the tableaux of shape $\lambda/\mu $ which satisfy
certain ``horizontal'' and ``vertical'' rules
similar to the rules of the semistandard tableaux
for $B_n$ \cite{KOS} (see Definition \ref{def:B-tab} for the rules).
For $C_n$ and $D_n$,
we still conjecture  (Conjecture \ref{conj:det}) that  $\chi_{\lambda/\mu,a}$
is the $q$-character of a certain, but not necessarily irreducible,
representation of $\Uqhg$.
However, the tableaux description for $\chi_{\lambda/\mu,a}$
is known only for
the basic cases,
$(\lambda,\mu)=((1^i),\phi)$ and $(\lambda,\mu)=((i),\phi)$
  \cite{KS,KOSY,FR2}.

The main purpose of the paper is to give the tableaux description
of $\chi_{\lambda/\mu,a}$  in the  $C_n$ case.

Let us preview our results and explain
what makes the tableaux description
more complicated for  $C_n$ and $D_n$
than $A_n$ and $B_n$.
To obtain the tableaux description of
$\chi_{\lambda/\mu,a}$,
we apply the paths method of  \cite{GV}.
The method was originally introduced to derive the
well-known semistandard tableaux description of the Schur function from 
the
(original) Jacobi-Trudi determinant;
but, the idea is applicable to our determinant
$\chi_{\lambda/\mu,a}$, too.
Roughly speaking, the method works as follows:
First, we express  the determinant by sequences of ``paths''.
Then, the contributions for the determinant from
the intersecting sequences of paths cancel,
and we obtain a positive sum expression of the determinant by
the nonintersecting sequences of paths.
Finally, we translate each
nonintersecting sequence of  paths into a ``tableau''; 
the definition of a path and the nonintersecting property turn into 
the horizontal and vertical rules, respectively. 
For $A_n$,
the method works perfectly, and
it immediately reproduces the result of \cite{BR} above.
For $B_n$,
though a slight modification is required, it works well, too,
and reproduces the result of \cite{KOS} above.
For $C_n$ and $D_n$,
however,
it turns out that
the contributions from the intersecting sequences of the paths
does not completely cancel out,
and we only get an {\it alternative\/} sum expression
by nonintersecting and intersecting sequences of paths.
Therefore, we need one more step to translate it into
a {\it positive\/} sum expression by tableaux,
and it can be done essentially by the inclusion-exclusion principle.
Then, due to the negative contribution in the alternative sum,
some additional rules emerge besides the horizontal and vertical rules,
which we call the {\it extra\/} rules 
(see the two-row diagram case in Section \ref{sec:bijective-T3}           
for the simplest example).
It turns out, however,
that these extra rules depend on the shape $\lambda/\mu$,
and have infinitely many variety.
This explains, at least in our point of view,
why the tableaux description for $C_n$ and $D_n$ has not been known
so far except for the basic cases.

The outline of the paper is as follows.
In Section \ref{sec:two},
we define the Jacobi-Trudi determinant
$\chig_{\lambda/\mu,a}$ (\ref{eq:det}) and
formulate our basic conjecture (Conjecture
\ref{conj:det}) that $\chig_{\lambda/\mu,a}$
is the $q$-character of an irreducible representation
of $\Uqhg$ (for $C_n$ and $D_n$, $\mu=\phi$).
In Sections \ref{sec:A} and \ref{sec:B}
we show how the Gessel-Viennot method  works well
to reproduce the results of \cite{BR} for $A_n$ and \cite{KOS}
for $B_n$.
In Section \ref{sec:C} we consider the $C_n$ case.
This is the main part of the paper.
As explained above,  the Gessel-Viennot method
only gives an alternative sum expression
  $ \chig_{\lambda/\mu,a}$ in terms of paths 
(Proposition \ref{prop:alternative-path}).
To apply the inclusion-exclusion principle,
we introduce the ``resolution'' of a transposed pair of paths,
and derive the extra rules explicitly for the skew diagram
of three rows (Theorem \ref{thm:3-row}) and
of two columns (Theorem \ref{thm:two-column}
and Conjecture \ref{conj:two-column}).

For general skew diagrams,
the extra rules have infinitely many variety,
and so far we have  not found a unified way to write them down explicitly.
However, the above examples  suggest that,
after all,
the extra rules are better described
 {\it in terms of  paths.}
We plan to study it in a separate publication.

The $D_n$ case is similar to $C_n$, and it
will be treated also in a separate publication. 

Let us briefly mention two possible applications of the results.
Firstly, the affine crystal for the Kirillov-Reshetikhin representations,
which are special cases of the representations treated here, are
highly expected but known only for basic cases (See \cite{SS}, for example).
It is interesting to study
if there is a natural affine crystal structure on our tableaux.
Secondly, our tableaux are quite compatible
with the conjectural algorithm of \cite{FM} to create the
$q$-character.
We hope that our tableaux help us
to prove the algorithm for these representations,
and also to prove Conjecture \ref{conj:det} itself.

We thank T.~Arakawa and S.~Okada for useful discussions.

\section{$q$-characters and the Jacobi-Trudi determinant}\label{sec:two}
In this section, we give
the conjecture of the Jacobi-Trudi type formula of
the $q$-characters. 
Throughout the paper, we assume that $q^k\ne 1$ for any $k \in \bZ$.

\subsection{The variable $Y_{i,a}^{\pm 1}$ and $z_{i,a}$}\label{sec:Y-T}
The $q$-character is originally described as a polynomial in 
$\bZ[Y_{i,a}^{\pm 1}]_{i=1, \dots, n; a \in \bC^{\times}}$
in \cite{FR1}, where $Y_{i,a}$ is the affinization of 
the formal exponential $y_i:=e^{\omega_i}$ of the fundamental weight $\omega_i$
in the character $\chi: \Rep \, \Uqg \to \bZ[y_i^{\pm 1}]_{i=1, \dots, n}$ 
of $\Uqg$, with the spectral parameter $a \in \bC^{\times}$. 
For simplicity, 
we write the variable $Y_{i,aq^k}^{\pm 1}$ 
\cite{FR1} in a ``logarithmic'' form as $Y_{i,a'+k}^{\pm 1}$,
where $k \in \bZ$, $a'= \log _q a\in \bC $ and $q\in \bC ^{\times}$.
In this subsection, we transform the variables 
$\{ Y_{i,a}^{\pm 1}\}_{i=1, \dots, n; a\in \bC}$
into new variables $\{ z_{i,a}\}_{i \in I; a\in \bC}$, 
which represent the monomials in the $q$-character of the first fundamental 
representation (see \eqref{eq:first}).

Set 
\begin{equation*}
\cY = 
\begin{cases}
\bZ[Y_{1,a}^{\pm 1}, Y_{2,a}^{\pm 1}, \dots , Y_{n,a}^{\pm 1}]_{a \in \bC}
, & (A_n, C_n)\\
\bZ[Y_{1,a}^{\pm 1}, Y_{2,a}^{\pm 1}, \dots , Y_{n-1,a}^{\pm 1},
Y_{n,a-1}^{\pm 1}Y_{n,a+1}^{\pm 1}]_{a \in \bC}
, & (B_n)\\
\bZ[Y_{1,a}^{\pm 1}, Y_{2,a}^{\pm 1}, \dots , Y_{n-2,a}^{\pm 1},
Y_{n-1,a}^{\pm 1}Y_{n,a}^{\pm 1}, Y_{n,a-1}^{\pm 1}Y_{n,a+1}^{\pm 1}]_{a \in \bC}
. & (D_n)\\
\end{cases}
\end{equation*}
Let $I$ be a set of letters, 
\begin{equation}\label{entries}
I=
\begin{cases}
\{1,2,\dots ,n,n+1\}, & \qquad \qquad (A_n)\\
\{1,2,\dots ,n,0,\overline{n},\dots 
,\overline{2},\overline{1}\}, 
& \qquad \qquad (B_n),\\
\{1,2,\dots ,n,\overline{n}, \dots ,\overline{2},\overline{1}\}, 
& \qquad \qquad (C_n,D_n).
\end{cases}
\end{equation}
Let $\cZ $ be the commutative ring over $\bZ$ generated by 
$\{z_{i,a}\}_{i \in I; a \in \bC}$,
with the following generating relations ($a\in \bC$) 
with $z_{0,a}=z_{\overline{0},a}=1$ in \eqref{gen} and \eqref{gen2}:
  \begin{align}
    &
    {\textstyle \prod _{k=1}^{n+1}}z_{k,a-2k} = 1, \ \ 
    & (A_n)\\
    &
    \left\{ 
      \begin{aligned}
        z_{i,a}z_{\overline{i},a-4n+4i-2} = 
        z_{i-1,a}z_{\overline{i-1},a-4n+4i-2}
        \ \ (i=2,\dots ,n), 
        &   \quad 
        \\ 
        z_{1,a}z_{\overline{1},a-4n+2} = 1, \quad 
        z_{0,a} = 
        {\textstyle \prod _{k=1}^n}z_{k,a+4n-4k}z_{\overline{k},a-4n+4k}, 
        &  \quad   
      \end{aligned}
    \right. 
    & (B_n)\\
    &
    z_{i,a}z_{\overline{i}, a-2n+2i-4} = 
    z_{i-1,a}z_{\overline{i-1}, a-2n+2i-4}
    \ \ (i=1,\dots ,n), \quad 
    & (C_n)\label{gen}\\
    &
    z_{i,a}z_{\overline{i}, a-2n+2i} = 
    z_{i-1,a}z_{\overline{i-1}, a-2n+2i}
    \ \ (i=1,\dots ,n). \quad  
    & (D_n)\label{gen2}
  \end{align}

We have
\begin{prop}\label{prop:T-Y}
$\cZ $ is isomorphic to $\cY$ as a ring.
\end{prop}

\begin{proof}
Let 
$f:\cZ \to \cY$ be a ring homomorphism
defined as follows, with $Y_{0,a}=1$, and in \eqref{eq:A-f}, $Y_{n+1,a}=1$:
  \begin{alignat}{3}
    &
    \quad \, z_{i,a} \mapsto Y_{i,a+i-1}Y_{i-1,a+i}^{-1},
    \qquad i = 1,\dots ,n+1,
    &
    \quad (A_n)   \label{eq:A-f}\\
    &
    \left\{
      \begin{aligned}
        & z_{i,a} \mapsto Y_{i,a+2i-2}Y_{i-1,a+2i}^{-1},
        \qquad i = 1,\dots ,n-1,\\
        & z_{n,a} \mapsto Y_{n,a+2n-3}Y_{n,a+2n-1}Y_{n-1,a+2n}^{-1},\\
        & z_{\overline{n},a} \mapsto Y_{n-1,a+2n-2}Y_{n,a+2n-1}^{-1}Y_{n,a+2n+1}^{-1},\\
        & z_{\overline{i},a} \mapsto Y_{i-1, a+4n-2i-2}Y_{i,a+4n-2i}^{-1},
        \qquad i = 1,\dots ,n-1,
      \end{aligned}
    \right.
    &
    \quad (B_n) \\
    & 
    \left\{
      \begin{aligned}
        & z_{i,a} \mapsto Y_{i,a+i-1}Y_{i-1,a+i}^{-1},
        \qquad i = 1,\dots ,n,\\
        & z_{\overline{i},a} \mapsto Y_{i-1, a+2n-i+2}Y_{i,a+2n-i+3}^{-1},
        \qquad i = 1,\dots ,n,
      \end{aligned}
    \right.
    &
    \quad (C_n)\\
    &
    \left\{
      \begin{aligned}
        & z_{i,a} \mapsto Y_{i,a+i-1}Y_{i-1,a+i}^{-1},
        \qquad i = 1,\dots ,n-2,\\
        & z_{n-1,a} \mapsto Y_{n,a+n-2}Y_{n-1,a+n-2}Y_{n-2,a+n-1}^{-1},\\
        & z_{n,a} \mapsto Y_{n,a+n-2}Y_{n-1,a+n}^{-1},\\
        & z_{\overline{n},a} \mapsto Y_{n-1, a+n-2}Y_{n,a+n}^{-1},\\
        & z_{\overline{n-1},a} \mapsto Y_{n-2,a+n-1}Y_{n-1,a+n}^{-1}Y_{n,a+n}^{-1},\\
        & z_{\overline{i},a} \mapsto Y_{i-1, a+2n-i-2}Y_{i,a+2n-i-1}^{-1},
        \qquad i = 1,\dots ,n-2.
      \end{aligned}
    \right.
    & \quad (D_n)
  \end{alignat}
Let $g:\cY \to \cZ$ be the ring homomorphism defined as follows ($a \in \bC$):
  \begin{align}
    & 
    \left\{
      \begin{aligned}
        & Y_{i,a}\mapsto {\textstyle \prod_{k=1}^i}z_{k,a+i-2k+1},
        \qquad i=1, \dots ,n,\\
        & Y_{i,a}^{-1}\mapsto {\textstyle \prod_{k=i+1}^{n+1}}z_{k,a+i-2k+1},
        \qquad i=1, \dots ,n,
      \end{aligned}
    \right.
    & 
    (A_n)\label{al:g-A}\\
    & 
    \left\{
      \begin{aligned}
        & Y_{i,a}\mapsto {\textstyle \prod_{k=1}^i}z_{k,a+2i-4k+2},
        \qquad i=1, \dots ,n-1,\\
        & Y_{n,a-1}Y_{n,a+1}\mapsto 
        {\textstyle \prod_{k=1}^n}z_{k,a+2n-4k+2},\\
        & Y_{n,a-1}^{-1}Y_{n,a+1}^{-1}\mapsto 
        {\textstyle \prod_{k=1}^n}z_{\overline{k},a-6n+4k},\\
        & Y_{i,a}^{-1}\mapsto 
        {\textstyle \prod_{k=1}^i}z_{\overline{k},a-4n-2i+4k},
        \qquad i=1, \dots ,n-1,
      \end{aligned}
    \right.
    & 
    (B_n)\label{al:g-B}\\
    &
    \left\{
      \begin{aligned}
        & Y_{i,a}\mapsto {\textstyle \prod_{k=1}^i}z_{k,a+i-2k+1},
        \qquad i=1, \dots ,n,\\
        & Y_{n,a}^{-1}\mapsto {\textstyle \prod_{k=1}^i}z_{\overline{k},a-2n-i+2k-3},
        \qquad i=1, \dots ,n,
      \end{aligned}
    \right.
    & 
    (C_n)\label{al:g-C}\\
    &
    \left\{
      \begin{aligned}
        & Y_{i,a}\mapsto {\textstyle \prod_{k=1}^i}z_{k,a+i-2k+1},
        \qquad i=1, \dots ,n-2,\\
        & Y_{n-1,a}Y_{n,a}\mapsto 
        {\textstyle \prod_{k=1}^{n-1}}z_{k,a+n-2k},\\
        & Y_{n,a-1}Y_{n,a+1}\mapsto 
        {\textstyle \prod_{k=1}^n}z_{k,a+n-2k+1},\\
        & Y_{n,a-1}^{-1}Y_{n,a+1}^{-1}\mapsto {\textstyle \prod_{k=1}^n}
        z_{\overline{k},a-3n+2k+1},\\
        & Y_{n-1,a}^{-1}Y_{n,a}^{-1}\mapsto {\textstyle \prod_{k=1}^{n-1}}
        z_{\overline{k},a-3n+2k+2},\\
        & Y_{i,a}^{-1}\mapsto {\textstyle \prod_{k=1}^i}z_{\overline{k},a-2n-i+2k+1},
        \qquad i=1, \dots ,n-2.
      \end{aligned}
    \right.
    & 
    (D_n)\label{al:g-D}
  \end{align}
It is easy to check that
each homomorphism is well defined and 
$f \circ g = g\circ f= \id$, so that
$f$ and $g$ are inverse to each other. 
\end{proof}

{}From now, we identify $\cY$ with $\cZ$ by the isomorphism $f$. 
Then, the $q$-character of the first fundamental representation 
$V_{\omega_1}(q^a)$ is given as \cite{FR1}
\begin{equation}\label{eq:first}
\chi_q(V_{\omega_1}(q^a))=\sum_{i \in I}z_{i,a}.
\end{equation}
\subsection{Partitions, Young diagrams, and tableaux}
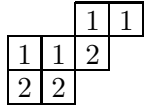
\begin{figure}[b]
{\setlength{\unitlength}{0.45mm}
\begin{picture}(40,30)
\put(0,0){\line(1,0){20}}
\put(0,10){\line(1,0){30}}
\put(0,20){\line(1,0){40}}
\put(20,30){\line(1,0){20}}
\put(0,20){\line(0,-1){20}}
\put(10,20){\line(0,-1){20}}
\put(20,30){\line(0,-1){30}}
\put(30,30){\line(0,-1){20}}
\put(40,30){\line(0,-1){10}}
\put(33,22){1}
\put(23,22){1}
\put(3,12){1}
\put(13,12){1}
\put(3,2){2}
\put(13,2){2}
\put(23,12){2}
\end{picture}
}
\caption{The highest weight tableau of $\lambda/\mu$ 
for $(\lambda, \mu) = ((4,3,2), (2))$.}\label{fig:diagrams}
\end{figure}

A {\it partition\/} 
is a sequence of weakly decreasing non-negative integers
$\lambda = (\lambda_1, \lambda_2, \dots )$
with finitely many non-zero terms
$\lambda_1\ge \lambda_2\ge \dots \ge \lambda_l > 0$.
The {\it length\/} $l(\lambda)$
of $\lambda $ is the number of the non-zero integers in $\lambda$.
The {\it conjugate\/} of $\lambda$ is denoted by 
$\lambda'=(\lambda'_1,\lambda'_2,\dots )$. 
As usual, we identify a partition $\lambda $ with a {\it Young diagram\/}
$\lambda = \{(i,j)\in \bN \times \bN \mid 1\le j\le \lambda_i\}$, 
and also identify a pair of partitions $(\lambda,\mu)$ such that $\mu \subset \lambda$,
i.e., $\lambda_i-\mu_i\ge 0$ for any $i$, 
with a {\it skew diagram\/} $\lambda/\mu = 
\{(i,j)\in\bN \times \bN \mid \mu_i+1\le j \le \lambda_i\}$.
If $\mu = \phi$, we write a skew diagram as 
a Young diagram $\lambda$ instead of $\lambda/\phi$. 
The {\it depth\/} $d(\lambda/\mu)$ of $\lambda/\mu$ is the length of its longest column,
i.e., $d(\lambda/\mu) = \max \{\lambda'_i-\mu'_i\}$.
A {\it tableau\/} $T$ of shape $\lambda/\mu$ is the skew diagram
$\lambda/\mu$ with each box filled by one entry of $I$
\eqref{entries}.

For a tableau $T$ and $a\in \bC$, we define
\begin{equation}\label{tab}
z^T_a:=\prod_{(i,j)\in \lambda/\mu}z_{T(i,j),a+2(j-i)\delta},
\end{equation}
where $T(i,j)$ is the entry of $T$ at $(i,j)$, namely, the 
entry at the $i$\,th row and the $j$\,th column, and 
$\delta$ is 
\begin{equation}\label{delta}
\delta =
\begin{cases}
1, & \qquad \qquad \qquad (A_n,C_n,D_n)\\
2. & \qquad \qquad \qquad (B_n)
\end{cases}
\end{equation}

For any skew diagram $\lambda/\mu$ with $d(\lambda/\mu)\le n$, let  
$T_+$ be the tableau of shape $\lambda/\mu$
such that $T(i,j)=i-\mu'_j$ for all $(i,j) \in \lambda/\mu$.
We call $T_+$ the {\it highest weight tableau\/} of $\lambda/\mu$. 
See Figure \ref{fig:diagrams} for example. 
Then we have 
\begin{equation}\label{eq:hw-tableau}
f(z^{T_+}_a)=\prod_{j=1}^{l(\lambda')}
Y_{\lambda'_j - \mu'_j,a(j)}^{1-\beta(j)}
Y_{n,a(j)}^{\alpha(j)}
Y_{n,a(j)-1}^{\beta(j)}Y_{n,a(j)+1}^{\beta(j)}, 
\end{equation}
where $f$ is the isomorphism in the proof of Proposition \ref{prop:T-Y} and  
\begin{align*}
a(j)= & a+(2j-\lambda'_j-\mu'_j -1)\delta, \\
\alpha(j)= & 
\begin{cases}
1, & \text{if $\fg$ is of type $D_n$ and $\lambda'_j-\mu'_j=n-1$}, \\
0, & \text{otherwise},
\end{cases}
\\
\beta(j)= & 
\begin{cases}
1, & \text{if $\fg$ is of type $B_n$ or $D_n$ and $\lambda'_j-\mu'_j=n$}, \\
0, & \text{otherwise}. 
\end{cases}
\end{align*}

\subsection{Representations of $\Uqhg$ associated to skew diagrams}
There is a bijection between the set of the isomorphism classes of
the finite dimensional irreducible representations of $\Uqhg$ 
and the set of  $n$-tuples of polynomials \cite{CP1, CP2} 
$$\bP(u) =(P_i(u))_{i=1, \dots ,n}, \qquad 
P_i(u)\in \bC[u] \text{ \ with constant term 1},$$
which are called the {\it Drinfel'd polynomials\/}. 
Let $V(\bP(u))$ be the representation associated to $\bP (u)$,
where
$$P_i(u)=\prod _{k=1}^{n_i}(1-uq^{a_{ik}}), \qquad i = 1, \dots ,n.$$
Then the $q$-character $\chi_q(V(\bP(u)))$ contains 
the {\it highest weight monomial\/} 
\begin{equation}\label{eq:hw-monomial}
m(\bP(u)):=\prod _{i=1}^n\prod_{k=1}^{n_i}Y_{i,a_{ik}}
\end{equation}
with multiplicity 1 \cite{FM}. 

For any skew diagram $\lambda/\mu$ with $d(\lambda/\mu)\le n$, 
one can uniquely associate a finite dimensional irreducible representation 
of $\Uqhg$ such that its highest weight monomial \eqref{eq:hw-monomial}
coincides with \eqref{eq:hw-tableau} for the highest weight tableau $T_+$ 
of $\lambda/\mu$. We write this representation as $V(\lambda/\mu,a)$. 
Namely, $V(\lambda/\mu,a)$ is the representation 
that corresponds to the Drinfel'd polynomial 
\begin{equation*}
\prod_{j=1}^{l(\lambda')}\bP_{\lambda'_j-\mu'_j,a(j)}^{1-\beta(j)}(u)
\bP_{n,a(j)}^{\alpha(j)}(u)
\bP_{n,a(j)-1}^{\beta(j)}(u)\bP_{n,a(j)+1}^{\beta(j)}(u), 
\end{equation*}
where $\bP\bQ:=(P_jQ_j)_{j=1, \dots , n}$ for any $\bP=(P_j)_{j=1,\dots,n}$ 
and $\bQ=(Q_j)_{j=1,\dots, n}$, 
and 
$\bP_{i,a}^{\gamma}(u)=(P_j(u))_{j=1, \dots ,n}$ is defined as 
\begin{equation*}
P_j(u)=
\begin{cases}
1-uq^a, & \text{if $j=i$ and $\gamma=1$}, \\
1, & \text{otherwise}.
\end{cases}
\end{equation*}

\subsection{The Jacobi-Trudi formula for the $q$-characters}
Let $\delta$ be the number in \eqref{delta}.
Let $\bZ [[X]]$ be the formal power series ring over $\bZ $
with variable $X$.
Let $\cA $ be the {\it non-commutative\/} ring generated by
$\cZ$ and $\bZ [[X]]$ with relations
\begin{equation}\label{gen:cA}
Xz_{i,a} = z_{i,a-2\delta}X, 
\qquad i\in I, a\in \bC.
\end{equation}
For any $a\in \bC$, we define $E_a(z,X)$, $H_a(z,X)\in \cA$ as follows:
\begin{equation}\label{eq:generating-E}
\begin{aligned}
& E_a(z,X) :=\\
& \begin{cases}
\rprod _{\scriptscriptstyle 1 \le k \le n+1}(1+z_{k,a}X) & \text{\ ($A_{n}$)}\\
\{\rprod _{\scriptscriptstyle 1 \le k \le n}(1+z_{k,a}X)\}
(1-z_{0,a}X)^{-1} 
\{\lprod _{\scriptscriptstyle 1 \le k \le n}(1+z_{\overline{k},a}X)\}
& \text{\ ($B_n$)}\\
\{\rprod _{\scriptscriptstyle 1 \le k \le n}(1+z_{k,a}X)\}
(1-z_{n,a}Xz_{\overline{n},a}X)
\{\lprod _{\scriptscriptstyle 1 \le k \le n}(1+z_{\overline{k},a}X)\}
& \text{\ ($C_n$)}\\
\{\rprod _{\scriptscriptstyle 1 \le k \le n}(1+z_{k,a}X)\}
(1-z_{\overline{n},a}Xz_{n,a}X)^{-1}
\{ \lprod _{1\scriptscriptstyle  \le k \le n}(1+z_{\overline{k},a}X)\}
& \text{\ ($D_n$)}
\end{cases}
\end{aligned}
\end{equation}
\begin{equation}\label{eq:generating-H}
\begin{aligned}
& H_a(z,X) :=\\
& 
\begin{cases}
\lprod _{\scriptscriptstyle 1 \le k \le n+1}(1-z_{k,a}X)^{-1} 
& (A_n)\\
\{\rprod _{\scriptscriptstyle 1 \le k \le n}(1-z_{\overline{k},a}X)^{-1}\}
(1+z_{0,a}X)
\{\lprod _{\scriptscriptstyle 1 \le k \le n}(1-z_{k,a}X)^{-1}\}
& (B_n)\\
\{\rprod _{\scriptscriptstyle 1 \le k \le n}(1-z_{\overline{k},a}X)^{-1}\}
(1-z_{n,a}Xz_{\overline{n},a}X)^{-1}
\{\lprod _{\scriptscriptstyle 1 \le k \le n}(1-z_{k,a}X)^{-1}\} \hspace{-5pt}
& (C_n)\\
\{\rprod _{\scriptscriptstyle 1 \le k \le n}(1-z_{\overline{k},a}X)^{-1}\}
(1-z_{\overline{n},a}Xz_{n,a}X)
\{\lprod _{\scriptscriptstyle 1 \le k \le n}(1-z_{k,a}X)^{-1}\}
& (D_n)
\end{cases}
\end{aligned}
\end{equation}
where $\rprod _{\scriptscriptstyle 1 \le k \le n}A_k=A_1\dots A_n$
and $\lprod _{\scriptscriptstyle 1 \le k \le n}A_k=A_n\dots A_1$.
Then we have
\begin{equation}\label{eq:HE}
H_a(z, X)E_a(z,-X) = E_a(z,-X)H_a(z,X)=1.
\end{equation}

For any $i \in \bZ$ and $a \in \bC$,  
we define $e_{i,a}$, $h_{i,a}\in \cZ$ as 
\begin{equation*}
E_a(z,X) = \sum _{i=0}^{\infty }e_{i,a}X^i, \qquad 
H_a(z,X) = \sum _{i=0}^{\infty}h_{i,a}X^i.
\end{equation*}
Set $e_{i,a}=h_{i,a}=0$ for $i<0$.
Note that 
$e_{i,a}=0$ if 
$i >n+1$ (resp.\ if $i >2n+2$ or $i=n+1$) 
for $A_n$ (resp.\ for $C_n$).

It has been observed in \cite{FR2, KOSY} (see also \cite{KOS,KS})
that $e_{i,a}$ is the $q$-character of the $i$\,th 
fundamental representation for $1 \le i \le n$
($i \ne n$ for $B_n$, $i \ne n-1, n$ for $D_n$), 
while 
$h_{i,a}$ is the $q$-character of the $i$\,th 
``symmetric'' power of the first fundamental representation
for any $i \ge 1$, though 
only a part of them are proven in the literature
(e.g. \cite{N}). 

Due to the relation \eqref{eq:HE}, 
it holds that \cite{M} 
\begin{equation}\label{eq:det}
\det (h_{\lambda_i-\mu_j-i+j,a+2(\lambda_i-i)\delta })_{1\le i,j \le l}
=\det (e_{\lambda'_i-\mu'_j-i+j,a-2(\mu'_j-j+1)\delta })_{1\le i,j \le l'}
\end{equation}
for any partitions $(\lambda, \mu)$, 
where $l$ and $l'$ are any non-negative integers such that
$l\ge l(\lambda),l(\mu)$ and $l'\ge l(\lambda'),l(\mu')$.
For any skew diagram $\lambda/\mu$, 
let $\chig_{\lambda /\mu ,a}$ denote the determinant 
on the left or right hand side of \eqref{eq:det}.
We call it the {\it Jacobi\/-\/Trudi determinant\/}
of $\Uqhg$ associated to $\lambda/\mu$ and $a \in \bC$. 
Note that $\chig_{(i),a}=h_{i,a}$ and $\chig_{(1^i),a}=e_{i,a}$.

\begin{conj}\label{conj:det} 
\begin{enumerate}
\item
If $\fg$ is of type $A_n$ or $B_n$ and $\lambda/\mu$ is a skew diagram 
of $d(\lambda/\mu) \le n$, then 
$\chig_{\lambda/\mu ,a}=\chi_q(V(\lambda/\mu, a))$.
\item
If $\fg$ is of type $C_n$ and 
$\lambda/\mu$ is a skew diagram of 
$d(\lambda/\mu)\le n$, 
then $\chig_{\lambda/\mu ,a}$ is the $q$-character of certain 
(not necessarily irreducible) representation $V$ of $\Uqhg$ 
which has $V(\lambda/\mu,a)$ as a subquotient; 
furthermore, if $\mu=\phi$, then $V=V(\lambda,a)$.
\item If $\fg$ is of type $D_n$ and $\lambda/\mu$ is a skew diagram of  
$d(\lambda/\mu)\le n$, 
then $\chig_{\lambda/\mu ,a}$ is the $q$-character of certain
(not necessarily irreducible) representation $V$ of $\Uqhg$
which has $V(\lambda/\mu,a)$ as a subquotient;
furthermore, if $\mu=\phi$ and $d(\lambda) \le n-1$, 
then $V=V(\lambda,a)$.
\end{enumerate}
\end{conj}

Several remarks on Conjecture \ref{conj:det} are in order. 
\newline 
1. For $C_n$,
we checked by computer that 
$\chi_{\lambda, a}$ agrees with the result obtained from
the conjectural algorithm of \cite{FM} to create the
$q$-character for several $\lambda$.
\newline
2. It is interesting that the determinant \eqref{eq:det} is simpler than the 
Jacobi-Trudi type formula for 
the characters of $\fg$ for the irreducible representations $V(\lambda)$ 
in \cite{KT}. 
\newline 
3. The determinant $\chig_{\lambda/\mu,a}$ appeared in \cite{BR} for $A_n$ 
and \cite{KOS} for $B_n$ in the context of the transfer matrices. 
\newline 
4. An analogue of 
Conjecture \ref{conj:det} is true for the representations of Yangian $Y(\fsl_n)$, 
which can be proved \cite{Ar1} using the results in 
\cite{Ar2, AS}.
\newline 
5. Conjecture \ref{conj:det} is an affinization of the conjecture 
of \cite{CK} (see Remark \ref{rem:CK} in Appendix  \ref{sec:g-ch}).
\newline 
6.  
For $C_n$ and $D_n$, we further expect that $V=V(\lambda/\mu,a)$
if $\lambda/\mu$ is connected.
But, if $\lambda/\mu$ is not connected,
 there are certainly counter-examples.
A counter-example for $C_2$ is as follows: 
Let $(\lambda, \mu) = ((3,1), (2))$.
By \eqref{eq:det}, we have
$\chig_{\lambda/\mu,a+2}
=h_{1,a}h_{1,a+6}
=\chi_q(V_{\omega_1}(q^{a})\otimes V_{\omega_1}(q^{a+6}))$.
On the other hand,
the $R$-matrix $R_{\omega_1, \omega_1}(u)$ has singularities
at $u=q^6$ (see \cite{AK} for example), which implies that 
$V_{\omega_1}(q^{a})\otimes V_{\omega_1}(q^{a+6})$ is not irreducible. 
The case $(\lambda, \mu) = ((3,1), (2))$ for $D_4$ 
is a similar counter-example.
\newline 

In the following sections,  
we study the explicit description of $\chig_{\lambda/\mu, a}$ by tableaux. 

\section{Tableaux description of type $A_n$}\label{sec:A}
In this section, we consider the case that $\fg$ is of type $A_n$.
The tableaux description of $\chi_{\lambda/\mu,a}$  \eqref{eq:det}
is given by \cite{BR}. 
We reproduce it 
by applying the ``paths'' method of \cite{GV} (see also \cite{S}).  
During this section, $I$ is of type $A_n$ in \eqref{entries}.

\subsection{Paths description}
\begin{figure}
\setlength{\unitlength}{1mm}
{\small
\begin{picture}(45,30)(-20,0)
\multiput(-20,20)(10,0){5}{\circle*{1}}
\multiput(-20,10)(10,0){5}{\circle*{1}}
\multiput(-20,0)(10,0){5}{\circle*{1}}
\multiput(-20,-10)(10,0){5}{\circle*{1}}
\multiput(-20,-20)(10,0){5}{\circle*{1}}
\put(-30,-23){$u=(0,0)$}
\put(19,22){$v$}
\put(-20,-20){\line(0,1){10}}
\put(-20,-10){\line(1,0){10}}
\put(-10,-10){\line(1,0){10}}
\put(0,-10){\line(0,1){10}}
\put(0,0){\line(1,0){10}}
\put(10,0){\line(1,0){10}}
\put(20,0){\line(0,1){10}}
\put(20,10){\line(0,1){10}}
\put(-30,0){$p=$}
\put(-24,-16){$s_1$}
\put(-17,-9){$s_2$}
\put(-7,-9){$s_3$}
\put(-4,-6){$s_4$}
\put(3,1){$s_5$}
\put(13,1){$s_6$}
\put(16,4){$s_7$}
\put(16,14){$s_8$}
\put(-30,-33){$L^2_a(s_i)$}
\put(28,25){$L^1_a(s_i)$ \qquad height}
\put(-15,-30){$\vdots$}
\put(-5,-30){$\vdots$}
\put(5,-30){$\vdots$}
\put(15,-30){$\vdots$}
{\small
\put(-15,-33){$a$}
\put(-8,-33){$a+2$}
\put(2,-33){$a+4$}
\put(12,-33){$a+6$}
}
\put(27,-21){$\cdots $1 \qquad \qquad 0}
\put(27,-11){$\cdots $2 \qquad \qquad 1}
\put(27,-1){$\cdots $3 \qquad \qquad 2}
\put(27,9){$\cdots $4 \qquad \qquad 3}
\put(27,19){$\cdots $5 \qquad \qquad 4}

\end{picture}
}
\vspace{3cm}
\caption{An example of a path $p$ and its $h$-labeling.}\label{fig:path}
\end{figure}
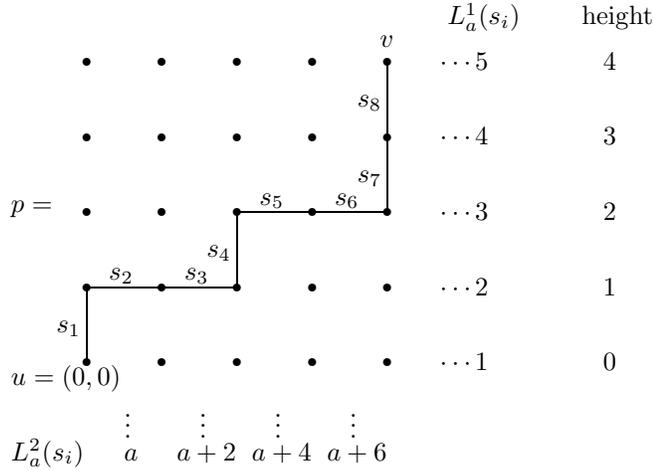

Consider the lattice $\bZ \times \bZ $.
A {\it path\/} $p$ in the lattice is a sequence of steps 
$(s_1, s_2, \dots )$ such that 
each step $s_i$ is of unit length with
the northward $(N)$ or eastward $(E)$ direction.
For example, see Figure \ref{fig:path}.
If $p$ starts at point $u$ and ends at point $v$, 
we write this by $u\overset{p}{\to}v$.
For any path $p$,
set $E(p):=\{s \in p \mid  \text{$s$ is an eastward step}\}$.

An {\it $h$-path of type $A_n$} is a path $u\overset{p}{\to}v$
such that the initial point $u$ is at height 0
and the final point $v$ is at height $n$,
where the {\it height\/} of the point $(x,y) \in \bZ \times \bZ $ is $y$.
Let $P(A_n)$ be the set of all the $h$-paths of type $A_n$. 
For any $a \in \bC$, 
the {\it $h$-labeling of type $A_n$} associated to $a\in \bC$ for a path $p\in P(A_n)$
is the pair of maps $L_a=(L_a^1, L_a^2)$,
\begin{align*}
L^1_a:E(p) \to I, \qquad 
L^2_a:E(p) \to \{a+2k\mid k \in \bZ \},
\end{align*}
defined as follows: If $s$ starts at the point $(x,y)$, then
$L_a(s)=(y+1, a+2x)$.
For example, $L_a(s_3)=(2, a+2)$
for $s_3$ in Figure \ref{fig:path}. 
Using these definitions, we define
\begin{equation}\label{eq:z-path}
z^p_a = \prod _{s\in E(p)}
z_{\scriptscriptstyle L^1_a(s),L^2_a(s)} \in \cZ
\end{equation}
for any $p\in P(A_n)$, where $\cZ$ is the ring defined in Section \ref{sec:two}.
For example, 
$z^p_a =z_{2,a}z_{2,a+2}z_{3,a+4}z_{3,a+6}$ for $p$ in Figure \ref{fig:path}.
By \eqref{eq:generating-H}, we have
\begin{equation}\label{eq:h-path-A}
h_{r,a+2k+2r-2}(z)
=\sum_{p}z^p_a,
\end{equation}
where the sum runs over all $p\in P(A_n)$ such that
$(k,0)\overset{p}{\to}(k+r,n)$.

For any $l$-tuples of initial points 
$\bu =(u_1, u_2, \dots ,u_l)$ 
and final points $\bv =(v_1, v_2, \dots , v_l)$, 
let $\fP(\pi;\bu ,\bv )$ be the set of $l$-tuples of paths
$\bp  = (p_1, \dots , p_l)$
such that
$u_i\overset{p_i}{\to}v_{\pi (i)}$ for any permutation $\pi \in \fS _l$.  
Set 
$${\textstyle{\fP(\bu ,\bv ):=\sum_{\pi \in \fS_l}\fP(\pi;\bu ,\bv )}}.$$
Then we define  
$$\fP(A_n; \bu, \bv):=\{\bp=(p_1, \dots, p_l) \in \fP(\bu, \bv) \mid p_i \in P(A_n)\}. $$

For any skew diagram $\lambda/\mu$, 
let $l=l(\lambda)$, and 
pick $\bu _{\mu}=(u_1, \dots ,u_l)$ and 
$\bv _{\lambda}=(v_1,\dots ,v_l)$
as $u_i=(\mu_i+1-i,0)$ and $v_i=(\lambda_i+1-i, n)$. 
In this case, 
we have $\fP(A_n; \bu_{\mu}, \bv_{\lambda})=\fP(\bu_{\mu} ,\bv_{\lambda})$. 
We define the {\it weight\/} $z^{\bp}_a$ 
and the {\it signature\/} $(-1)^{\bp}$ for any 
$\bp \in \fP(A_n; \bu_{\mu} ,\bv_{\lambda})$ by 
\begin{equation}\label{eq:z-bp}
z^{\bp}_a =\prod_{i=1}^lz^{p_i}_a\qquad \text{and }\qquad 
(-1)^{\bp}  =\sgn \,\pi \text{\quad if $\bp \in \fP(\pi;\bu_{\mu} ,\bv_{\lambda} )$}.
\end{equation}

Then, the determinant \eqref{eq:det} can be written as
\begin{equation}\label{eq:det-path}
\chi_{\lambda/\mu,a}
=\sum_{\bp  \in \fP(A_n; \bu _{\mu},\bv _{\lambda})}(-1)^{\bp} z^{\bp}_a,  
\end{equation}
by \eqref{eq:h-path-A}.
Applying the method of \cite{GV}, we have
\begin{prop}\label{prop:path}
For any skew diagram $\lambda/\mu$, 
\begin{equation}\label{eq:A-path}
\chi_{\lambda/\mu,a}
=\sum_{\bp \in P(A_n; \mu, \lambda)}z^{\bp}_a,
\end{equation}
where $P(A_n; \mu, \lambda)$ is the set of all 
$\bp = (p_1, \dots , p_l) \in \fP(A_n; \bu _{\mu},\bv _{\lambda})$
which do not have any intersecting pair of paths $(p_i,p_j)$. 
\end{prop}

\begin{proof}
Let $P^c(A_n; \mu, \lambda):=
\{ \bp\in \fP(A_n; \bu_{\mu}, \bv_{\lambda})\mid  \bp \not\in P(A_n; \mu, \lambda)\}$.
The idea of \cite{GV} is to consider an involution 
\begin{equation*}
\iota : P^c(A_n; \mu, \lambda) \to P^c(A_n; \mu, \lambda)
\end{equation*}
defined as follows:
For $\bp = (p_1, \dots , p_l)$,
let $(p_i, p_j)$ be the first intersecting pair of paths, i.e., 
$i$ is the minimal number such that 
$p_i$ intersects with another path 
and $j(\ne i)$ is the minimal number such that 
$p_j$ intersects with $p_i$. 
Let $v_0$ be the first intersecting point of $p_i$ and $p_j$. 
If $u_i\overset{p_i}{\to}v_{\pi(i)}$ ($i = 1, \dots , l$),
then $\iota(\bp)=(p'_1, \dots ,p'_l)$
is given by $p'_k:= p_k$ ($k \ne i,j$) and 
$$p'_i:u_i\overset{p_i}{\to} v_0 \overset{p_j}{\to}v_{\pi(j)}, \qquad 
p'_j:u_j \overset{p_j}{\to} v_0 \overset{p_i}{\to} v_{\pi(i)}.$$ 
Then $\iota$ preserves the weights and inverts the signature, 
i.e., $z^{\iota(\bp)}_a = z^{\bp}_a$ and $(-1)^{\iota(\bp)}=-(-1)^{\bp}$.
Therefore,  the contributions of 
all $\bp \in P^c(A_n; \mu, \lambda)$ to the right hand side of \eqref{eq:det-path}
are canceled with each other.
The signature of any $\bp \in P(A_n; \mu, \lambda)$
is $(-1)^{\bp}= 1$, and we obtain the proposition.
\end{proof}
\subsection{Tableaux description}
\begin{defn}
A tableau $T$ with entries $T(i,j) \in I$ is called an
{\it $A_n$-tableau\/} if it satisfies the following conditions:
\newline
\begin{tabular}{lll}
$(\bH )$ & \text{ horizontal rule } & $T(i,j)\le T(i,j+1)$.\\
$(\bV )$ & \text{ vertical rule }   & $T(i,j) < T(i+1,j)$.
\end{tabular}
\newline
Namely,
an $A_n$-tableau is nothing but a semistandard tableau.
We write the set of all the $A_n$-tableaux of shape $\lambda/\mu$  by 
$\Tab(A_n,\lambda/\mu)$. 
\end{defn}

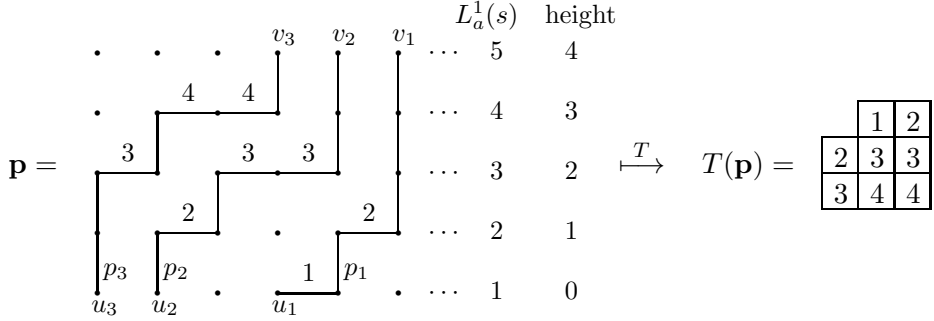
\begin{figure}
{\setlength{\unitlength}{0.8mm}
$\bp =$
\begin{picture}(90,30)(-25,0)
\multiput(-20,20)(10,0){6}{\circle*{1}}
\multiput(-20,10)(10,0){6}{\circle*{1}}
\multiput(-20,0)(10,0){6}{\circle*{1}}
\multiput(-20,-10)(10,0){6}{\circle*{1}}
\multiput(-20,-20)(10,0){6}{\circle*{1}}
{\small 
\put(35,-21){$\cdots \quad 1 \qquad \ 0$}
\put(35,-11){$\cdots \quad 2 \qquad \ 1$}
\put(35,-1){$\cdots \quad 3 \qquad \ 2$}
\put(35,9){$\cdots \quad 4 \qquad \ 3$}
\put(35,19){$\cdots \quad 5 \qquad \ 4$}
\put(35,25){\quad $L^1_a(s)$ \ \ height}
\put(9,-23){$u_1$}
\put(-11,-23){$u_2$}
\put(-21,-23){$u_3$}
\put(29,22){$v_1$}
\put(19,22){$v_2$}
\put(9,22){$v_3$}
\put(14,-18){1}
\put(24,-8){2}
\put(-6,12){4}
\put(4,12){4}
\put(-16,2){3}
\put(-6,-8){2}
\put(4,2){3}
\put(14,2){3}
\put(21,-17){$p_1$}
\put(-9,-17){$p_2$}
\put(-19,-17){$p_3$}
}
\put(10,-20){\line(1,0){10}}
\put(20,-20){\line(0,1){10}}
\put(20,-10){\line(1,0){10}}
\put(30,-10){\line(0,1){10}}
\put(30,0){\line(0,1){10}}
\put(30,10){\line(0,1){10}}
\put(-10,-20){\line(0,1){10}}
\put(-10,-10){\line(1,0){10}}
\put(-10,0){\line(0,1){10}}
\put(-10,10){\line(1,0){10}}
\put(0,10){\line(1,0){10}}
\put(10,10){\line(0,1){10}}
\put(-20,-20){\line(0,1){10}}
\put(-20,-10){\line(0,1){10}}
\put(-20,0){\line(1,0){10}}
\put(0,-10){\line(0,1){10}}
\put(0,0){\line(1,0){10}}
\put(10,0){\line(1,0){10}}
\put(20,0){\line(0,1){10}}
\put(20,10){\line(0,1){10}}
\end{picture}
$\overset{T}{\longmapsto} \quad 
T(\bp )=\  
\text{
\begin{picture}(18,12)(0,0)
\put(0,-6){\line(1,0){18}}
\put(0,0){\line(1,0){18}}
\put(0,6){\line(1,0){18}}
\put(6,12){\line(1,0){12}}
\put(0,-6){\line(0,1){12}}
\put(6,-6){\line(0,1){18}}
\put(12,-6){\line(0,1){18}}
\put(18,-6){\line(0,1){18}}
\put(2,-5){$3$}
\put(8,-5){$4$}
\put(14,-5){$4$}
\put(2,1){$2$}
\put(8,1){$3$}
\put(14,1){$3$}
\put(8,7){$1$}
\put(14,7){$2$}
\end{picture}
}
$
}
\vspace{1.7cm}
\caption{An example of $\bp$ and the tableau $T(\bp)$ for 
$(\lambda, \mu)=((3^3),(1))$.}\label{fig:p-T-correspondence}
\end{figure}

For any $\bp=(p_1, \dots , p_l) \in P(A_n; \mu, \lambda)$, 
we associate a tableau $T(\bp )$ of shape $\lambda/\mu$  
such that the $i$\,th row of $T(\bp )$ 
is given by $\{ L^1_a(s)\mid  s \in E(p_i) \}$ 
listed in the increasing order.
See Figure \ref{fig:p-T-correspondence} for an example.
Clearly, $T(\bp )$ satisfies the horizontal rule 
because of the $h$-labeling rule of $\bp$, and
$T(\bp )$ satisfies the vertical rule 
since $\bp \in P(A_n; \mu, \lambda)$ does not have any intersecting pair of paths.
Therefore, we obtain a map
$$T:P(A_n; \mu, \lambda)\ni \bp \mapsto T(\bp) \in \Tab (A_n, \lambda/\mu)$$
for any skew diagram $\lambda/\mu$. 
In fact, 
\begin{prop}\label{prop:tableau}
The map $T$ is a weight-preserving bijection.
\end{prop}
By Proposition \ref{prop:path} and \ref{prop:tableau}, we reproduce 
the result of \cite{BR}.
\begin{thm}[\cite{BR}]
If $\lambda/\mu$ is a skew diagram, then
$$\chig _{\lambda/\mu, a}
= \sum_{T\in \Tab(A_n,\lambda/\mu)}z^T_a.$$
\end{thm}

\section{Tableaux description of type $B_n$}\label{sec:B}
In this section, we consider the case that $\fg$ is of type $B_n$.
The tableaux description of $\chi_{\lambda/\mu,a}$ \eqref{eq:det} 
is given by \cite{KOS}. We reproduce it
using the path method of \cite{GV}.
During this section, $I$ is of type $B_n$ in \eqref{entries}.

\subsection{Paths description}
\begin{figure}[t]
\begin{picture}(125,108)(-50,-40)
\setlength{\unitlength}{0.8mm}
\multiput(-20,20)(10,0){6}{\circle*{1}}
\multiput(-20,10)(10,0){6}{\circle*{1}}
\multiput(-20,0)(10,0){6}{\circle*{1}}
\multiput(-20,-10)(10,0){6}{\circle*{1}}
\multiput(-20,-20)(10,0){6}{\circle*{1}}
\put(10,-20){\line(1,0){10}}
\put(20,-20){\line(0,1){10}}
\put(20,-10){\line(1,0){10}}
\put(30,-10){\line(0,1){10}}
\put(30,0){\line(0,1){10}}
\put(30,10){\line(0,1){10}}
\put(-10,-20){\line(0,1){10}}
\put(-10,-10){\line(1,0){10}}
\put(-10,0){\line(0,1){10}}
\put(-10,10){\line(1,0){10}}
\put(0,10){\line(1,0){10}}
\put(10,10){\line(0,1){10}}
\put(-20,-20){\line(0,1){10}}
\put(-20,-10){\line(0,1){10}}
\put(-20,0){\line(1,0){10}}
\put(0,-10){\line(0,1){10}}
\put(0,0){\line(1,0){10}}
\put(10,0){\line(1,0){10}}
\put(20,0){\line(0,1){10}}
\put(20,10){\line(0,1){10}}
\put(15,-16){$p_1$}
\put(-15,-16){$p_2$}
\put(-25,-16){$p_3$}
{\small 
\put(37,-21){$\cdots -2$}
\put(37,-11){$\cdots -1$}
\put(37,-1){$\cdots \quad 0$}
\put(37,9){$\cdots \quad  1$}
\put(37,19){$\cdots \quad  2$}
\put(37,25){\quad \ height}
}
\end{picture}
\vspace{0.5cm}
\caption{An example of $h$-paths of type $B_n$ and $C_n$.}\label{fig:B-path-exmp}
\end{figure}
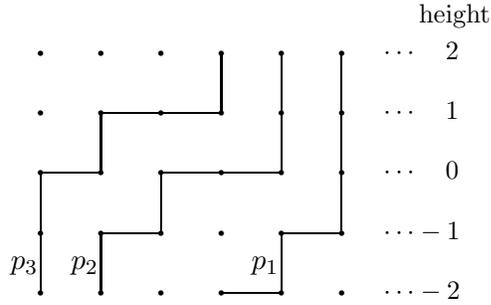

In view of the definition of the generating function of $H_a(z,X)$ in
\eqref{eq:generating-H}, we define an $h$-path and its $h$-labeling as follows:
\begin{defn}\label{def:B-path}
Consider the lattice $\bZ \times \bZ$.
An {\it $h$-path of type $B_n$} is a path $u\overset{p}{\to}v$
such that the initial point $u$ is at height $-n$ and the final point 
$v$ is at height $n$,
and an eastward step at height 0 occurs at most once.
We write the set of all the $h$-paths of type $B_n$ by $P(B_n)$. 
\end{defn}
For example, $p_1$ and $p_3$ in Figure \ref{fig:B-path-exmp}
are the $h$-paths of type $B_n$, but $p_2$ is not.

The {\it $h$-labeling of type $B_n$} associated to $a\in \bC$ for any 
$p\in P(B_n)$ is the pair of maps $L_a=(L_a^1, L_a^2)$,
\begin{align*}
L^1_a :E(p)\to I, \qquad 
L^2_a :E(p)\to \{a+4k\mid k\in\bZ\},
\end{align*}
defined as follows: If $s$ starts at $(x,y)$ then
\begin{equation*}
L^1_a(s)=
\begin{cases}
n+1+y, & \text{\qquad if $y<0$},\\
0, & \text{\qquad if $y=0$},\\
\overline{n+1-y}, & \text{\qquad if $y>0$},
\end{cases}
\end{equation*}
and $L^2_a(s)=a+4x$.
Then, we define $z^p_a$ as in \eqref{eq:z-path}.  
By \eqref{eq:generating-H}, we have
\begin{equation}\label{eq:h-path-B}
h_{r, a+4k+4r-4}(z)=\sum_{p}z^p_a,
\end{equation}
where the sum runs over all $p \in P(B_n)$
such that $(k,-n)\overset{p}{\to}(k+r,n)$.

For any $l$-tuples of initial and final points
$\bu =(u_1, u_2, \dots ,u_l)$, $\bv =(v_1, v_2, \dots , v_l)$, set 
$$\fP(B_n;\bu ,\bv ):=\{\bp=(p_1, \dots, p_l)  \in \fP(\bu ,\bv )\mid 
p_i \in P(B_n) \}.$$

\begin{figure}
{\setlength{\unitlength}{0.8mm}
$\bp = $
\begin{picture}(80,30)(-25,0)
\multiput(-20,20)(10,0){5}{\circle*{1}}
\multiput(-20,10)(10,0){5}{\circle*{1}}
\multiput(-20,0)(10,0){5}{\circle*{1}}
\multiput(-20,-10)(10,0){5}{\circle*{1}}
\multiput(-20,-20)(10,0){5}{\circle*{1}}
\put(-15,-5){$p_2$}
\put(-20,-20){\line(1,0){10}}
\put(-10,-20){\line(0,1){20.3}}
\put(-10,0.2){\line(1,0){9}}
\put(-0.2,1){\line(0,1){19}}
\put(-0.2,20){\line(1,0){10.2}}
\put(2,-15){$p_1$}
{\thicklines
\put(0.2,-20){\line(0,1){10}}
\put(0.2,-10){\line(0,1){9.8}}
\put(0.2,-0.2){\line(1,0){9.8}}
\put(10,-0.2){\line(0,1){10.2}}
\put(10,10){\line(1,0){10}}
\put(20,10){\line(0,1){10}}
}
{\small 
\put(-2,-23){$u_1$}
\put(-22,-23){$u_2$}
\put(19,22){$v_1$}
\put(9,22){$v_2$}
\put(23,-21){$\cdots \quad 1$}
\put(23,-11){$\cdots \quad 2$}
\put(23,-1){$\cdots \quad 0$}
\put(23,9){$\cdots \quad \overline{2}$}
\put(23,19){$\cdots \quad \overline{1}$}
\put(47,-21){$-2$}
\put(47,-11){$-1$}
\put(46,-1){$\quad 0$}
\put(46,9){$\quad 1$}
\put(46,19){$\quad 2$}
\put(27,25){\ \ $L^1_a(s)$ \quad height}
\put(-17,-19){1}
\put(-7,1){0}
\put(3,1){0}
\put(13,11){$\overline{2}$}
\put(3,21){$\overline{1}$}
}
\end{picture}
\quad $\overset{T}{\longmapsto} \quad  
T(\bp )= \ 
\text{
\begin{picture}(18,8)(0,4)
\put(0,0){\line(1,0){18}}
\put(0,6){\line(1,0){18}}
\put(6,12){\line(1,0){12}}
\put(0,0){\line(0,1){6}}
\put(6,0){\line(0,1){12}}
\put(12,0){\line(0,1){12}}
\put(18,0){\line(0,1){12}}
\put(2,1){$1$}
\put(8,1){$0$}
\put(14,1){$\overline{1}$}
\put(8,7){$0$}
\put(14,7){$\overline{2}$}
\end{picture}
}
$
}
\vspace{1.5cm}
\caption{An example of $\bp = (p_1, p_2)\in P(B_n;\lambda,\mu)$ 
which is {\it specially\/} intersecting, 
and their $h$-labelings for $n=2$, 
$(\lambda,\mu)=((3^2),(1))$.
If we set $\iota$ as in the $A_n$ case, 
then $\iota(\bp)\not\in P(B_n; \mu, \lambda)$, because 
the paths in $\iota(\bp)$ are not the $h$-paths of type $B_n$.}
\label{fig:B-path}
\end{figure}
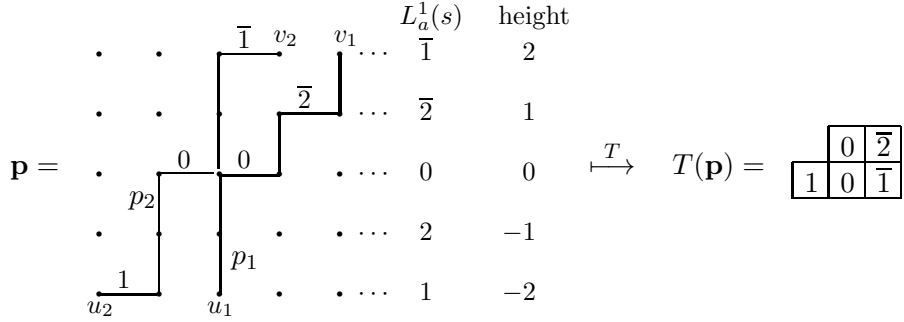

Let $\lambda/\mu$ be a skew diagram and let $l=l(\lambda)$. 
Pick $\bu _{\mu}=(u_1, \dots ,u_l)$ and 
$\bv _{\lambda}=(v_1, \dots ,v_l)$
as $u_i=(\mu_i+1-i, -n)$ and $v_i=(\lambda_i+1-i,n)$.
We define the weight $z^{\bp}$ and its signature $(-1)^{\bp}$
for any $\bp \in \fP(B_n;\bu_{\mu},\bv_{\lambda})$, 
as in the $A_n$ case in \eqref{eq:z-bp}. 
Then, 
the determinant \eqref{eq:det} can be written as 
$$\chi_{\lambda/\mu,a}
= \sum_{\bp \in \fP(B_n;\bu_{\mu},\bv_{\lambda})}(-1)^{\bp}z^{\bp}_a, $$
by \eqref{eq:h-path-B}.
The difference from the $A_n$ case is that 
the involution $\iota$ is not defined on any
$\bp =(p_1, \dots , p_l)\in \fP(B_n;\bu_{\mu},\bv_{\lambda})$ 
that possesses an intersecting pair $(p_i,p_j)$
(see Figure \ref{fig:B-path}).
To define an involution for the $B_n$ case, we give the following definition.

\begin{defn}
An intersecting pair $(p,p')$ of $h$-paths of type $B_n$ is called 
{\it specially intersecting\/} 
(resp.\ {\it ordinarily intersecting\/})
if the intersection of 
$p$ and $p'$ occurs only at height 0
(resp.\ otherwise).
\end{defn}
For example, the pair $(p_1,  p_2)$ given in Figure \ref{fig:B-path}
is specially intersecting.

Applying the method of \cite{GV} as in the $A_n$ case, we have

\begin{prop}\label{prop:B-path}
For any skew diagram $\lambda/\mu$,
\begin{equation}\label{eq:B-path}
\chi_{\lambda/\mu,a}
=\sum_{\bp \in P(B_n; \mu, \lambda)}z^{\bp}_a,
\end{equation}
where $P(B_n; \mu, \lambda)$ is the set of all 
$\bp = (p_1, \dots , p_l)\in \fP(B_n; \bu_{\mu}, \bv_{\lambda})$ 
which do not have any ordinarily intersecting pairs of paths $(p_i, p_j)$. 
\end{prop}

\begin{proof}
Let $P^c(B_n; \mu, \lambda)
:= \{ \bp\in \fP(B_n;\bu_{\mu}, \bv_{\lambda})\mid  \bp \not\in P(B_n; \mu, \lambda)\}$.
Consider an involution 
\begin{equation*}
\iota : P^c(B_n; \mu, \lambda) \to P^c(B_n; \mu, \lambda)
\end{equation*}
defined as follows: For $\bp =(p_1, \dots p_l)$, 
there exists $(p_i, p_j)$ which is ordinarily intersecting.
Let $(p_i, p_j)$ be the first such pair and 
let $v_0$ be the first intersecting point 
whose height is not 0.  
Then set $\iota(\bp)$ as in the proof of Proposition \ref{prop:path}.
Then $\iota$ is weight-preserving and sign-inverting, 
which implies that all $\bp \in P^c(B_n; \mu, \lambda)$ 
will be canceled as in the $A_n$ case.
The signature of any $\bp \in P(B_n; \mu, \lambda)$
is $(-1)^{\bp}=1$, and we obtain the proposition.
\end{proof}

\subsection{Tableaux description}
Define a total ordering in $I$ \eqref{entries} by
$$1 \prec 2 \prec \dots \prec n \prec 0 \prec \overline{n}
\prec \dots  \prec \overline{2} \prec \overline{1}.$$

\begin{defn}[\cite{KOS}]\label{def:B-tab}
A tableau $T$ with  entries $T(i,j)\in I$
is called a {\it $B_n$-tableau\/} if it satisfies the following conditions:
\newline
\begin{tabular}{llcl}
$(\bH )$ & $T(i,j)\preceq T(i,j+1)$ & \text{and} & $(T(i,j),T(i,j+1))\ne (0,0)$. \\
$(\bV )$ & $T(i,j)\prec T(i+1,j)$   & \text{or}  & $(T(i,j),T(i+1,j))=(0,0)$.
\end{tabular}
\newline
We write the set of all the $B_n$-tableaux of shape $\lambda/\mu$ by 
$\Tab(B_n,\lambda/\mu)$. 
\end{defn}

For any $\bp\in P(B_n; \mu, \lambda)$, 
let $T(\bp)$ be the tableau of shape $\lambda/\mu$ defined 
by assigning the $h$-labeling of each path in $\bp$ 
to the corresponding rows, as in the $A_n$ case
(see Figure \ref{fig:B-path}).
Then $T(\bp)$ satisfies the rule $(\bH)$ in Definition \ref{def:B-tab}
because of 
the rule for the $h$-labeling of $\bp$, 
and it satisfies the rule $(\bV)$
since $\bp$ does not have any ordinarily intersecting pairs of paths.
Therefore, we obtain a map
$$T:P(B_n; \mu, \lambda)\ni \bp  \longmapsto T(\bp) \in \Tab(B_n,\lambda/\mu)$$
for any skew diagram $\lambda/\mu$. In fact, 
\begin{prop}
The map $T$ is a weight-preserving bijection.
\end{prop}

Thus, we obtain 
\begin{thm}[\cite{KOS}]
If $\lambda/\mu$ is a skew diagram 
, then
$$\chig _{\lambda/\mu, a}
= \sum_{T\in \Tab(B_n,\lambda/\mu)}z^T_a.$$
\end{thm}

\section{Tableaux description of type $C_n$}\label{sec:C}
In this section, we consider the case that $\fg$ is of type $C_n$.
We determine the tableaux description by the horizontal, vertical and 
``extra'' rules for skew diagrams of at most three rows 
and of at most two columns. 
The one-row and one-column cases are already given by \cite{KS}. 

\subsection{Paths description}
In view of the definition of the generating function of $H_a(z,X)$ in 
\eqref{eq:generating-H}, we define an $h$-path and its $h$-labeling as follows:
\begin{defn}
Consider the lattice $\bZ \times \bZ$.
An {\it $h$-path of type $C_n$} is a path $u\overset{p}{\to}v$
such that the initial point $u$ is at height $-n$ and the final point 
$v$ is at height $n$,
and the number of the eastward steps 
at height 0 is even. 
We write the set of all the $h$-paths of type $C_n$ by $P(C_n)$.
\end{defn}
For example, $p_1$ and $p_2$ in Figure \ref{fig:B-path-exmp}
are the $h$-paths of type $C_n$, but $p_3$ is not.

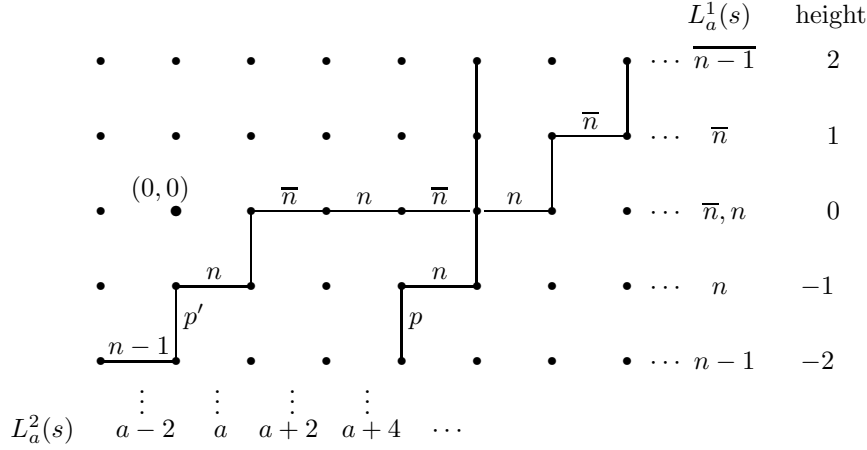
\begin{figure}
{\setlength{\unitlength}{1.0mm}
\begin{picture}(90,30)(-20,0)
\multiput(-20,20)(10,0){8}{\circle*{1}}
\multiput(-20,10)(10,0){8}{\circle*{1}}
\multiput(-20,0)(10,0){8}{\circle*{1}}
\multiput(-20,-10)(10,0){8}{\circle*{1}}
\multiput(-20,-20)(10,0){8}{\circle*{1}}
\put(-20,-20){\line(1,0){10}}
\put(-10,-20){\line(0,1){10}}
\put(-10,-10){\line(1,0){10}}
\put(0,-10){\line(0,1){10}}
\put(0,0){\line(1,0){10}}
\put(10,0){\line(1,0){10}}
\put(20,0){\line(1,0){9}}
\put(31,0){\line(1,0){9}}
\put(40,0){\line(0,1){10}}
\put(40,10){\line(1,0){10}}
\put(50,10){\line(0,1){10}}
{\thicklines
\put(20,-20){\line(0,1){10}}
\put(20,-10){\line(1,0){10}}
\put(30,-10){\line(0,1){10}}
\put(30,0){\line(0,1){10}}
\put(30,10){\line(0,1){10}}
}
{\small
\put(-16,2){$(0,0)$}
\put(-9,-15){$p'$}
\put(21,-15){$p$}
\put(-10,0){\circle*{1.5}}
%
\put(24,-9){$n$}
\put(-19,-19){$n-1$}
\put(-6,-9){$n$}
\put(4,1){$\overline{n}$}
\put(14,1){$n$}
\put(24,1){$\overline{n}$}
\put(34,1){$n$}
\put(44,11){$\overline{n}$}
\put(53,-21){$\cdots \  n-1$}
\put(53,-11){$\cdots \quad n$}
\put(53,-1){$\cdots \ \  \overline{n}, n$}
\put(53,9){$\cdots \quad \overline{n}$}
\put(53,19){$\cdots \ \overline{n-1}$}
\put(73,-21){$-2$}
\put(73,-11){$-1$}
\put(73,-1){$\quad 0$}
\put(73,9){$\quad  1$}
\put(73,19){$\quad  2$}
\put(58,25){$L^1_a(s)$ \quad \ height}
\put(-15,-27){$\vdots$}
\put(-5,-27){$\vdots$}
\put(5,-27){$\vdots$}
\put(15,-27){$\vdots$}
%
\put(-32,-30){$L^2_a(s)$}
\put(-18,-30){$a-2$}
\put(-5,-30){$a$}
\put(1,-30){$a+2$}
\put(12,-30){$a+4$}
\put(24,-30){$\cdots$}
}
\end{picture}
}
\vspace{3cm}
\caption{An example of $h$-paths of type $C_n$ 
and their $h$-labelings.}\label{fig:c-label}
\end{figure}

For a path $p =(s_1, s_2, \dots ) \in P(C_n)$,
let $E_0(p)=(s_{j},s_{j+1}, \dots )$
be the sequence of all the eastward steps at height 0 in $p$.
Let $E^1_0(p)$ and $E^2_0(p)$ 
be the subsequence of $E_0(p)$ defined by
$E^1_0(p)=(s_{j}, s_{j+2}, s_{j+4}, \dots )$ and 
$E^2_0(p)=(s_{j+1}, s_{j+3}, s_{j+5}, \dots )$. 
The {\it $h$-labeling of type $C_n$} associated to $a\in \bC$ 
for any $p\in P(C_n)$ is a pair of maps $L_a=(L_a^1,L_a^2)$,
\begin{equation}\label{eq:h-label-C}
L^1_a :E(p)\to I, \qquad L^2_a :E(p)\to \{a+2k\mid k\in\bZ\},
\end{equation}
defined as follows: If $s$ starts at $(x,y)$, then
\begin{equation*}
L^1_a(s)=
\begin{cases}
n+1+y,            & \text{\qquad if $y<0$},\\
\overline{n+1-y}, & \text{\qquad if $y>0$},\\
\overline{n},     & \text{\qquad if $s \in E^1_0(p)$},\\
n,                & \text{\qquad if $s \in E^2_0(p)$},
\end{cases}
\end{equation*}
and $L^2_a(s)=a+2x$.
See Figure \ref{fig:c-label} for an example. 

Define $z^p_a$ as in \eqref{eq:z-path}. 
By \eqref{eq:generating-H}, we have
\begin{equation}\label{eq:h-path-C}
h_{r, a+2k+2r-2}(z)=\sum_{p}z^p_a,
\end{equation}
where the sum runs over all $p\in P(C_n)$
such that $(k,-n)\overset{p}{\to}(k+r,n)$.

For any $l$-tuples of initial and final points
$\bu =(u_1, u_2, \dots ,u_l)$, $\bv =(v_1, v_2, \dots , v_l)$, set 
$$\fP(C_n;\bu ,\bv ):=\{\bp=(p_1, \dots, p_l)  \in \fP(\bu ,\bv )\mid 
\text{$p_i \in P(C_n)$} \}.$$

Let $\lambda/\mu$ be a skew diagram and let 
$l=l(\lambda)$. 
Pick $\bu _{\mu}=(u_1, \dots ,u_l)$ and
$\bv _{\lambda}=(v_1, \dots ,v_l)$
as $u_i=(\mu_i+1-i, -n)$ and $v_i=(\lambda_i+1-i,n)$.
We define the weight $z^{\bp}$ 
and the signature $(-1)^{\bp}$ for any $\bp \in \fP(C_n;\bu ,\bv )$
by the $h$-labeling of type $C_n$ as in \eqref{eq:z-bp}. Then, 
the determinant \eqref{eq:det} can be written as
$$\chi_{\lambda/\mu,a}
=\sum_{\bp  \in \fP(C_n;\bu _{\mu},\bv _{\lambda})}(-1)^{\bp} z^{\bp}_a, $$
by \eqref{eq:h-path-C}.

As in the $B_n$ case, the involution $\iota $ is not defined 
on any $\bp =(p_1, \dots, p_l) \in \fP(C_n;\bu _{\mu},\bv _{\lambda})$
which possesses an intersecting pair of paths $(p_i,p_j)$.
To define the involution for the $C_n$ case, we give the 
definition of the specially (resp.\ ordinarily) intersecting pair of paths.

Consider two paths $p, p'$ which is intersecting at height 0.
Let $(x, 0)$ (resp.\ $(x', 0)$)
be the leftmost point on $p$ (resp.\ $p'$) at height 0.
Then set $[p,p']:=|x-x'|$.

\begin{defn}
An intersecting pair $(p,p')$ of $h$-paths of type $C_n$ is called 
{\it specially intersecting\/} (resp.\ {\it ordinarily intersecting\/}) if 
the intersection of $p$ and $p'$ occurs only at height 0 
and $[p,p']$ is odd (resp.\ otherwise).
\end{defn}
For example, $[p,p']=3$
for $(p,p')$ in Figure \ref{fig:c-label}, and therefore, 
it is specially intersecting.

Applying the method of \cite{GV}, we have
\begin{prop}\label{prop:alternative-path}
For any skew diagram $\lambda/\mu$,
\begin{equation}\label{eq:pre-transposed}
\chi_{\lambda/\mu,a}
=\sum_{\bp \in P(C_n; \mu, \lambda)}(-1)^{\bp}z^{\bp}_a,
\end{equation}
where $P(C_n; \mu, \lambda)$ is the set of all 
$\bp =(p_1, \dots , p_l)\in \fP(C_n; \bu_{\mu}, \bv_{\lambda})$ 
which do not have any ordinarily intersecting 
pair of paths $(p_i, p_j)$.
\end{prop}

\begin{proof}
Let $P^c(C_n; \mu, \lambda):=
\{\bp \in \fP(C_n;\bu _{\mu},\bv _{\lambda})\mid  \bp \not\in P(C_n; \mu, \lambda)\}$.
Consider a weight-preserving involution
$$\iota :P^c(C_n; \mu, \lambda) \to P^c(C_n; \mu, \lambda)$$
defined as follows:
For $\bp =(p_1, \dots ,p_l)$,
let $(p_i, p_j)$ be the first ordinarily intersecting pair of paths
and let $v_0$ be the first intersecting point.
Set $\iota (\bp )$ as in the $A_n$ case (Proposition \ref{prop:path}).
Then $\iota $ is weight-preserving and sign-inverting, 
as in the $A_n$ and $B_n$ cases.
\end{proof}

Consider any two $h$-paths of type $C_n$, 
$(x_0,y_0)\overset{p}{\to}(x_1,y_1)$, $(x'_0,y'_0)\overset{p'}{\to}(x'_1,y'_1)$, 
which are not ordinarily intersecting.
We say that $(p, p')$ is {\it transposed\/} 
if $(x_0-x'_0)(x_1-x'_1)<0$. For example, 
the pair $(p,p')$ in Figure \ref{fig:c-label} is transposed.

Let $P_k(C_n; \mu, \lambda)$ be the set of all
$\bp \in P(C_n; \mu, \lambda)$ which possess exactly $k$ transposed pairs
of paths. Then we have $(-1)^{\bp} = (-1)^k$.
Note that if $\bp = (p_1, \dots p_l) \in P(C_n; \mu, \lambda)$, 
then each triplet $(p_i, p_j, p_k)$ is not intersecting simultaneously
at one point. Therefore, 
$P(C_n; \mu, \lambda)=\sum_{k=0}^{l-1}P_k(C_n; \mu, \lambda)$ and 
the sum \eqref{eq:pre-transposed} is rewritten as
\begin{equation}\label{eq:transposed-paths}
\chi_{\lambda/\mu,a}
=\sum_{k=0}^{l-1}(-1)^k
\sum_{\bp \in P_k(C_n; \mu, \lambda)}z^{\bp}_a.
\end{equation}

The right hand side of \eqref{eq:transposed-paths} is not as simple as 
that of \eqref{eq:A-path} for $A_n$ and that of \eqref{eq:B-path} for $B_n$. 
This is the main reason why the description of $C_n$ becomes 
more complicated than that of $A_n$ and $B_n$. 
(The $D_n$ case which is not dealt with in this paper is similar 
to the $C_n$ case.)

\subsection{Tableaux description}\label{subsec:C-Tab}
To formulate the tableaux description of \eqref{eq:transposed-paths} 
we introduce a certain set of tableaux (called $HV$-tableaux) 
$\widetilde{\Tab}(C_n, \lambda/\mu)$ and the corresponding set of paths 
$\tilde{P}(C_n; \mu, \lambda)$.

Define a total ordering in $I$ \eqref{entries} by
$$1 \prec 2 \prec \dots \prec n \prec \overline{n}
\prec \dots  \prec \overline{2} \prec \overline{1}.$$

\begin{figure}
{\setlength{\unitlength}{0.8mm}
$\bp = $
\begin{picture}(80,30)(-25,0)
\multiput(-20,20)(10,0){5}{\circle*{1}}
\multiput(-20,10)(10,0){5}{\circle*{1}}
\multiput(-20,0)(10,0){5}{\circle*{1}}
\multiput(-20,-10)(10,0){5}{\circle*{1}}
\multiput(-20,-20)(10,0){5}{\circle*{1}}
\put(-19,-15){$p_3$}
\put(-20,-20){\line(0,1){20.6}}
\put(-20,0.6){\line(1,0){19.7}}
\put(-0.3,0.6){\line(0,1){19.4}}
{\thicklines
\put(-9,-15){$p_2$}
\put(-10,-20){\line(0,1){20}}
\put(-10,0){\line(1,0){20}}
\put(10,0){\line(0,1){20}}
}
\put(1.3,-15){$p_1$}
\put(0.3,-20){\line(0,1){19.4}}
\put(0.3,-0.6){\line(1,0){19.7}}
\put(20,-0.6){\line(0,1){20.6}}
{\small 
\put(-2,-24){$u_1$}
\put(-12,-24){$u_2$}
\put(-22,-24){$u_3$}
\put(19,21){$v_1$}
\put(9,21){$v_2$}
\put(-1,21){$v_3$}
\put(23,-21){$\cdots \ n-1$}
\put(23,-11){$\cdots \quad n$}
\put(23,-1){$\cdots \ \ \overline{n}, n$}
\put(23,9){$\cdots \quad \overline{n}$}
\put(23,19){$\cdots \ \overline{n-1}$}
\put(47,-21){$-2$}
\put(47,-11){$-1$}
\put(46,-1){$\quad 0$}
\put(46,9){$\quad 1$}
\put(46,19){$\quad 2$}
\put(27,25){\ \ $L^1_a(s)$ \quad height}
\put(-17,1){$\overline{n}$}
\put(-7,1){$n$}
\put(-7,-5){$\overline{n}$}
\put(3,1){$n$}
\put(3,-5){$\overline{n}$}
\put(13,-5){$n$}
}
\end{picture}
\quad $\overset{T}{\longmapsto} \quad  
T(\bp )= \ 
\text{
\begin{picture}(12,9)(0,9)
\put(0,0){\line(1,0){12}}
\put(0,6){\line(1,0){12}}
\put(0,12){\line(1,0){12}}
\put(0,18){\line(1,0){12}}
\put(0,0){\line(0,1){18}}
\put(6,0){\line(0,1){18}}
\put(12,0){\line(0,1){18}}
\put(2,1){$\overline{n}$}
\put(2,7){$\overline{n}$}
\put(2,13){$\overline{n}$}
\put(8,1){$n$}
\put(8,7){$n$}
\put(8,13){$n$}
\end{picture}
}
$
}
\vspace{1.5cm}
\caption{An example of $\bp \in \tilde{P}(C_n; \mu, \lambda)$
for $(\lambda,\mu)=((2^3), \phi)$
and its $h$-labeling. 
The pair $(p_1,p_3)$ is ordinarily intersecting, 
and therefore, $\bp \not\in P(C_n; \mu, \lambda)$.}
\label{fig:C-path}
\end{figure}
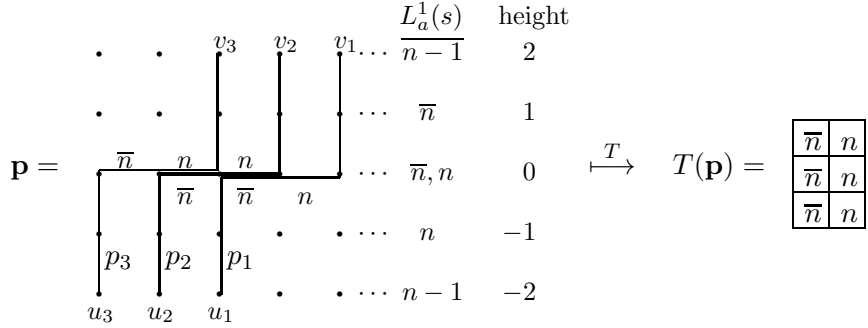

\begin{defn}\label{def:HV-tab}
A tableau $T$ (of shape $\lambda/\mu$) with entries $T(i,j)\in I$
is called an {\it $HV$-tableau\/} if it satisfies the following conditions: 
\newline
\begin{tabular}{ll}
$(\bH )$ & Each $(i,j)\in \lambda/\mu$ satisfies both of the following conditions:
\end{tabular}
\begin{itemize}
\item $T(i,j)\preceq T(i,j+1)$ or $(T(i,j), T(i,j+1))=(\overline{n}, n)$.
\item $(T(i,j-1), T(i,j), T(i,j+1))\ne
    (\overline{n}, \overline{n}, n), (\overline{n}, n, n)$.
\end{itemize}
\begin{tabular}{ll}
$(\bV )$ & Each $(i,j)\in \lambda/\mu$ satisfies 
at least one of the following conditions:
\end{tabular}
\begin{itemize}
\item $T(i,j)\prec T(i+1,j)$.
\item $T(i,j)=T(i+1,j)=n$,  
$(i+1,j-1)\in \lambda/\mu$ and $T(i+1,j-1)=\overline{n}$.
\item $T(i,j)=T(i+1,j)=\overline{n}$, 
$(i,j+1)\in \lambda/\mu$ and $T(i,j+1)=n$.
\end{itemize}
We write the set of the all the $HV$-tableaux of shape $\lambda/\mu$ 
by  $\widetilde{\Tab}(C_n,\lambda/\mu)$. 
\end{defn}
Let $\tilde{P}(C_n; \mu, \lambda)$ be the set of all 
$\bp =(p_1, \dots ,p_l)\in  \fP(C_n; \bu_{\mu}, \bv_{\lambda})$
which do not have any {\it adjacent\/} pair $(p_i,p_{i+1})$ which is either
ordinarily intersecting or transposed. 
We remark that
$P_k(C_n; \mu, \lambda)\cap \tilde{P}(C_n; \mu, \lambda) = \phi$ 
($k \ge 1$) and 
\begin{gather*}
P_0(C_n; \mu, \lambda)=\tilde{P}(C_n; \mu, \lambda) 
\quad \text{if $l(\lambda)\le 2$},  \\
P_0(C_n; \mu, \lambda)\subsetneqq \tilde{P}(C_n; \mu, \lambda) \quad 
\text{if $l(\lambda) \ge 3$}.  
\end{gather*}
For any $\bp \in \tilde{P}(C_n; \mu, \lambda)$, 
let $T(\bp)$ be the tableau of shape $\lambda/\mu$ 
defined by assigning the $h$-labeling of each path in $\bp$ 
to the corresponding row, as in the $A_n$ and $B_n$ cases
(see Figure \ref{fig:C-path}). 
Then $T(\bp)$ is an $HV$-tableau; 
it satisfies the rule $(\bH)$ in Definition \ref{def:HV-tab}
because of
the rule for the $h$-labeling of $\bp$,
and it satisfies the rule $(\bV)$
since $\bp$ does not have any adjacent pairs of paths which is 
either ordinarily intersecting or transposed.
Therefore, we obtain a map
$$T:\tilde{P}(C_n; \mu, \lambda) \ni \bp \mapsto T(\bp)
\in \widetilde{\Tab}(C_n, \lambda/\mu)$$
for any skew diagram $\lambda/\mu$ as similar as the previous cases. 
Moreover, 
\begin{prop}\label{prop:bijective-T}
The map $T$ is a weight-preserving bijection.
\end{prop}

We expect that the alternative sum \eqref{eq:transposed-paths} 
can be translated into the following positive sum by tableaux, 
\begin{equation}\label{eq:C-tab}
\chig_{\lambda/\mu,a}
= \sum_{T\in \Tab(C_n, \lambda/\mu)}z^T_a,
\end{equation}
where $\Tab (C_n, \lambda/\mu)$ is a certain subset of
$\widetilde{\Tab}(C_n, \lambda/\mu)$. 
Thus, the tableaux in $\Tab(C_n, \lambda/\mu)$ are described by the horizontal 
rules $(\bH)$ and the vertical rules $(\bV)$ in Definition \ref{def:HV-tab}, 
{\it and\/} the {\it extra\/} rules which select them out of 
$\widetilde{\Tab}(C_n, \lambda/\mu)$. 

In the following subsections, we show how the tableaux description 
\eqref{eq:C-tab} is naturally obtained from 
\eqref{eq:transposed-paths} for 
the skew diagrams $\lambda/\mu$ of at most three rows 
and of at most two columns. 

Roughly speaking, the idea is as follows 
(see \eqref{eq:C-maps} and \eqref{eq:lambda-4}): 
We introduce the weight-preserving 
maps $f_k$ which ``resolve'' the intersection of a transposed pair 
of $\bp \in P_k(C_n; \mu, \lambda)$ in \eqref{eq:transposed-paths},  
and show that the contributions for \eqref{eq:transposed-paths} from 
$P_k(C_n; \mu,\lambda)$ ($k \ge 1$) almost cancel with each other. 
Then, the remaining positive contributions fill the difference 
$\tilde{P}(C_n; \mu, \lambda) \setminus P_0(C_n; \mu, \lambda)$, 
while the remaining negative contributions 
turn into the extra rules. 
We remark that the relation \eqref{gen} plays a crucial role in the 
weight-preserving property of the maps $f_k$.

\subsection{Skew diagrams of at most three rows}\label{sec:bijective-T3}
In this subsection, we consider the tableaux description 
for skew diagrams of at most three rows. 

{\it The case of one-row\/}. 
Let $\lambda/\mu$ be a one-row diagram, i.e., $l(\lambda)=1$. 
Then there does not exist any $\bp \in P(C_n; \mu, \lambda)$ which possesses 
a transposed pair of paths, and therefore, 
$P(C_n; \mu, \lambda)=P_0(C_n; \mu, \lambda)=\tilde{P}(C_n; \mu, \lambda)$. 
Thus, $\Tab(C_n, \lambda/\mu)$ 
in the equality \eqref{eq:C-tab} 
is exactly the set $\widetilde{\Tab}(C_n, \lambda/\mu)$.  

{\it The case of two-row\/}. 
Let $\lambda/\mu$ be a skew diagram of two rows, i.e., $l(\lambda)=2$. 
Let $\Tab(C_n,\lambda/\mu)$ be the set of all the $HV$-tableaux $T$ with the following
extra condition: 

$(\bEtwoR )$ \quad If $T$ contains
a subtableau 
(excluding $a$ and $b$)
\begin{equation}\label{eq:arr-2}
\overbrace{%
\begin{picture}(65,13)(0,13)
\put(0,0){\line(1,0){65}}
\put(0,13){\line(1,0){65}}
\put(0,26){\line(1,0){65}}
\put(0,0){\line(0,1){26}}
\put(13,0){\line(0,1){26}}
\put(26,0){\line(0,1){26}}
\put(52,0){\line(0,1){26}}
\put(65,0){\line(0,1){26}}
\text{
\put(3,3){$\overline{n}$}
\put(3,16){$n$}
\put(16,3){$\overline{n}$}
\put(16,16){$n$}
\put(55,3){$\overline{n}$}
\put(55,16){$n$}
\put(33,3){$\cdots$}
\put(33,16){$\cdots$}
\put(-8,3){$b$}
\put(68,16){$a$}
}
\end{picture}%
}^{k}
\vspace{0.4cm}
\end{equation}
where $k$ is an odd number, 
then at least one of the following conditions holds:
\begin{enumerate}
\item \label{item:con-1} 
Let $(i_1,j_1)$ be the position of the top-right corner 
of the subtableau \eqref{eq:arr-2}. Then 
$(i_1,j_1+1) \in \lambda/\mu$ and $a:=T(i_1, j_1+1)=n$.
\item \label{item:con-2} 
Let $(i_2,j_2)$ be the position of the bottom-left corner 
of the subtableau \eqref{eq:arr-2}. Then 
$(i_2, j_2-1)\in \lambda/\mu$ and $b:=T(i_2, j_2-1)=\overline{n}$.
\end{enumerate}
Then  
\begin{thm}\label{thm:2-row}
For any skew diagram $\lambda /\mu$ with $l(\lambda)=2$, 
the equality \eqref{eq:C-tab} holds.
\end{thm}

\begin{proof}
For $\lambda /\mu$ is a skew diagram of $l(\lambda)=2$, 
there does not exist any $\bp \in P(C_n; \mu, \lambda)$ that possesses 
more than one transposed pair of paths, we have 
$P(C_n; \mu, \lambda) = P_0(C_n; \mu, \lambda) \sqcup P_1(C_n; \mu, \lambda)$.
We also have $\tilde{P}(C_n; \mu, \lambda)=P_0(C_n; \mu, \lambda)$.
Define $f_1:P_1(C_n; \mu, \lambda)\to P_0(C_n; \mu, \lambda)$ 
by $f_1=r_0$, 
where $r_0$ is a weight-preserving injection
defined as in Appendix \ref{sec:maps}. See also Figure \ref{fig:r}. 
Roughly speaking, $r_0$ is a map which resolves the intersection of
specially intersecting paths. From 
\eqref{eq:transposed-paths} and Proposition \ref{prop:bijective-T}, 
we have 
\begin{equation*}
\chi_{\lambda/\mu,a}
=\sum_{T\in \widetilde{\Tab}(C_n, \lambda/\mu)}z^{T}_a
- \sum_{\bp \in \Image \, f_1}z^{T(\bp)}_a.
\end{equation*}
The set $\{ T(\bp)\mid  \bp \in \Image \, f_1\}$ consists of all the tableaux in
$\widetilde{\Tab}(C_n, \lambda/\mu)$ prohibited by 
the extra rule $(\bEtwoR )$.
\end{proof}




{\it The case of three-row\/}. 
Let $\lambda/\mu$ be a skew diagram of three rows, i.e., $l(\lambda)= 3$. 
Let $\Tab(C_n, \lambda/\mu)$ be all the $HV$-tableaux $T$ 
which satisfy $(\bEtwoR )$  
and the following conditions: 
\newline
$(\bEthreeR)$
\begin{enumerate}
\item If $T$ contains a subtableau
(excluding $a$ and $b$)
\begin{equation}\label{eq:arr-3-1}
\overbrace{
\begin{picture}(47,6)(1,20)
{\thicklines
\put(0,0){\line(1,0){49}}
\put(0,13){\line(1,0){49}}
\put(0,26){\line(1,0){49}}
}
\put(0,0){\line(0,1){26}}
\text{
\put(-8,3){$\scriptstyle b$}
\put(3,3){$\scriptstyle \overline{n}$}
\put(3,16){$\scriptstyle n$}
\put(39,3){$\scriptstyle \overline{n}$}
\put(18,3){$\cdots$}
\put(18,16){$\cdots$}
}
\end{picture}
}^{k_1}
\overbrace{
\begin{picture}(49,19)(1,20)
{\thicklines
\put(0,0){\line(1,0){51}}
\put(0,13){\line(1,0){51}}
\put(0,26){\line(1,0){51}}
\put(0,39){\line(1,0){51}}
\put(0,0){\line(0,1){13}}
}
\put(0,26){\line(0,1){13}}
\multiput(0,13)(0,4){4}{\line(0,1){1}}
\text{
\put(1,3){${\scriptstyle \overline{n-1}}$}
\put(1,29){${\scriptstyle n-1}$}
\put(41,16){$\scriptstyle n$}
\put(21,3){$\cdots$}
\put(21,16){$\cdots$}
\put(21,29){$\cdots$}
}
\end{picture}
}^{k_2}
\overbrace{
\begin{picture}(104,19)(1,20)
{\thicklines
\put(0,0){\line(1,0){106}}
\put(0,13){\line(1,0){106}}
\put(0,26){\line(1,0){106}}
\put(0,39){\line(1,0){106}}
\put(0,13){\line(0,1){13}}
}
\multiput(0,0)(0,4){4}{\line(0,1){1}}
\multiput(0,26)(0,4){4}{\line(0,1){1}}
\text{
\put(3,16){$\scriptstyle \overline{n}$}
\put(16,16){$\scriptstyle n$}
\put(29,16){$\scriptstyle \overline{n}$}
\put(42,16){$\scriptstyle n$}
\put(83,16){$\scriptstyle \overline{n}$}
\put(96,16){$\scriptstyle n$}
\put(60,3){$\cdots$}
\put(60,16){$\cdots$}
\put(60,29){$\cdots$}
}
\end{picture}
}^{2k_3}
\overbrace{
\begin{picture}(49,19)(1,20)
{\thicklines
\put(0,0){\line(1,0){51}}
\put(0,13){\line(1,0){51}}
\put(0,26){\line(1,0){51}}
\put(0,39){\line(1,0){51}}
\put(0,13){\line(0,1){13}}
}
\multiput(0,0)(0,4){4}{\line(0,1){1}}
\multiput(0,26)(0,4){4}{\line(0,1){1}}
\text{
\put(3,16){$\scriptstyle \overline{n}$}
\put(34,3){${\scriptstyle \overline{n-1}}$}
\put(34,29){${\scriptstyle n-1}$}
\put(15,3){$\cdots$}
\put(15,16){$\cdots$}
\put(15,29){$\cdots$}
}
\end{picture}
}^{k_4}
\overbrace{
\begin{picture}(47,19)(1,7)
{\thicklines 
\put(0,0){\line(1,0){49}}
\put(0,13){\line(1,0){49}}
\put(0,26){\line(1,0){49}}
\put(0,13){\line(0,1){13}}
}
\put(0,-13){\line(0,1){13}}
\multiput(0,0)(0,4){4}{\line(0,1){1}}
\put(49,0){\line(0,1){26}}
\text{
\put(3,16){$\scriptstyle n$}
\put(39,3){$\scriptstyle \overline{n}$}
\put(39,16){$\scriptstyle n$}
\put(18,3){$\cdots$}
\put(18,16){$\cdots$}
\put(52,16){$\scriptstyle a$}
}
\end{picture}
}^{k_5}
\vspace{0.6cm}
\end{equation}
where $k_i \in \bZ_{\ge 0}$ and $k_1+k_2+k_4+k_5$ is an odd number
with $k_2\ne 0$ or $k_4 \ne 0$, 
then at least one of the following conditions holds.
\begin{enumerate}
\item \label{item:length-three-a}
Let $(i_1,j_1)$ be the position of the top-right corner of 
the subtableau \eqref{eq:arr-3-1}. Then 
$(i_1, j_1+1)\in \lambda/\mu$ and 
$a:=T(i_1,j_1+1)\prec T(i_1+1, j_1)$.
\item \label{item:length-three-b}
Let $(i_2,j_2)$ be the position of the bottom-left corner of the 
subtableau \eqref{eq:arr-3-1}. Then 
$(i_2, j_2-1)\in \lambda/\mu$ and $b:=T(i_2,j_2-1)\succ T(i_2-1, j_2)$.
\end{enumerate}
\item If $T$ contains the subtableau
(excluding $a$)
\begin{equation}\label{eq:arr-3-2}
\overbrace{
\begin{picture}(29,19)(1,20)
{\thicklines 
\put(0,0){\line(1,0){31}}
\put(0,13){\line(1,0){31}}
\put(0,26){\line(1,0){31}}
\put(13,39){\line(1,0){18}}
\put(13,0){\line(0,1){13}}
}
\put(0,0){\line(0,1){26}}
\multiput(13,13)(0,4){4}{\line(0,1){1}}
\put(13,26){\line(0,1){13}}
\text{
\put(3,3){$\scriptstyle \overline{n}$}
\put(3,16){$\scriptstyle \overline{n}$}
\put(14,3){$\scriptstyle{\overline{n-1}}$}
\put(19,16){$\scriptstyle n$}
\put(14,29){$\scriptstyle n-1$}
}
\end{picture}
}^2
\overbrace{
\begin{picture}(78,19)(1,20)
{\thicklines
\put(0,0){\line(1,0){80}}
\put(0,13){\line(1,0){80}}
\put(0,26){\line(1,0){80}}
\put(0,39){\line(1,0){80}}
}
\multiput(0,0)(0,4){10}{\line(0,1){1}}
\text{
\put(3,16){$\scriptstyle \overline{n}$}
\put(16,16){$\scriptstyle n$}
\put(57,16){$\scriptstyle \overline{n}$}
\put(70,16){$\scriptstyle n$}
\put(34,3){$\cdots$}
\put(34,16){$\cdots$}
\put(34,29){$\cdots$}
}
\end{picture}
}^{2k_3}
\overbrace{
\begin{picture}(52,19)(1,20)
{\thicklines
\put(0,0){\line(1,0){54}}
\put(0,13){\line(1,0){54}}
\put(0,26){\line(1,0){54}}
\put(0,39){\line(1,0){54}}
\put(0,13){\line(0,1){13}}
}
\multiput(0,0)(0,4){4}{\line(0,1){1}}
\multiput(0,26)(0,4){4}{\line(0,1){1}}
\text{
\put(3,16){$\scriptstyle \overline{n}$}
\put(37,3){${\scriptstyle \overline{n-1}}$}
\put(37,29){${\scriptstyle n-1}$}
\put(18,3){$\cdots$}
\put(18,16){$\cdots$}
\put(18,29){$\cdots$}
}
\end{picture}
}^{k_4}
\overbrace{
\begin{picture}(50,19)(1,7)
{\thicklines 
\put(0,0){\line(1,0){52}}
\put(0,13){\line(1,0){52}}
\put(0,26){\line(1,0){52}}
\put(0,13){\line(0,1){13}}
}
\put(0,-13){\line(0,1){13}}
\multiput(0,0)(0,4){4}{\line(0,1){1}}
\put(52,0){\line(0,1){26}}
\text{
\put(3,16){$\scriptstyle n$}
\put(42,3){$\scriptstyle \overline{n}$}
\put(42,16){$\scriptstyle n$}
\put(20,3){$\cdots$}
\put(20,16){$\cdots$}
\put(55,16){$\scriptstyle a$}
}
\end{picture}
}^{k_5}
\vspace{0.6cm}
\end{equation}
where
$k_i \in \bZ_{\ge 0}$ and $k_4+k_5$ is an odd number
with $k_4\ne 0$, then the following holds: 
Let $(i,j)$ be the position of the top-right corner of 
the subtableau \eqref{eq:arr-3-2}. Then 
$(i,j+1)\in \lambda/\mu$ and $a:=T(i,j+1)\prec T(i+1, j)$. 
\item If $T$ contains the subtableau
(excluding $b$)
\begin{equation}\label{eq:arr-3-3}
\overbrace{
\begin{picture}(50,6)(1,20)
{\thicklines
\put(0,0){\line(1,0){52}}
\put(0,13){\line(1,0){52}}
\put(0,26){\line(1,0){52}}
}
\put(0,0){\line(0,1){26}}
\text{
\put(-8,3){$\scriptstyle b$}
\put(3,3){$\scriptstyle \overline{n}$}
\put(3,16){$\scriptstyle n$}
\put(42,3){$\scriptstyle \overline{n}$}
\put(20,3){$\cdots$}
\put(20,16){$\cdots$}
}
\end{picture}
}^{k_1}
\overbrace{
\begin{picture}(52,19)(1,20)
{\thicklines
\put(0,0){\line(1,0){54}}
\put(0,13){\line(1,0){54}}
\put(0,26){\line(1,0){54}}
\put(0,39){\line(1,0){54}}
\put(0,0){\line(0,1){13}}
}
\put(0,26){\line(0,1){13}}
\multiput(0,13)(0,4){4}{\line(0,1){1}}
\text{
\put(1,3){${\scriptstyle \overline{n-1}}$}
\put(1,29){${\scriptstyle n-1}$}
\put(44,16){$\scriptstyle n$}
\put(23,3){$\cdots$}
\put(23,16){$\cdots$}
\put(23,29){$\cdots$}
}
\end{picture}
}^{k_2}
\overbrace{
\begin{picture}(73,19)(1,20)
{\thicklines
\put(0,0){\line(1,0){75}}
\put(0,13){\line(1,0){75}}
\put(0,26){\line(1,0){75}}
\put(0,39){\line(1,0){75}}
\put(0,13){\line(0,1){13}}
}
\multiput(0,0)(0,4){4}{\line(0,1){1}}
\multiput(0,26)(0,4){4}{\line(0,1){1}}
\text{
\put(3,16){$\scriptstyle \overline{n}$}
\put(16,16){$\scriptstyle n$}
\put(52,16){$\scriptstyle \overline{n}$}
\put(65,16){$\scriptstyle n$}
\put(31,3){$\cdots$}
\put(31,16){$\cdots$}
\put(31,29){$\cdots$}
}
\end{picture}
}^{2k_3}
\overbrace{
\begin{picture}(29,19)(1,20)
{\thicklines
\put(0,0){\line(1,0){18}}
\put(0,13){\line(1,0){31}}
\put(0,26){\line(1,0){31}}
\put(0,39){\line(1,0){31}}
\put(18,26){\line(0,1){13}}
}
\multiput(0,0)(0,4){10}{\line(0,1){1}}
\put(18,0){\line(0,1){13}}
\multiput(18,13)(0,4){4}{\line(0,1){1}}
\put(31,13){\line(0,1){26}}
\text{
\put(1,3){$\scriptstyle \overline{n-1}$}
\put(5,16){$\scriptstyle \overline{n}$}
\put(1,29){$\scriptstyle \scriptstyle n-1$}
\put(21,16){$\scriptstyle n$}
\put(21,29){$\scriptstyle n$}
}
\end{picture}
}^2
\vspace{0.6cm}
\end{equation}
where $k_i \in \bZ_{\ge 0}$ and $k_1+k_2$ is an odd number 
with $k_2\ne 0$, then the following holds: 
Let $(i,j)$ be the position of the bottom-left corner of the 
subtableau \eqref{eq:arr-3-3}. Then 
$(i,j-1)\in T$ and $b:=T(i,j-1) \succ T(i-1,j)$.
\end{enumerate}
Then 
\begin{thm}\label{thm:3-row}
For any skew diagram $\lambda/\mu$ of $l(\lambda)=3$, 
the equality \eqref{eq:C-tab} holds.
\end{thm}

\begin{proof}
In this proof, we use some maps which are defined in detail in Appendix \ref{sec:maps}. 
For a summary of this proof, see the maps and their relations in the following
diagram: 
\begin{equation}\label{eq:C-maps}
\text{
\begin{picture}(180,50)(0,50)
\put(0,93){$P_2^{\times}$}
\put(40,93){$\sqcup$} 
\put(85,93){$P_2^{\circ}$} 
\put(130,93){$=$}
\put(150,93){$P_2$}
\put(7,83){\vector(0,-1){63}}
\put(85,83){\vector(-1,-1){20}}
\put(95,83){\vector(1,-1){20}}
\put(55,48){$P^{12}_1$}
\put(85,48){$\sqcup$}
\put(110,48){$P^{23}_1$}
\put(140,48){$=$}
\put(160,48){$P_1$}
\put(65,40){\vector(1,-1){20}}
\put(116,40){\vector(-1,-1){20}}
\put(0,5){$\Image \, g$}
\put(40,5){$\sqcup$}
\put(85,5){$P_0$}
\put(130,5){$=$}
\put(150,5){$\tilde{P}$}
%
\put(0,50){$\scriptstyle g$}
\put(60,78){$\scriptstyle f_2^{23}$}
\put(110,78){$\scriptstyle f_2^{13}$}
\put(55,25){$\scriptstyle f_1^{12}$}
\put(110,25){$\scriptstyle f_1^{23}$}
\end{picture}
}
\vspace*{1.3cm}
\end{equation}
Here, $P_2^{\times}$ denotes $P_2(C_n; \mu, \lambda)^{\times}$, for instance.

For $\lambda/\mu$ is a skew diagram of $l(\lambda)=3$, 
there does not exist any $\bp\in P(C_n; \mu, \lambda)$ which possesses
more than two transposed pair of paths. Therefore, 
we have 
$P(C_n; \mu, \lambda) = P_0(C_n; \mu, \lambda) \sqcup P_1(C_n; \mu, \lambda) 
\sqcup P_2(C_n; \mu, \lambda)$.
Let $P_2(C_n; \mu, \lambda)^{\times}$ be the subset 
of $P_2(C_n; \mu, \lambda)$ (see Figure \ref{fig:C-PP}) 
which consists of all $\bp \in P_2(C_n; \mu, \lambda)$
such that the point $u':= u+(-1,1)$ and $v' := v+(1,-1)$
are on $\bp=(p_1,p_2,p_3)$, where $u$ (resp.\ $v$) is 
the leftmost intersecting point of $(p_1,p_3)$ 
(resp.\ the rightmost intersecting point of $(p_2,p_3)$). 
Let $P_1^{ij}(C_n; \mu, \lambda)$ ($1 \le i< j \le 3$) be the set of all 
$\bp=(p_1, p_2, p_3)\in P_1(C_n; \mu, \lambda)$ such that $(p_i,p_j)$ is 
transposed. Let 
$P_2(C_n; \mu, \lambda)^{\circ}:=
P_2(C_n; \mu, \lambda)\setminus P_2(C_n; \mu, \lambda)^{\times}$. 
Let
\begin{align*}
f_2^{ij} & : P_2(C_n; \mu, \lambda)^{\circ}
\to P_1(C_n; \mu, \lambda), \quad \text{$(i,j)=(1,3)$, $(2,3)$}, \\
f_1^{ij} & : P_1^{ij}(C_n; \mu, \lambda)
\to P_0(C_n; \mu, \lambda), \quad \text{$(i,j)=(1,2)$, $(2,3)$}
\end{align*}
be the maps that resolve the transposed pair $(p_i,p_j)$
in $\bp =(p_1,p_2,p_3)\in P_k(C_n; \mu, \lambda)$, which are 
defined in Appendix \ref{sec:3-row} (see also \eqref{eq:C-maps}). 
These maps are weight-preserving injections
(Lemmas \ref{lem:f2-map} and \ref{lem:f1-map}).
We remark that the set 
$P_2(C_n; \mu, \lambda)^{\circ}$ consists of all 
$\bp$ in $P_2(C_n; \mu, \lambda)$ 
such that $f_2^{13}$ or $f_2^{23}$ is well-defined
(in fact, both of them are well-defined), while 
$P_2(C_n; \mu, \lambda)^{\times}$ consists of all 
$\bp \in P_2(C_n; \mu, \lambda)$ 
such that both $f_2^{13}$ and $f_2^{23}$ are not well-defined.

By Lemma \ref{lem:ff}, 
we have 
$\Image \, f^{12}_1 \cap \Image \, f^{23}_1 = \Image \, (f^{23}_1 \circ f^{13}_2)
= \Image \, (f^{12}_1 \circ f^{23}_2)$, 
and therefore, 
\begin{equation}\label{eq:minus}
-\sum_{\bp \in P_1(C_n; \mu, \lambda)}z^{\bp}_a 
+ \sum_{\bp \in P_2(C_n; \mu, \lambda)^{\circ}}
z^{\bp}_a
= -\sum_{\bp \in \Image \, f^{12}_1 \cup \Image \, f^{23}_1}z^{\bp}_a.
\end{equation}

Let $g:P_2(C_n; \mu, \lambda)^{\times} \to \tilde{P}(C_n; \mu, \lambda)$ 
be the weight-preserving injection defined in Section \ref{sec:3-row} 
(see also Figure \ref{fig:C-PP}). 
By Lemma \ref{lem:g-map} \eqref{item:g-map}, we have 
\begin{equation}\label{eq:plus}
\sum_{\bp \in P_0(C_n; \mu, \lambda)\sqcup P_2(C_n; \mu, \lambda)^{\times}}z^{\bp}_a
= \sum_{\bp \in \tilde{P}(C_n; \mu, \lambda)}z^{\bp}_a.
\end{equation}
Combining \eqref{eq:minus} and \eqref{eq:plus}, we obtain 
\begin{align*}
\chi_{\lambda/\mu,a}
& = \sum_{\bp \in \tilde{P}(C_n; \mu, \lambda)} z^{\bp}_a - 
\sum_{\bp \in \Image \, f^{12}_1 \cup \Image \, f^{23}_1}z^{\bp}_a\\
& = \sum_{T \in \widetilde{\Tab}(C_n, \lambda/\mu)} z^{T}_a -
\sum_{\bp \in \Image \, f^{12}_1 \cup \Image \, f^{23}_1}z^{T(\bp)}_a. 
\end{align*}
By Lemma \ref{lem:p-T}, 
the set $\{ T(\bp) \mid  \bp \in \Image \, f^{12}_1 \cup \Image \, f^{23}_1 \}$ 
consists of all the tableaux in 
$\widetilde{\Tab}(C_n, \lambda/\mu)$ prohibited by the extra rules 
$(\bEtwoR )$ and $(\bEthreeR)$.
\end{proof}

\subsection{Skew diagrams of at most two columns}\label{sec:two-column}
In this subsection, we conjecture the tableaux description  
for skew diagrams $\lambda/\mu$ 
of at most two columns, and prove it for $l(\lambda)\le 4$. 
We assume that $l(\lambda) \le n+1$. 


{\it The case of one-column\/}.
Let $\lambda/\mu$ be a skew diagram of one column (i.e., $l(\lambda')=1$). 
Let $\Tab(C_n, \lambda/\mu)$ be the set of all the $HV$-tableaux $T$
(actually, the horizontal rule $(\bH )$ is not required) 
with the following condition: 
\newline
\begin{tabular}{l}
$(\bEoneC)$ \quad  
\end{tabular}
If $T$ contains a subtableau 
\begin{equation}\label{eq:arr-V1}
\text{
\begin{picture}(10,23)(0,22)
\multiput(0,0)(14,0){2}{\line(0,1){45}}
\multiput(0,0)(0,14){2}{\line(1,0){14}}
\multiput(0,31)(0,14){2}{\line(1,0){14}}
\put(3,4){$c_l$}
\put(5,17){$\scriptstyle \vdots$}
\put(3,35){$c_1$}
\end{picture}
}
\vspace{19pt}
\end{equation}
such that $l\ge 2$, $c_1=c$ and $c_l=\overline{c}$
for some $1 \le c \le n$, then $l-1 \le n-c$. 

The following theorem is due to \cite{KS}. 
We reproduce it using the paths description. 

\begin{thm}[\cite{KS}]\label{thm:one-column}
For any skew diagram $\lambda/\mu$ of $l(\lambda')=1$ and $l(\lambda) \le n+1$, 
the equality \eqref{eq:C-tab} holds. 
\end{thm}

\begin{proof}
By $l(\lambda') = 1$, there does not exist any 
$\bp\in P(C_n; \mu, \lambda)$ which contains 
more than one transposing pair of paths, and therefore, 
we have $P(C_n; \mu, \lambda) = P_0(C_n; \mu, \lambda)\sqcup P_1(C_n; \mu, \lambda)$. 
We can define a weight-preserving, sign-inverting injection 
(which is well-defined if $l(\lambda) \le n+1$)
\begin{equation}\label{eq:one-column-f}
f_1 : P_1(C_n; \mu, \lambda) \to P_0(C_n; \mu, \lambda),
\end{equation}
using the maps $r_y^{ij}$ in Appendix \ref{sec:r}
(see also Figure \ref{fig:one-column}), 
and show that the set $\{ T(\bp) \mid  \bp \in \Image \, f_1 \}$ consists of all the 
tableaux in $\widetilde{\Tab}(C_n, \lambda/\mu)$ prohibited by the $(\bEoneC)$ rule. 
\end{proof}

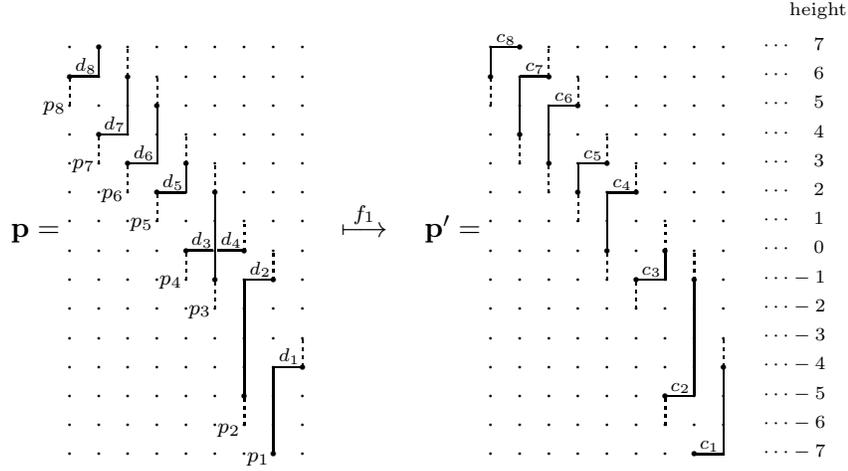
\begin{figure}
$\bp=$ 
\begin{picture}(88,90)(0,82)
\put(-10,130){$\scriptstyle p_8$}
\put(1,108){$\scriptstyle p_7$}
\put(12,97){$\scriptstyle p_6$}
\put(23,86){$\scriptstyle p_5$}
\put(34,64){$\scriptstyle p_4$}
\put(45,53){$\scriptstyle p_3$}
\put(56,8){$\scriptstyle p_2$}
\put(67,-3){$\scriptstyle p_1$}
\multiput(0,0)(0,11){15}{\circle*{0.5}}
\multiput(11,0)(0,11){15}{\circle*{0.5}}
\multiput(22,0)(0,11){15}{\circle*{0.5}}
\multiput(33,0)(0,11){15}{\circle*{0.5}}
\multiput(44,0)(0,11){15}{\circle*{0.5}}
\multiput(55,0)(0,11){15}{\circle*{0.5}}
\multiput(66,0)(0,11){15}{\circle*{0.5}}
\multiput(77,0)(0,11){15}{\circle*{0.5}}
\multiput(88,0)(0,11){15}{\circle*{0.5}}
%
%
\multiput(0,132)(0,3){4}{\line(0,1){1}}
\put(0,143){\line(1,0){11}}
\put(11,143){\line(0,1){11}}
\multiput(11,110)(0,3){4}{\line(0,1){1}}
\put(11,121){\line(1,0){11}}
\put(22,121){\line(0,1){22}}
\multiput(22,143)(0,3){4}{\line(0,1){1}}
\multiput(22,99)(0,3){4}{\line(0,1){1}}
\put(22,110){\line(1,0){11}}
\put(33,110){\line(0,1){22}}
\multiput(33,132)(0,3){4}{\line(0,1){1}}
\multiput(33,88)(0,3){4}{\line(0,1){1}}
\put(33,99){\line(1,0){11}}
\put(44,99){\line(0,1){11}}
\multiput(44,110)(0,3){4}{\line(0,1){1}}
\multiput(44,66)(0,3){4}{\line(0,1){1}}
\put(44,77){\line(1,0){10}}
\put(55,77){\line(0,1){22}}
\multiput(55,99)(0,3){4}{\line(0,1){1}}
\multiput(55,55)(0,3){4}{\line(0,1){1}}
\put(55,66){\line(0,1){11}}
\put(56,77){\line(1,0){10}}
\multiput(66,77)(0,3){4}{\line(0,1){1}}
\multiput(66,11)(0,3){4}{\line(0,1){1}}
\put(66,22){\line(0,1){44}}
\put(66,66){\line(1,0){11}}
\multiput(77,66)(0,3){4}{\line(0,1){1}}
\put(77,0){\line(0,1){33}}
\put(77,33){\line(1,0){11}}
\multiput(88,33)(0,3){4}{\line(0,1){1}}
\put(0,143){\circle*{1.5}}
\put(11,154){\circle*{1.5}}
\put(11,121){\circle*{1.5}}
\put(22,143){\circle*{1.5}}
\put(22,110){\circle*{1.5}}
\put(33,132){\circle*{1.5}}
\put(33,99){\circle*{1.5}}
\put(44,110){\circle*{1.5}}
\put(44,77){\circle*{1.5}}
\put(55,99){\circle*{1.5}}
\put(55,66){\circle*{1.5}}
\put(66,77){\circle*{1.5}}
\put(66,22){\circle*{1.5}}
\put(77,66){\circle*{1.5}}
\put(77,0){\circle*{1.5}}
\put(88,33){\circle*{1.5}}
{\tiny 
\put(79,35){$d_1$}
\put(68,68){$d_2$}
\put(57,79){$d_4$}
\put(46,79){$d_3$}
\put(35,101){$d_5$}
\put(24,112){$d_6$}
\put(13,123){$d_7$}
\put(2,145){$d_8$}
}
\end{picture}
$\quad \overset{f_1}{\longmapsto} \quad \bp'=$
\begin{picture}(150,90)(0,82)
\text{\tiny  
\put(103,-1){$\cdots -7$}
\put(103,10){$\cdots -6$}
\put(103,21){$\cdots -5$}
\put(103,32){$\cdots -4$}
\put(103,43){$\cdots -3$}
\put(103,54){$\cdots -2$}
\put(103,65){$\cdots -1$}
\put(103,76){$\cdots \quad 0$}
\put(103,87){$\cdots \quad 1$}
\put(103,98){$\cdots \quad 2$}
\put(103,109){$\cdots \quad 3$}
\put(103,120){$\cdots \quad 4$}
\put(103,131){$\cdots \quad 5$}
\put(103,142){$\cdots \quad 6$}
\put(103,153){$\cdots \quad 7$}
\put(113,166){height}
\put(79,2){$c_1$}
\put(68,24){$c_2$}
\put(57,68){$c_3$}
\put(46,101){$c_4$}
\put(35,112){$c_5$}
\put(24,134){$c_6$}
\put(13,145){$c_7$}
\put(2,156){$c_8$}
%
}
\multiput(0,0)(0,11){15}{\circle*{0.5}}
\multiput(11,0)(0,11){15}{\circle*{0.5}}
\multiput(22,0)(0,11){15}{\circle*{0.5}}
\multiput(33,0)(0,11){15}{\circle*{0.5}}
\multiput(44,0)(0,11){15}{\circle*{0.5}}
\multiput(55,0)(0,11){15}{\circle*{0.5}}
\multiput(66,0)(0,11){15}{\circle*{0.5}}
\multiput(77,0)(0,11){15}{\circle*{0.5}}
\multiput(88,0)(0,11){15}{\circle*{0.5}}
\multiput(0,132)(0,3){4}{\line(0,1){1}}
\put(0,143){\line(0,1){11}}
\put(0,154){\line(1,0){11}}
\multiput(11,110)(0,3){4}{\line(0,1){1}}
\put(11,121){\line(0,1){22}}
\put(11,143){\line(1,0){11}}
\multiput(22,143)(0,3){4}{\line(0,1){1}}
\multiput(22,99)(0,3){4}{\line(0,1){1}}
\put(22,110){\line(0,1){22}}
\put(22,132){\line(1,0){11}}
\multiput(33,132)(0,3){4}{\line(0,1){1}}
\multiput(33,88)(0,3){4}{\line(0,1){1}}
\put(33,99){\line(0,1){11}}
\put(33,110){\line(1,0){11}}
\multiput(44,110)(0,3){4}{\line(0,1){1}}
\multiput(44,66)(0,3){4}{\line(0,1){1}}
\put(44,77){\line(0,1){22}}
\put(44,99){\line(1,0){11}}
\multiput(55,99)(0,3){4}{\line(0,1){1}}
\multiput(55,55)(0,3){4}{\line(0,1){1}}
\put(55,66){\line(1,0){11}}
\put(66,66){\line(0,1){11}}
\multiput(66,77)(0,3){4}{\line(0,1){1}}
\multiput(66,11)(0,3){4}{\line(0,1){1}}
\put(66,22){\line(1,0){11}}
\put(77,22){\line(0,1){44}}
\multiput(77,66)(0,3){4}{\line(0,1){1}}
\put(77,0){\line(1,0){11}}
\put(88,0){\line(0,1){33}}
\multiput(88,33)(0,3){4}{\line(0,1){1}}
\put(0,143){\circle*{1.5}}
\put(11,154){\circle*{1.5}}
\put(11,121){\circle*{1.5}}
\put(22,143){\circle*{1.5}}
\put(22,110){\circle*{1.5}}
\put(33,132){\circle*{1.5}}
\put(33,99){\circle*{1.5}}
\put(44,110){\circle*{1.5}}
\put(44,77){\circle*{1.5}}
\put(55,99){\circle*{1.5}}
\put(55,66){\circle*{1.5}}
\put(66,77){\circle*{1.5}}
\put(66,22){\circle*{1.5}}
\put(77,66){\circle*{1.5}}
\put(77,0){\circle*{1.5}}
\put(88,33){\circle*{1.5}}
\end{picture}
\vspace{2.8cm}
\caption{An Example of $\bp\in P_1(C_n; \mu, \lambda)$ 
for one-column $\lambda/\mu$. 
For this $\bp$, the map 
$f_1$ in \eqref{eq:one-column-f} is given by 
$f_1=r_6^{18}\circ r_5^{17}\circ r_4^{16}\circ 
r_4^{27}\circ r_3^{26}\circ 
r_2^{25}\circ r_1^{24}\circ r_0^{34}$.
The tableau $T(\bp')$ does not satisfy $(\bEoneC)$.
}\label{fig:one-column}
\end{figure}

{\it The case of two-column\/}.
Let $\lambda/\mu$ be a skew diagram of two columns, i.e., $l(\lambda') = 2$.
Let $T\in \widetilde{\Tab}(C_n, \lambda/\mu)$ be a tableau 
that contains a subtableau 
\begin{equation}\label{eq:subtableau}
T'=
\text{
\begin{picture}(10,23)(0,22)
\multiput(0,0)(14,0){2}{\line(0,1){45}}
\multiput(0,0)(0,14){2}{\line(1,0){14}}
\multiput(0,31)(0,14){2}{\line(1,0){14}}
\put(3,4){$c_l$}
\put(5,17){$\scriptstyle \vdots$}
\put(3,35){$c_1$}
\end{picture}
}
\subset T
\vspace{22pt}
\end{equation}
such that $l \ge 2$, 
$c_1=n+2-l$, $c_l=\overline{n+2-l}$  
and every proper subtableau of $T'$ 
satisfies $(\bEoneC)$.
Let $\tilde{\lambda}/\tilde{\mu}$ be the one-column shape of $T'$. 
Then we can pick $\bp \in P_0(C_n; \tilde{\mu}, \tilde{\lambda})$ such that 
$T(\bp) = T'$. 
For $T'$ does not satisfy the extra rule $(\bEoneC)$, 
we have $\bp \in \Image \, f_1$, where $f_1$ is the injection \eqref{eq:one-column-f}
in the proof of Theorem \ref{thm:one-column}. 
Let $f_1^{-1}(\bp)=(p_1, \dots, p_l)\in P_1(C_n; \tilde{\mu}, \tilde{\lambda})$ be 
the inverse image of $\bp$. 
Then set (see Figure \ref{fig:one-column})
\begin{equation}\label{eq:c-lab}
d_i=d_i(T'):=
\begin{cases}
L_a^1(s^i), & \qquad i=1,\dots, l, \quad i \ne k, k+1, \\
\overline{n}, & \qquad i=k, \\
n, & \qquad i=k+1, 
\end{cases}
\end{equation}
where $L_a^1(s^i)$ is the $h$-label (defined in \eqref{eq:h-label-C}) of 
the unique eastward step $s^i$ in $p_i$, and 
$k$ is the number such that $c_k \preceq n$ and $c_{k+1} \succeq \overline{n}$. 
Then, one can show that
\begin{align}
\{c_1,\dots,c_k\}\cup\{\overline{d}_{k+2},\dots,\overline{d}_{l}\}
&=\{n,n-1,\dots,n+2-l\},\\
\{\overline{d}_1,\dots,\overline{d}_{k-1}\}
\cup\{{c}_{k+1},\dots,{c}_{l}\}
&=\{\overline{n},\overline{n-1},\dots,\overline{n+2-l}\}.
\end{align}
For example, these elements for all $T'$ as in \eqref{eq:subtableau}
of $l\le 4$
are given in Table \ref{tab:d}.

Now we define 
$\Tab(C_n, \lambda/\mu)$ as the set of all the $HV$-tableaux $T$ 
with the following condition: 
\newline
\begin{tabular}{l}
$(\bEtwoC )$ \quad 
\end{tabular}
Let $T'$ be any subtableau of $T$ 
(excluding $a_1, \dots , a_k, b_{k+1}, \dots , b_l$)
\begin{equation}\label{eq:arr-V2}
T'=
\text{
\begin{picture}(55,44)(-20,44)
\multiput(0,0)(20,0){2}{\line(0,1){88}}
\multiput(0,0)(0,14){2}{\line(1,0){20}}
\multiput(0,30)(0,14){3}{\line(1,0){20}}
\multiput(0,74)(0,14){2}{\line(1,0){20}}
\put(-13,4){$b_l$}
\put(-11,17){$\vdots$}
\put(-19,35){$b_{\scriptscriptstyle k+1}$}
\put(5,4){$c_l$}
\put(7,17){$\vdots$}
\put(1,35){$c_{\scriptscriptstyle k+1}$}
\put(5,48){$c_k$}
\put(7,61){$\vdots$}
\put(5,78){$c_1$}
\put(24,48){$a_k$}
\put(26,61){$\vdots$}
\put(24,78){$a_1$}
\end{picture}
}
\subset T
\vspace{45pt}
\end{equation}
such that $l \ge 2$,
$c_1=n+2-l$, $c_l=\overline{n+2-l}$, 
$c_k \preceq n$, $c_{k+1}\succeq \overline{n}$, and 
every proper subtableau in $T'$ 
satisfies the extra condition $(\bEoneC)$.   
Let $(i_1,j_1)$ be the position of the top of the subtableau $T'$
in \eqref{eq:arr-V2} (i.e., the position of $c_1$). 
Then one of the following conditions holds:
\begin{enumerate}
\item
$(i_1+i-1, j_1-1) \in \lambda/\mu$ and 
$a_i:=T(i_1+i-1, j_1+1) \prec d_i(T')$ for some $1 \le i \le k$.
\item
$(i_1+i-1, j_1+1) \in \lambda/\mu$ and 
$b_i:=T(i_1+i-1, j_1-1) \succ d_i(T')$ for some $k+1 \le i \le l$. 
\end{enumerate}

We remark that the extra rule $(\bEtwoC)$
is reduced to the extra rule $(\bEoneC)$, if $l(\lambda')=1$.   

We conjecture that 
\begin{conj}\label{conj:two-column}
For any skew diagram $\lambda/\mu$ of $l(\lambda')=2$ and $l(\lambda) \le n+1$,
the equality \eqref{eq:C-tab} holds. 
\end{conj}

\begin{thm}\label{thm:two-column}
Conjecture \ref{conj:two-column} is true for $l(\lambda) \le 4$. 
\end{thm}

\begin{proof}
If $l(\lambda) \le 3$, then 
the extra rule $(\bEtwoC)$ for $l(\lambda)\le 3$ 
coincides with the extra rule $(\bEtwoR)$ with 
$(\bEthreeR)$. 
For a summary of the proof for $l(\lambda) =4$, which is parallel to 
that of Theorem \ref{thm:3-row}, 
see the maps and their relations in the following diagram: 
\begin{equation}\label{eq:lambda-4}
\text{
\begin{picture}(255,65)(35,65)
\put(30,123){$P_2^{\times}$}
\put(60,123){$\sqcup $}
\put(85,123){$(P^{\scriptscriptstyle 13;23}_2)^{\circ}$}
\put(135,123){$\sqcup$}
\put(153,123){$(P^{\scriptscriptstyle 12;34}_2)^{\circ}$}
\put(203,123){$\sqcup$}
\put(220,123){$(P^{\scriptscriptstyle 24;34}_2)^{\circ}$}
\put(270,123){$=$}
\put(290,123){$P_2$}
\put(102,62){$P^{12}_1$}
\put(138,62){$\sqcup$}
\put(165,62){$P^{23}_1$}
\put(202,62){$\sqcup$}
\put(227,62){$P^{34}_1$}
\put(270,62){$=$}
\put(290,62){$P_1$}
\put(25,0){$\Image \, g$}
\put(90,0){$\sqcup$}
\put(165,0){$P_0$}
\put(230,0){$=$}
\put(270,0){$\tilde{P}$}
\put(35,115){\vector(0,-1){100}}
\put(40,55){$\scriptstyle g$}
\put(105,115){\vector(0,-1){40}}
\put(120,115){\vector(1,-1){40}}
\put(160,115){\vector(-1,-1){40}}
\put(180,115){\vector(1,-1){40}}
\put(220,115){\vector(-1,-1){40}}
\put(235,115){\vector(0,-1){40}}
\put(117,55){\vector(1,-1){42}}
\put(170,55){\vector(0,-1){42}}
\put(223,55){\vector(-1,-1){42}}
\put(90,95){$\scriptstyle f_2^{23}$}
\put(115,100){$\scriptstyle f_2^{13}$}
\put(150,100){$\scriptstyle f_2^{34}$}
\put(175,100){$\scriptstyle f_2^{12}$}
\put(210,100){$\scriptstyle f_2^{34}$}
\put(240,95){$\scriptstyle f_2^{24}$}
\put(110,35){$\scriptstyle f_1^{12}$}
\put(175,35){$\scriptstyle f_1^{23}$}
\put(215,35){$\scriptstyle f_1^{34}$}
\end{picture}
}
\vspace*{2.2cm}
\end{equation}
Here, $(P_2^{\scriptscriptstyle 13;23})^{\circ}$ denotes 
$P_2^{\scriptscriptstyle 13;23}(C_n; \mu, \lambda)^{\circ}$, for instance.

There does not exist 
$\bp\in P(C_n; \mu, \lambda)$ 
that contains more than two transposed pair 
of paths, and therefore, 
\begin{equation*}
P(C_n; \mu, \lambda)=P_0(C_n; \mu, \lambda) \sqcup 
P_1(C_n; \mu, \lambda) \sqcup P_2(C_n; \mu, \lambda).
\end{equation*} 
As in the proof of Theorem \ref{thm:3-row}, 
we define $P_1^{ij}(C_n; \mu, \lambda)$ 
($1\le i<j \le 4$) as the set of all 
$\bp = (p_1, \dots, p_4) \in P_1(C_n; \mu, \lambda)$ 
such that $(p_i, p_j)$ is transposed. 
Then we have  
\begin{equation*}
P_1(C_n; \mu, \lambda)=P_1^{12}(C_n; \mu, \lambda)\sqcup 
P_1^{23}(C_n; \mu, \lambda)\sqcup P_1^{34}(C_n; \mu, \lambda). 
\end{equation*}
Similarly, we define 
$P_2^{\scriptscriptstyle ij;km}(C_n; \mu, \lambda)$ (
$1\le i<j \le 4$, $1\le k<m \le 4$) as the set of all 
$\bp = (p_1, \dots, p_4) \in P_2(C_n; \mu, \lambda)$ 
such that $(p_{i},p_{j})$ and $(p_{k},p_{m})$ are transposed. 
Then we have 
\begin{equation*}
P_2(C_n; \mu, \lambda)=P_2^{\scriptscriptstyle 13;23}(C_n; \mu, \lambda)\sqcup 
P_2^{\scriptscriptstyle 12;34}(C_n; \mu, \lambda)\sqcup 
P_2^{\scriptscriptstyle 24;34}(C_n; \mu, \lambda). 
\end{equation*}
Let $P_2^{\scriptscriptstyle 13;23}(C_n; \mu, \lambda)^{\times}$ be the set 
that consists of all $\bp =(p_1, \dots, p_4)\in 
P_2^{\scriptscriptstyle 13;23}(C_n; \mu, \lambda)$ 
which satisfy one of the following conditions, where 
$u$ (resp. $v=u-(2,0)$) is the unique intersecting point of $(p_1,p_3)$
(resp. $(p_2,p_3)$) at height 0:
\begin{enumerate} 
\item Both points $u+(-1,1)$ and $v+(1,-1)$ are on $\bp$. 
\item All four points $u+(-1,1)$, $v+(1,-2)$, $v+(-1,1)$ and $v+(0,2)$ are on $\bp$.   
\end{enumerate}
Let $P_2^{\scriptscriptstyle 24;34}(C_n; \mu, \lambda)^{\times}$ be the set
of all $\bp \in P_2^{\scriptscriptstyle 24;34}(C_n; \mu, \lambda)$ such that 
$\omega(\bp) \in 
P_2^{\scriptscriptstyle 13;23}(C_n; \tilde{\mu}, \tilde{\lambda})^{\times}$, 
where $\omega$ is a map that rotates $\bp$ by 180 degrees 
defined as in \eqref{eq:omega}. 
Let $P_2^{\scriptscriptstyle 12;34}(C_n; \mu, \lambda)^{\times}$ be the set
that consists of all $\bp =(p_1, \dots, p_4)\in P_2^{12;34}(C_n; \mu, \lambda)$
such that all four points $u+(-1,1)$, $u+(-2,2)$, $w+(1,-1)$ and $w+(2,-2)$ 
are on $\bp$, 
where
$u$ and (resp. $w=u-(3,0)$) is the unique intersecting point of $(p_1,p_2)$
(resp. $(p_3,p_4)$) at height 0. 
Let $P_2^{\scriptscriptstyle 13;23}(C_n; \mu, \lambda)^{\circ}:=
P_2^{\scriptscriptstyle 13;23}(C_n; \mu, \lambda) \setminus 
P_2^{\scriptscriptstyle 13;23}(C_n; \mu, \lambda)^{\times}$, 
etc. 
We can define a weight-preserving, sign-inverting injection 
\begin{align*}
& f_2^{ij}:P_2^{\scriptscriptstyle km;k'm'}(C_n; \mu, \lambda)^{\circ}
\to P_1(C_n; \mu, \lambda), \quad (i,j) =(k,m), (k',m'), \\
& f_1^{ij}:P_1^{ij}(C_n; \mu, \lambda) \to P_0(C_n; \mu, \lambda),
\end{align*}
which resolve the transposed pair $(p_i,p_j)$ 
of $\bp =(p_1, \dots, p_4)\in P_k(C_n; \mu, \lambda)$ as  
the maps $f_2^{ij}$, $f_1^{ij}$ in the proof of Theorem \ref{thm:3-row}. 
These maps can also be defined as the composition of the maps $r^{ij}_y$ 
given in Section \ref{sec:r}. 
We remark that $P_2^{\scriptscriptstyle km;k'm'}(C_n; \mu, \lambda)^{\circ}$ 
consists of all $\bp \in P_2^{\scriptscriptstyle km;k'm'}(C_n; \mu, \lambda)$ 
such that 
$f_2^{ij}$ for some $1 \le i<j\le 4$ is well-defined
(in fact, all $f_2^{ij}$ are well-defined), while 
$P_2^{\scriptscriptstyle km;k'm'}(C_n; \mu, \lambda)^{\times}$ consists of all 
$\bp \in P_2^{\scriptscriptstyle km;k'm'}(C_n; \mu, \lambda)$ such that $f_2^{ij}$ 
for any $1 \le i<j \le 4$ is not well-defined. 
These maps satisfy 
\begin{gather*}
\Image \, (f_1^{34}\circ f_2^{12}) 
= \Image \, (f_1^{12}\circ  f_2^{34})
= \Image \, f_1^{12} \cap \Image \, f_1^{34}, \\
\Image \, (f_1^{23}\circ f_2^{13})  
= \Image \, (f_1^{12}\circ f_2^{23})
= \Image \, f_1^{12} \cap \Image \, f_1^{23}, \\
\Image \, (f_1^{34}\circ f_2^{24})  
= \Image \, (f_1^{23}\circ f_2^{34})
= \Image \, f_1^{23} \cap \Image \, f_1^{34},\\
\Image \, f_1^{12} \cap \Image \, f_1^{23} \cap \Image \, f_1^{34} = \phi, 
\end{gather*}
which can be proved by using the forms of the subtableaux in $T(\bp)$ 
of $\bp \in \Image \, f_1^{ij}$ for $(i,j)=(1,2),(2,3),(3,4)$ 
(see Table \ref{tab:d} and Figure 
\ref{fig:length-4-two-column-1}). 
We can also define a weight-preserving, sign-preserving injection 
\begin{equation*}
g: P_2(C_n; \mu, \lambda)^{\times}\to \tilde{P}(C_n; \mu, \lambda)
\end{equation*}
on $P_2(C_n; \mu, \lambda)^{\times}:=
P_2^{\scriptscriptstyle 13;23}(C_n; \mu, \lambda)^{\times}\sqcup 
P_2^{\scriptscriptstyle 12;34}(C_n; \mu, \lambda)^{\times}\sqcup 
P_2^{\scriptscriptstyle 24;34}(C_n; \mu, \lambda)^{\times}$,
which satisfies
$\Image \, g \sqcup P_0(C_n; \mu, \lambda) = \tilde{P}(C_n; \mu, \lambda)$, 
as in the proof of Theorem \ref{thm:3-row}.
Then we similarly obtain the equality \eqref{eq:C-tab} by the following lemma. 
\end{proof}

\begin{table}
\begin{center}
{\setlength{\unitlength}{0.45mm}
\begin{tabular}{c|cccccccc}
\noalign{\hrule height0.8pt}
\begin{picture}(7,50)(0,-20)
\put(0,0){$T'$}
\end{picture}
& 
\begin{picture}(15,50)(-5,-15)
\multiput(0,0)(0,10){3}{\line(1,0){10}}
\multiput(0,0)(10,0){2}{\line(0,1){20}}
\put(3,2){$\scriptstyle \overline{n}$}
\put(3,12){$\scriptstyle n$}
\end{picture}
&
\begin{picture}(21,50)(-5,-10)
\multiput(0,0)(0,10){4}{\line(1,0){16}}
\multiput(0,0)(16,0){2}{\line(0,1){30}}
\put(2,2){$\scriptstyle \overline{n-1}$}
\put(6,13){$\scriptstyle n$}
\put(2,23){$\scriptstyle n-1$}
%
\end{picture}
& 
\begin{picture}(21,50)(-5,-10)
\multiput(0,0)(0,10){4}{\line(1,0){16}}
\multiput(0,0)(16,0){2}{\line(0,1){30}}
\put(2,2){$\scriptstyle \overline{n-1}$}
\put(6,12){$\scriptstyle \overline{n}$}
\put(2,23){$\scriptstyle n-1$}
%
\end{picture}
& 
\begin{picture}(21,50)(-5,-5)
\multiput(0,0)(0,10){5}{\line(1,0){16}}
\multiput(0,0)(16,0){2}{\line(0,1){40}}
\put(2,2){$\scriptstyle \overline{n-2}$}
\put(6,13){$\scriptstyle n$}
\put(2,23){$\scriptstyle n-1$}
\put(2,33){$\scriptstyle n-2$}
%
\end{picture}
& 
\begin{picture}(21,50)(-5,-5)
\multiput(0,0)(0,10){5}{\line(1,0){16}}
\multiput(0,0)(16,0){2}{\line(0,1){40}}
\put(2,2){$\scriptstyle \overline{n-2}$}
\put(6,12){$\scriptstyle \overline{n}$}
\put(2,23){$\scriptstyle n-1$}
\put(2,33){$\scriptstyle n-2$}
%
\end{picture}
& 
\begin{picture}(21,50)(-5,-5)
\multiput(0,0)(0,10){5}{\line(1,0){16}}
\multiput(0,0)(16,0){2}{\line(0,1){40}}
\put(2,2){$\scriptstyle \overline{n-2}$}
\put(2,12){$\scriptstyle \overline{n-1}$}
\put(2,23){$\scriptstyle n-1$}
\put(2,33){$\scriptstyle n-2$}
%
\end{picture}
&
\begin{picture}(21,50)(-5,-5)
\multiput(0,0)(0,10){5}{\line(1,0){16}}
\multiput(0,0)(16,0){2}{\line(0,1){40}}
\put(2,2){$\scriptstyle \overline{n-2}$}
\put(2,12){$\scriptstyle \overline{n-1}$}
\put(6,23){$\scriptstyle n$}
\put(2,33){$\scriptstyle n-2$}
%
\end{picture}
& 
\begin{picture}(21,50)(-5,-5)
\multiput(0,0)(0,10){5}{\line(1,0){16}}
\multiput(0,0)(16,0){2}{\line(0,1){40}}
\put(2,2){$\scriptstyle \overline{n-2}$}
\put(2,12){$\scriptstyle \overline{n-1}$}
\put(6,22){$\scriptstyle \overline{n}$}
\put(2,33){$\scriptstyle n-2$}
%
\end{picture}
\\
\hline
\hfil $d_1$ & \hfil $\overline{n}$ \hfil & \hfil $n$ \hfil 
& \hfil $\overline{n}$ \hfil & \hfil $n-1$
& \hfil $n-1$ & \hfil $n$ \hfil & \hfil $n$ \hfil 
& \hfil $\overline{n}$ \hfil \\
\hfil $d_2$ & \hfil $n$ \hfil & \hfil $\overline{n}$ \hfil 
& \hfil $n$ \hfil & \hfil $n$ \hfil & \hfil $\overline{n}$ \hfil 
& \hfil $\overline{n}$ \hfil & \hfil $\overline{n}$ \hfil & \hfil $n$ \hfil \\
\hfil $d_3$ &  & \hfil $n$ \hfil & \hfil $\overline{n}$ \hfil 
& \hfil $\overline{n}$ \hfil & \hfil $n$ \hfil 
& \hfil $n$ \hfil & \hfil $n$ \hfil & \hfil $\overline{n}$ \hfil \\
\hfil $d_4$ &  &  &  & \hfil $n$ \hfil & \hfil $\overline{n}$ \hfil & 
\hfil $\overline{n}$ \hfil & \hfil $\overline{n-1}$  
& \hfil $\overline{n-1}$ \\
\noalign{\hrule height0.8pt}
\end{tabular}
}
\end{center}
\caption{The table of $(d_1, \dots, d_l)$ for one-column tableaux $T'$
as in \eqref{eq:subtableau} of $l\le 4$.}\label{tab:d}
\end{table}

\begin{figure}
{\setlength{\unitlength}{0.4mm}
\begin{tabular}{lclccl}
{\small \qquad \qquad \quad $\bp$} & & 
{\small $\bp'=f_1^{12}(\bp)$} & & 
{\small $T(\bp')$}\\
\hspace{1.3cm}
\begin{picture}(30,20)(0,10)
%
\multiput(-10,0)(0,10){3}{\circle*{0.5}}
\multiput(0,0)(0,10){3}{\circle*{0.5}}
\multiput(10,0)(0,10){3}{\circle*{0.5}}
\multiput(20,0)(0,10){3}{\circle*{0.5}}
\multiput(30,0)(0,10){3}{\circle*{0.5}}
{\tiny
\put(12,3){$p_1$}
\put(-8,8){$p_2$}
\put(-30,-2){$-1$}
\put(-29,8){$\ \ 0$}
\put(-29,18){$\ \ 1$}
\put(-31,28){height}
}
\put(0,0){\circle*{1.5}}
\put(10,0){\circle*{1.5}}
\put(20,20){\circle*{1.5}}
\put(10,20){\circle*{1.5}}
%
%
\put(0,0){\line(0,1){10}}
\put(0,10){\line(1,0){9}}
\put(11,10){\line(1,0){9}}
\put(20,10){\line(0,1){10}}
%
\put(10,0){\line(0,1){20}}
%
\multiput(-10,0)(3,0){4}{\line(1,0){1}}
\put(-7,-7){$\scriptstyle b_2$}
\multiput(20,20)(3,0){4}{\line(1,0){1}}
\put(23,22){$\scriptstyle a_1$}
\end{picture}
&
$\overset{f_1^{12}}{\mapsto}$
& 
\begin{picture}(30,20)(-10,10)
\multiput(-10,0)(0,10){3}{\circle*{0.5}}
\multiput(0,0)(0,10){3}{\circle*{0.5}}
\multiput(10,0)(0,10){3}{\circle*{0.5}}
\multiput(20,0)(0,10){3}{\circle*{0.5}}
\multiput(30,0)(0,10){3}{\circle*{0.5}}
\put(0,0){\circle*{1.5}}
\put(10,0){\circle*{1.5}}
\put(20,20){\circle*{1.5}}
\put(10,20){\circle*{1.5}}
%
%
\put(0,0){\line(0,1){20}}
\put(0,20){\line(1,0){10}}
%
\put(10,0){\line(1,0){10}}
\put(20,0){\line(0,1){20}}
\multiput(-10,0)(3,0){4}{\line(1,0){1}}
\put(-7,-7){$\scriptstyle b_2$}
\multiput(20,20)(3,0){4}{\line(1,0){1}}
\put(23,22){$\scriptstyle a_1$}
\end{picture}
& 
$\overset{T}{\mapsto}$ 
& 
\begin{picture}(24,16)(-10,10)
\multiput(0,0)(0,10){3}{\line(1,0){10}}
\multiput(0,0)(10,0){2}{\line(0,1){20}}
\put(3,3){$\scriptstyle \overline{n}$}
\put(3,13){$\scriptstyle n$}
\put(-10,3){$\scriptstyle b_2$}
\put(13,13){$\scriptstyle a_1$}
\end{picture}
& \quad 
{\small 
\begin{tabular}{l}
$a_1 \succeq \overline{n}$, \\
$b_2 \preceq n$
\end{tabular}
}
\\
\hspace{1.3cm}
\begin{picture}(30,50)(0,20)
%
\multiput(-10,0)(0,10){5}{\circle*{0.5}}
\multiput(0,0)(0,10){5}{\circle*{0.5}}
\multiput(10,0)(0,10){5}{\circle*{0.5}}
\multiput(20,0)(0,10){5}{\circle*{0.5}}
\multiput(30,0)(0,10){5}{\circle*{0.5}}
\multiput(40,0)(0,10){5}{\circle*{0.5}}
{\tiny
\put(22,5){$p_1$}
\put(3,18){$p_2$}
\put(2,36){$p_3$}
\put(-30,-2){$-2$}
\put(-30,8){$-1$}
\put(-29,18){$\ \ 0$}
\put(-29,28){$\ \ 1$}
\put(-29,38){$\ \ 2$}
}
\put(0,30){\circle*{1.5}}
\put(10,10){\circle*{1.5}}
\put(20,0){\circle*{1.5}}
\put(30,30){\circle*{1.5}}
\put(20,30){\circle*{1.5}}
\put(10,40){\circle*{1.5}}
%
%
\put(0,30){\line(1,0){10}}
\put(10,30){\line(0,1){10}}
\put(10,10){\line(0,1){10}}
\put(10,20){\line(1,0){9}}
\put(21,20){\line(1,0){9}}
\put(30,20){\line(0,1){10}}
\put(20,0){\line(0,1){30}}
\multiput(-10,30)(3,0){4}{\line(1,0){1}}
\put(-9,23){$\scriptstyle b_3$}
\multiput(0,10)(3,0){4}{\line(1,0){1}}
\put(1,3){$\scriptstyle b_2$}
\multiput(30,30)(3,0){4}{\line(1,0){1}}
\put(33,32){$\scriptstyle a_1$}
\end{picture}
&
$\overset{f_1^{12}}{\mapsto}$
& 
\begin{picture}(30,50)(-10,20)
\multiput(-10,0)(0,10){5}{\circle*{0.5}}
\multiput(0,0)(0,10){5}{\circle*{0.5}}
\multiput(10,0)(0,10){5}{\circle*{0.5}}
\multiput(20,0)(0,10){5}{\circle*{0.5}}
\multiput(30,0)(0,10){5}{\circle*{0.5}}
\multiput(40,0)(0,10){5}{\circle*{0.5}}
\put(0,30){\circle*{1.5}}
\put(10,10){\circle*{1.5}}
\put(20,0){\circle*{1.5}}
\put(30,30){\circle*{1.5}}
\put(20,30){\circle*{1.5}}
\put(10,40){\circle*{1.5}}
%
%
\put(0,30){\line(0,1){10}}
\put(0,40){\line(1,0){10}}
\put(10,10){\line(0,1){20}}
\put(10,30){\line(1,0){10}}
\put(20,0){\line(1,0){10}}
\put(30,0){\line(0,1){30}}
\multiput(-10,30)(3,0){4}{\line(1,0){1}}
\put(-9,23){$\scriptstyle b_3$}
\multiput(0,10)(3,0){4}{\line(1,0){1}}
\put(1,3){$\scriptstyle b_2$}
\multiput(30,30)(3,0){4}{\line(1,0){1}}
\put(33,32){$\scriptstyle a_1$}
\end{picture}
& 
$\overset{T}{\mapsto}$   
& 
\begin{picture}(30,30)(-10,15)
\multiput(0,0)(0,10){4}{\line(1,0){16}}
\multiput(0,0)(16,0){2}{\line(0,1){30}}
\put(1,2){$\scriptstyle \overline{n-1}$}
\put(6,13){$\scriptstyle \overline{n}$}
\put(1,23){$\scriptstyle n-1$}
\put(19,23){$\scriptstyle a_1$}
\put(-11,13){$\scriptstyle b_2$}
\put(-11,3){$\scriptstyle b_3$}
\end{picture}
& \quad 
{\small
\begin{tabular}{l}
$a_1 \succeq \overline{n}$, \\
$b_2 \preceq n$, \\
$b_3 \preceq \overline{n}$
\end{tabular}
}
\\
\hspace{1.3cm}
\begin{picture}(40,60)(10,30)
%
\multiput(-10,0)(0,10){7}{\circle*{0.5}}
\multiput(0,0)(0,10){7}{\circle*{0.5}}
\multiput(10,0)(0,10){7}{\circle*{0.5}}
\multiput(20,0)(0,10){7}{\circle*{0.5}}
\multiput(30,0)(0,10){7}{\circle*{0.5}}
\multiput(40,0)(0,10){7}{\circle*{0.5}}
\multiput(50,0)(0,10){7}{\circle*{0.5}}
{\tiny
\put(32,5){$p_1$}
\put(13,26){$p_2$}
\put(12,46){$p_3$}
\put(2,56){$p_4$}
\put(-30,-2){$-3$}
\put(-30,8){$-2$}
\put(-30,18){$-1$}
\put(-29,28){$\ \ 0$}
\put(-29,38){$\ \ 1$}
\put(-29,48){$\ \ 2$}
\put(-29,58){$\ \ 3$}
}
\put(0,50){\circle*{1.5}}
\put(10,40){\circle*{1.5}}
\put(20,20){\circle*{1.5}}
\put(30,0){\circle*{1.5}}
\put(40,40){\circle*{1.5}}
\put(30,40){\circle*{1.5}}
\put(20,50){\circle*{1.5}}
\put(10,60){\circle*{1.5}}
%
%
\put(0,50){\line(1,0){10}}
\put(10,50){\line(0,1){10}}
\put(10,40){\line(1,0){10}}
\put(20,40){\line(0,1){10}}
\put(20,20){\line(0,1){10}}
\put(20,30){\line(1,0){9}}
\put(31,30){\line(1,0){9}}
\put(40,30){\line(0,1){10}}
\put(30,0){\line(0,1){40}}
\multiput(-10,50)(3,0){4}{\line(1,0){1}}
\put(-9,43){$\scriptstyle b_4$}
\multiput(0,40)(3,0){4}{\line(1,0){1}}
\put(1,33){$\scriptstyle b_3$}
\multiput(10,20)(3,0){4}{\line(1,0){1}}
\put(11,13){$\scriptstyle b_2$}
\multiput(40,40)(3,0){4}{\line(1,0){1}}
\put(43,42){$\scriptstyle a_1$}
\end{picture}
&
$\overset{f_1^{12}}{\mapsto}$ 
&
\begin{picture}(55,60)(-10,30)
\multiput(-10,0)(0,10){7}{\circle*{0.5}}
\multiput(0,0)(0,10){7}{\circle*{0.5}}
\multiput(10,0)(0,10){7}{\circle*{0.5}}
\multiput(20,0)(0,10){7}{\circle*{0.5}}
\multiput(30,0)(0,10){7}{\circle*{0.5}}
\multiput(40,0)(0,10){7}{\circle*{0.5}}
\multiput(50,0)(0,10){7}{\circle*{0.5}}
\put(0,50){\circle*{1.5}}
\put(10,40){\circle*{1.5}}
\put(20,20){\circle*{1.5}}
\put(30,0){\circle*{1.5}}
\put(40,40){\circle*{1.5}}
\put(30,40){\circle*{1.5}}
\put(20,50){\circle*{1.5}}
\put(10,60){\circle*{1.5}}
%
%
\put(0,50){\line(0,1){10}}
\put(0,60){\line(1,0){10}}
\put(10,40){\line(0,1){10}}
\put(10,50){\line(1,0){10}}
\put(20,20){\line(0,1){20}}
\put(20,40){\line(1,0){10}}
\put(30,0){\line(1,0){10}}
\put(40,0){\line(0,1){40}}
\multiput(-10,50)(3,0){4}{\line(1,0){1}}
\put(-9,43){$\scriptstyle b_4$}
\multiput(0,40)(3,0){4}{\line(1,0){1}}
\put(1,33){$\scriptstyle b_3$}
\multiput(10,20)(3,0){4}{\line(1,0){1}}
\put(11,13){$\scriptstyle b_2$}
\multiput(40,40)(3,0){4}{\line(1,0){1}}
\put(43,42){$\scriptstyle a_1$}
\end{picture}
&
$\overset{T}{\mapsto}$  
&
\begin{picture}(30,30)(-10,15)
\multiput(0,0)(0,10){5}{\line(1,0){16}}
\multiput(0,0)(16,0){2}{\line(0,1){40}}
\put(1,2){$\scriptstyle \overline{n-2}$}
\put(1,12){$\scriptstyle \overline{n-1}$}
\put(6,23){$\scriptstyle \overline{n}$}
\put(1,33){$\scriptstyle n-2$}
\put(19,33){$\scriptstyle a_1$}
\put(-11,23){$\scriptstyle b_2$}
\put(-11,13){$\scriptstyle b_3$}
\put(-11,3){$\scriptstyle b_4$}
\end{picture}
&
\quad 
{\small 
\begin{tabular}{l}
$a_1 \succeq \overline{n}$, \\
$b_2 \preceq n$, \\
$b_3 \preceq \overline{n}$, \\
$b_4 \preceq \overline{n-1}$
\end{tabular}
}
\end{tabular}
}
\newline
\vspace{0.5cm}
\caption{
Examples of $\bp =(p_1, \dots, p_4)\in P_1^{12}(C_n; \mu, \lambda)$, 
the map 
$f_1^{12}:P_1^{12}(C_n; \mu, \lambda)\to P_0(C_n; \mu, \lambda)$, 
and the subtableaux of $T(f_1^{23}(\bp))$. 
If the $(E)$ step for $a_i$ (resp.\ $b_i$) exists, 
then $a_i$ (resp.\ $b_i$) satisfies the condition as above, which implies 
the corresponding tableau is prohibited by $(\bEtwoC)$.
}\label{fig:length-4-two-column-1}
\end{figure}
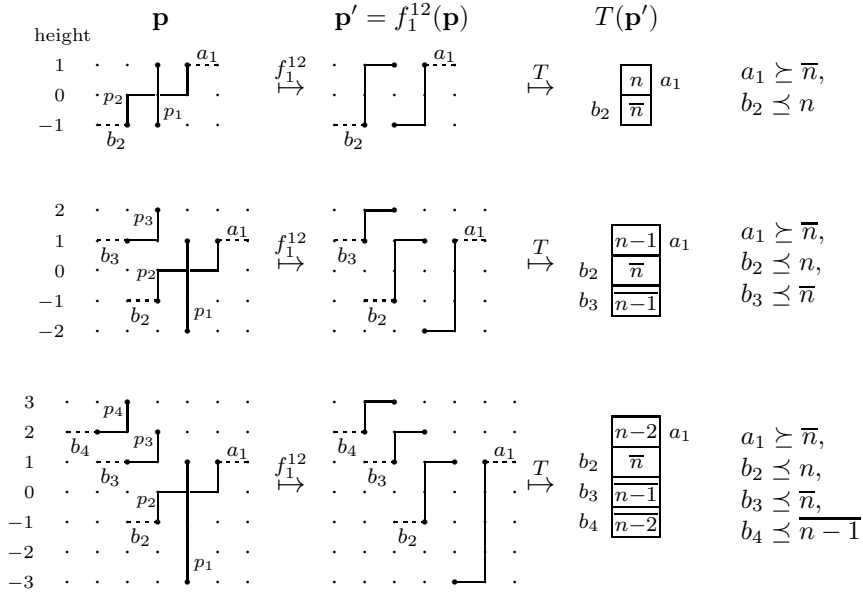

\begin{lem}\label{lem:2-column}
Let $\lambda/\mu$ be a skew diagram of $l(\lambda') \le 2$ 
and $l(\lambda)=n+1$. 
For $\bp \in \tilde{P}(C_n; \mu, \lambda)$, 
$\bp \in \Image \, f^{12}_1 \cup \Image \, f^{23}_1 \cup \Image \, f^{34}_1$ 
if and only if $T(\bp)$ is prohibited by $(\bEtwoC)$.  
(See Figure \ref{fig:length-4-two-column-1} for example.  )
\end{lem}

\appendix

\section{Classical projection of $\chi_{\lambda,a}$}\label{sec:g-ch}
In this section, we give a ``classical projection'' of 
the determinant $\chi_{\lambda,a}$ in \eqref{eq:det}, 
the one obtained by dropping the spectral parameters $a \in \bC$. 
We prove that the classical projection of $\chi_{\lambda,a}$ 
coincides with the character for the representation of $\Uqg$ 
defined in \cite{CK}. 

Let $\beta :\bZ[Y_{i,a}^{\pm 1}]_{i=1, \dots ,n; a\in \bC} 
\to \bZ[y_i^{\pm 1}]_{i=1, \dots ,n}$ be the classical projection, i.e.,
the algebra homomorphism defined by $\beta(Y_{i,a})=y_i$. 
Identifying $\cZ$ with $\cY$ by the isomorphism in the proof of 
Proposition \ref{prop:T-Y}, $\beta |_{\cZ}$ is the map 
$\cZ \to \bZ [z_i^{\pm 1}]_{i=1, \dots, N}$ 
such that 
$\beta(z_{i,a})=z_i$, 
$\beta(z_{0,a})=1$ (for $B_n$)
and $\beta(z_{\overline{i},a}) = z_i^{-1}$ (for $B_n$, $C_n$ and $D_n$). 
The homomorphism $\beta$ sends the $q$-character $\chi_q(V)$ 
for any finite dimensional representation $V$ of $\Uqhg$ 
to the character $\chi(V)$ of $\Uqg$ 
for $V$ as a $\Uqg$-module \cite{FR1}.

Let $\chi_{\lambda,a}$ be the determinant \eqref{eq:det}
with $\mu = \phi$. 
Let $\chi_{\lambda}\in \bZ[z_i^{\pm 1}]_{i=1, \dots, N}$ be the character of $\fg$
for the irreducible representation with highest weight $\lambda$. 
For any partitions $\mu ,\nu ,\lambda$,
let $c^{\lambda}_{\mu \nu}$ be the Littlewood-Richardson coefficient
\cite{M}.
Then we have 

%


\begin{thm}\label{thm:classical}
For any $\lambda$ such that $l(\lambda)\le n$,
\begin{equation}\label{eq:decomposition}
\beta(\chig_{\lambda,a})=
\begin{cases}
\chi_{\lambda}, 
& \text{\qquad if $\fg $ is of type $A_n$},\\
\sum_{\kappa ,\mu}c^{\lambda}_{(2\kappa)', \mu}\chi_{\mu},
& \text{\qquad \phantom{if $\fg $ is of type} $B_n$},\\
\sum_{\kappa ,\mu}c^{\lambda}_{2\kappa, \mu}\chi_{\mu},
& \text{\qquad \phantom{if $\fg $ is of type} $C_n$},\\
\sum_{\kappa ,\mu}c^{\lambda}_{(2\kappa)', \mu}\tilde{\chi}_{\mu},
& \text{\qquad \phantom{if $\fg $ is of type} $D_n$},
\end{cases}
\end{equation}
where
\begin{equation*}
\tilde{\chi}_{\lambda}:=
\begin{cases}
\chi_{\lambda}, 
& \text{if $1\le l(\lambda) \le n-1$},
\\
\chi_{\lambda} + \chi_{\sigma (\lambda)},
& \text{if $l(\lambda) = n$},
\end{cases}
\end{equation*}
for $D_n$,
where $\sigma$ is induced from the automorphism of the Dynkin diagram.
%
%
\end{thm}

\begin{proof}
Let $\Lambda$ be the graded ring of
symmetric functions with countable many variables $z_1, z_2, \dots $, 
and let $S_{\lambda}\in \Lambda$ be the {\it Schur function\/}. 
It is well-known that 
$S_{\lambda}$ satisfies the Jacobi-Trudi identity
\begin{equation*}
S_{\lambda}= \det (h_{\lambda_i -i + j})_{1 \le i,j \le l(\lambda)}
=\det (e_{\lambda'_i -i + j})_{1 \le i,j \le l(\lambda')}, 
\end{equation*}
where $h_i, e_i\in \Lambda$ are defined as 
\begin{equation*}
\prod_{k=1}^{\infty}(1-z_kx)^{-1} = \sum_{i=0}^{\infty}h_ix^i, \qquad  
\prod_{k=1}^{\infty}(1+z_kx) = \sum_{i=0}^{\infty}e_ix^i. 
\end{equation*}
Let $\varphi: \Lambda \to \Lambda $
be the algebra automorphism defined by
\begin{alignat*}{2}
& \varphi(e_i) = e_i-e_{i-2}, 
& \qquad 
& \varphi^{-1}(e_i) = \textstyle\sum_{m=0}^{\infty}e_{i-2m}, \\
& \varphi(h_i) = \textstyle\sum_{m=0}^{\infty}h_{i-2m}, 
& \qquad 
& \varphi^{-1}(h_i) = h_i-h_{i-2}.
\end{alignat*}
(All the four conditions are equivalent to each other.)

Set $\Lambda_n:=\bZ [z_1,\dots ,z_{n}]^{\fS_{n}}$. Let 
$\pi_{n+1}:\Lambda \to \Lambda_{n+1} /
\langle z_1\cdots z_{n+1}-1 \rangle $ be the map induced from 
the natural projection $\Lambda  \to \Lambda_{n+1}$
and let
$\pi_{\scriptscriptstyle Sp(2n)}$ and $\pi_{\scriptscriptstyle O(N)}$ 
be the {\it specialization homomorphisms\/}
in \cite{KT}. 
We define $\rho_n$ as 
\begin{equation*}
\rho_n = 
\begin{cases}
\pi_{n+1}, & \qquad (A_n)\\
\pi_{\scriptscriptstyle O(2n+1)} \circ \varphi^{-1}, & \qquad (B_n)\\
\pi_{\scriptscriptstyle Sp(2n)} \circ \varphi, & \qquad (C_n)\\
\pi_{\scriptscriptstyle O(2n)} \circ \varphi^{-1}. & \qquad (D_n)
\end{cases}
\end{equation*}
By the properties of $\pi_{\scriptscriptstyle O(N)}$ and 
$\pi_{\scriptscriptstyle Sp(2n)}$ \cite{KT} and the definitions of 
$h_{i,a}$ and $e_{i,a}$ in \eqref{eq:generating-H} and \eqref{eq:generating-E}, 
we have $\beta(\chig_{\lambda,a}) 
= \rho_n(S_{\lambda})$ for any Young diagram $\lambda$ of 
$l(\lambda) \le n$. 
Therefore, 
for $A_n$, \eqref{eq:decomposition} is obvious by the fact that 
$\pi_{n+1}(S_{\lambda})=\chi_{\lambda}$, 
while for $B_n$, $C_n$ and $D_n$, 
\eqref{eq:decomposition} 
are obtained by the equalities \cite{W, KT} 
\begin{gather*}
\prod_{i,j\ge 1}\frac{1}{(1-z_i\tilde{z}_j)} 
= \sum_{\lambda}S_{\lambda}(z)S_{\lambda}(\tilde{z}),
\\
\frac{\prod_{1\le i \le j}(1-\tilde{z}_i\tilde{z}_j)}
{\prod_{i,j\ge 1}(1-z_i\tilde{z}_j)}
=\sum_{\lambda}\chi_{\scriptscriptstyle O}(\lambda)(z)S_{\lambda}(\tilde{z}),\\
\frac{\prod_{1\le i < j}(1-\tilde{z}_i\tilde{z}_j)}
{\prod_{i,j\ge 1}(1-z_i\tilde{z}_j)}
=\sum_{\lambda}\chi_{\scriptscriptstyle Sp}(\lambda)(z)S_{\lambda}(\tilde{z}),
\end{gather*}
where 
$\chi_{\scriptscriptstyle Sp}(\lambda), 
\chi_{\scriptscriptstyle O}(\lambda)\in \Lambda$ 
are the {\it universal character\/}
of $Sp$ and $O$ \cite{KT}, 
and the Littlewood's Lemma \cite{L}
$$\prod_{1\le i \le j}\frac{1}{(1-z_iz_j)}
=\sum_{\kappa}S_{2\kappa}(z), \qquad 
\prod_{1\le i < j}\frac{1}{(1-z_iz_j)}
=\sum_{\kappa}S_{(2\kappa)'}(z).$$
\end{proof}

\begin{rem}\label{rem:CK}
The right hand side of \eqref{eq:decomposition}
is the character of the representation $W_G(\lambda)$
defined in \cite{CK}. 
Therefore, by Theorem \ref{thm:classical}, under the classical projection, 
Conjecture \ref{conj:det} reduces to 
Conjecture 2 in \cite{CK}
of the existence of an irreducible representation of $\Uqhg$, 
which is proved by \cite{C} 
for $\lambda = (i^m)$  such that $m\ge 1$ and 
$1\le i \le n$ ($A_n$ and $B_n$), 
$1\le i \le n-1$ ($C_n$),
$1\le i \le n-2$ ($D_n$).

\end{rem}

\section{The weight-preserving maps for $C_n$ Case}\label{sec:maps}
In this section, we define some weight-preserving
maps and give their properties 
which we use in the proof of Theorems \ref{thm:2-row} and \ref{thm:3-row}.

\subsection{The map $r_y$}\label{sec:r}
In this subsection, we give weight-preserving maps for 
a pair of $h$-paths of type $C_n$.  
These maps are used to define the maps 
in Section \ref{sec:3-row}. 
First, we define the map $r_y$ for $y=0,1, \dots , n-1$, 
which is defined on all $(p_1,p_2)\in P(C_n) \times P(C_n)$ 
that satisfy certain condition $(\bR_y)$. 

Set 
$(x,y)\pm (x',y'):= (x\pm x',y\pm y')$.

\begin{figure}
{\setlength{\unitlength}{1.1mm}
\begin{picture}(35,15)(-15,0)
\multiput(-15,-10)(5,0){8}{\circle*{0.5}}
\multiput(-15,-5)(5,0){8}{\circle*{0.5}}
\multiput(-15,0)(5,0){8}{\circle*{0.5}}
\multiput(-15,5)(5,0){8}{\circle*{0.5}}
\multiput(-15,10)(5,0){8}{\circle*{0.5}}
\put(-15,-10){\line(0,1){15}}
\put(-15,5){\line(1,0){20}}
\put(5,5){\line(0,1){5}}
\put(0,-10){\line(0,1){5}}
\put(0,-5){\line(1,0){20}}
\put(20,-5){\line(0,1){15}}
{\small
\put(7,-7){$p_1$}
\put(-14,-7){$p_2$}
\put(-11.2,4.2){$\times$}
\put(5,5){\circle*{1}}
\put(-10,6){$w_1^*$}
\put(5,6){$w_2$}
\put(0,-5){\circle*{1}}
\put(13.8,-5.8){$\times$}
\put(0,-4){$w_1$}
\put(15,-4){$w_2^*$}
}
\end{picture}
\quad $\overset{r_1}{\longmapsto}$ \quad 
\begin{picture}(35,18)(-15,0)
\multiput(-15,-10)(5,0){8}{\circle*{0.5}}
\multiput(-15,-5)(5,0){8}{\circle*{0.5}}
\multiput(-15,0)(5,0){8}{\circle*{0.5}}
\multiput(-15,5)(5,0){8}{\circle*{0.5}}
\multiput(-15,10)(5,0){8}{\circle*{0.5}}
\put(-15,-10){\line(0,1){15}}
\put(-15,5){\line(1,0){5}}
\put(-10,5){\line(0,1){5}}
\put(-10,10){\line(1,0){15}}
\put(0,-10){\line(1,0){15}}
\put(15,-10){\line(0,1){5}}
\put(15,-5){\line(1,0){5}}
\put(20,-5){\line(0,1){15}}
{\small
\put(6,-9){$p'_1$}
\put(-14,-7){$p'_2$}
\put(22,-11){$\cdots -2$}
\put(22,-6){$\cdots -1$}
\put(22,-1){$\cdots \ \phantom{ - }0$}
\put(22,4){$\cdots  \ \phantom{ - }1$}
\put(22,9){$\cdots  \ \phantom{ - }2$}
\put(22,13){\quad \ height}
\put(-11.2,4.2){$\times$}
\put(5,5){\circle*{1}}
\put(-10,6){$w_1^*$}
\put(5,6){$w_2$}
\put(0,-5){\circle*{1}}
\put(13.8,-5.8){$\times$}
\put(0,-4){$w_1$}
\put(15,-4){$w_2^*$}
}
\end{picture}
}
\\
\vspace{1cm}
\caption{An example of two paths and the map
$r_y$ for $y=1$.}\label{fig:C-relation}
\end{figure}
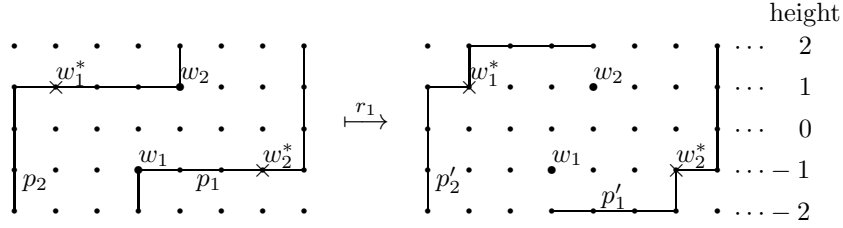

\begin{figure}[b]
{\setlength{\unitlength}{1mm}
\begin{picture}(45,18)(-15,3)
{\small
\multiput(-15,-5)(0,5){5}{\circle*{0.6}}
\multiput(-10,-5)(0,5){5}{\circle*{0.6}}
\multiput(-5,-5)(0,5){5}{\circle*{0.6}}
\multiput(0,-5)(0,5){5}{\circle*{0.6}}
\multiput(5,-5)(0,5){5}{\circle*{0.6}}
\multiput(10,-5)(0,5){5}{\circle*{0.6}}
\multiput(15,-5)(0,5){5}{\circle*{0.6}}
\multiput(20,-5)(0,5){5}{\circle*{0.6}}
\multiput(25,-5)(0,5){5}{\circle*{0.6}}
\multiput(30,-5)(0,5){5}{\circle*{0.6}}
%
\put(-10,-5){\line(0,1){9.7}}
\put(-10,4.7){\line(1,0){14}}
\put(6,4.7){\line(1,0){24}}
\put(30,4.7){\line(0,1){10.3}}
{
\put(5,-5){\line(0,1){10.3}}
\put(5,5.3){\line(1,0){10}}
\put(15,5.3){\line(0,1){9.7}}
}
\put(-1.2,9.2){$\times$}
\put(18.8,-0.8){$\times$}
\put(5,5){\circle*{1}}
\put(15,5){\circle*{1}}
\put(5,6){$u$}
\put(13,6){$v$}
\put(-10,-7){$u_2$}
\put(5,-7){$u_1$}
\put(15,16){$v_2$}
\put(30,16){$v_1$}
\put(-14,0){$p_2$}
\put(1,0){$p_1$}
}
\end{picture}
}
\ $\overset{r_0}{\longmapsto}$ \ 
{\setlength{\unitlength}{1mm}
\begin{picture}(60,18)(-15,3)
{\small
\multiput(-15,-5)(0,5){5}{\circle*{0.6}}
\multiput(-10,-5)(0,5){5}{\circle*{0.6}}
\multiput(-5,-5)(0,5){5}{\circle*{0.6}}
\multiput(0,-5)(0,5){5}{\circle*{0.6}}
\multiput(5,-5)(0,5){5}{\circle*{0.6}}
\multiput(10,-5)(0,5){5}{\circle*{0.6}}
\multiput(15,-5)(0,5){5}{\circle*{0.6}}
\multiput(20,-5)(0,5){5}{\circle*{0.6}}
\multiput(25,-5)(0,5){5}{\circle*{0.6}}
\multiput(30,-5)(0,5){5}{\circle*{0.6}}
%
\put(35,-5){$\cdots -2$}
\put(35,0){$\cdots -1$}
\put(35,5){$\cdots \quad 0$}
\put(35,10){$\cdots \quad 1$}
\put(35,15){$\cdots \quad 2$}
\put(40,20){height}
\put(-10,-5){\line(0,1){10}}
\put(-10,5){\line(1,0){10}}
\put(0,5){\line(0,1){5}}
\put(0,10){\line(1,0){15}}
\put(15,10){\line(0,1){5}}
\put(5,-5){\line(0,1){5}}
\put(5,0){\line(1,0){15}}
\put(20,0){\line(0,1){5}}
\put(20,5){\line(1,0){10}}
\put(30,5){\line(0,1){10}}
\put(-1.2,9.2){$\times$}
\put(18.8,-0.8){$\times$}
\put(5,5){\circle*{1}}
\put(15,5){\circle*{1}}
\put(5,6){$u$}
\put(13,6){$v$}
\put(-10, -7){$u_2$}
\put(5,-7){$u_1$}
\put(15,16){$v_2$}
\put(30,16){$v_1$}
\put(-14,0){$p'_2$}
\put(6,1){$p'_1$}
}
\end{picture}
}
\vspace{0.9cm}
\caption{An example of specially intersecting, transposed paths and
the map $r_0$}
\label{fig:r}
\end{figure}
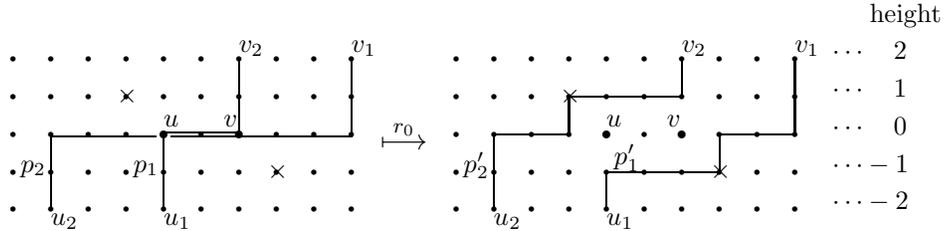

Let $y = 1, \dots , n-1$. For any $p_1, p_2 \in P(C_n)$, 
let $w_1$ (resp.\ $w_2$) be the leftmost point of height $-y$ on $p_1$
(resp.\ the rightmost point of height $y$ on $p_2$), 
i.e., if $w_1=(x_1,-y)$ and $w_2=(x_2,y)$, then 
$$x_1 = \min \{ x \mid  \text{$(x, -y)$ is on $p_1$}\},\qquad  
x_2 = \max \{ x \mid  \text{$(x, y)$ is on $p_2$}\}.$$
See Figure \ref{fig:C-relation} for example. 
Note that $w_1-(0,1)$ is on $p_1$ and $w_2+(0,1)$ is on $p_2$. 
We define the condition $(\bR_y)$ for any $p_1, p_2 \in P(C_n)$ 
as follows: 
\vspace{5pt}
\newline 
$(\bR_y)$ \  $w_1^*:=w_1+(-y-1, 2y)$ 
is on $p_2$ and $w_2^*:=w_2+(y+1, -2y)$ is on $p_1$. 
\vspace*{1pt}
\newline
For any $p_1, p_2 \in P(C_n)$ which satisfy 
$(\bR_y)$, we define $r_y(p_1, p_2)=(p'_1, p'_2)$ ($y=1, \dots , n-1$) 
as (see Figure \ref{fig:C-relation})
\begin{equation}\label{eq:ry}
\begin{aligned}
p'_1 & :u_1\overset{p_1}{\longrightarrow}w_1-(0,1)
\longrightarrow w_2^*-(0,1)
\longrightarrow w_2^*
\overset{p_1}{\longrightarrow}v_1,\\
p'_2 & :u_2\overset{p_2}{\longrightarrow}w_1^*
\longrightarrow w_1^*+(0,1)
\longrightarrow w_2+(0,1)\overset{p_2}{\longrightarrow}v_2. 
\end{aligned}
\end{equation}
For the $y=0$ case, 
let $p_1, p_2\in P(C_n)$ satisfy the following condition: \vspace{5pt}
\newline
$(\bR_0)$ \quad 
$(p_1,p_2)$ is specially intersecting at height 0. \vspace{5pt}
\newline
Then we define $r_0(p_1,p_2)=(p'_1,p'_2)$ as follows: 
If $p_1$ and $p_2$ are not transposed, then 
let $w_1$ (resp.\ $w_2$) be the leftmost point of height 0 on $p_1$
(resp.\ the rightmost point of height 0 on $p_2$), and set 
$w_1^*$ and $w_2^*$ as in $(\bR_y)$ by putting $y=0$. 
Then set $(p'_1,p'_2)$ as in \eqref{eq:ry}.
If $(p_1, p_2)$ is transposed (see Figure \ref{fig:r} for example), then
let $u$ (resp.\ $v$) be the leftmost 
(resp.\ rightmost) intersecting point of $p_1$ and $p_2$ at height 0. 
We assume that 
$u-(0,1)$ and $v+(0,1)$ is on $p_1$ while
$u-(1,0)$ and $v+(1, 0)$ is on $p_2$.
Set $r_0(p_1, p_2)=(p'_1,p'_2)$ by 
\begin{align*}
& p'_1:u_1\overset{p_1}{\longrightarrow}u-(0,1) 
\longrightarrow v+(1,-1) 
\longrightarrow v+(1,0) 
\overset{p_2}{\longrightarrow}v_1,\\
& p'_2:u_2\overset{p_2}{\longrightarrow}u-(1,0) 
\longrightarrow u+(-1,1) 
\longrightarrow v+(0,1)
\overset{p_1}{\longrightarrow}v_2.
\end{align*}
(Roughly speaking, $r_0$ ``resolves'' the transposed pair $(p_1,p_2)$.)
By \eqref{gen} and the definition of the $h$-label of type $C_n$, we have 
\begin{lem}\label{lem:weight}
$r_y$ ($y=0, \dots, n-1$) preserves the weight of $(p_1, p_2)$.
\end{lem}

Let $1\le i<j\le l$. 
For any $\bp=(p_1, \dots , p_l)$ such that the pair $(p_i,p_j)$
satisfies $(\bR_y)$ ($0 \le y \le n-1$), 
we define $r^{ij}_y(\bp)=(p'_1, \dots ,p'_l)$ by 
\begin{equation}\label{eq:r-map}
(p'_i, p'_j) := r_y(p_i, p_j), \qquad  p'_k:=p_k, \quad  (k\ne i,j).
\end{equation}
By Lemma \ref{lem:weight}, it is obvious that 
\begin{prop}\label{prop:r-y}
$r^{ij}_y$ preserves the weight for any  
$0 \le y \le n-1$ and $1\le i<j \le l$.
\end{prop}
Remark that $r^{ij}_y(\bp)$ for $\bp \in P(C_n; \mu, \lambda)$ 
may include an ordinarily intersecting pair of paths, which implies that
$r^{ij}_y(\bp)$ is not necessarily an element of $P(C_n; \mu, \lambda)$.

\subsection{The maps in the proof of Theorem \ref{thm:3-row}}\label{sec:3-row}
In this subsection, 
$\lambda/\mu$ is a skew diagram of $l(\lambda)= 3$. 
In this case, we have 
$P(C_n; \mu, \lambda)=
P_0(C_n; \mu, \lambda)\sqcup P_1(C_n; \mu, \lambda)\sqcup P_2(C_n; \mu, \lambda)$.
We define some maps
which we use in Section \ref{sec:bijective-T3} and show their properties.

{\it The map $g$}.
For any $\bp = (p_1,p_2,p_3) \in P_2(C_n; \mu, \lambda)$, 
let $u$ be the leftmost intersecting point of
$p_1$ and $p_3$, and let 
$v$ be the rightmost intersecting point of
$p_2$ and $p_3$ (see Figure \ref{fig:C-PP}).
Then set $u':=u+(-1,1)$ and $v':=v+(1,-1)$.
Let $P_2(C_n; \mu, \lambda)^{\times}$ be the set of all 
$\bp \in P_2(C_n; \mu, \lambda)$ such that both $u'$ and $v'$ are on some $p_i$
(actually, $u'$ is on $p_2$ and $v'$ is on $p_1$).
For example, $\bp$ in Figure \ref{fig:C-PP}
is an element of $P_2(C_n; \mu, \lambda)^{\times}$. 
Let $\tilde{P}(C_n; \mu, \lambda)$ be the subset of $\fP(C_n; \bu_{\mu}, \bv_{\lambda})$ 
defined in Section \ref{subsec:C-Tab}.
We define a map
\begin{equation*}
g:P_2(C_n; \mu, \lambda)^{\times} \to \tilde{P}(C_n; \mu, \lambda)
\end{equation*}
as follows: 
For $\bp=(p_1, p_2, p_3) \in P_2(C_n; \mu, \lambda)^{\times}$,
set $g(\bp) = (p'_1,p'_2,p'_3)$ as (see Figure \ref{fig:C-PP}) 
\begin{align*}
p'_1 & :u_1 \overset{p_1}{\longrightarrow}v'
\longrightarrow v'+(0,1) 
\overset{p_3}{\longrightarrow}v_1,\\
p'_2 & :u_2 \overset{p_2}{\longrightarrow}v
\longrightarrow u\overset{p_1}{\longrightarrow} v_2,\\
p'_3 & :u_3 \overset{p_3}{\longrightarrow}u'-(0,1)
\longrightarrow u'
\overset{p_2}{\longrightarrow} v_3.
\end{align*}
Then we easily show that 
\begin{lem}\label{lem:g-map}
\begin{enumerate}
\item \label{item:g-map}
$g$ is a weight-preserving, sign-preserving injection 
and $\tilde{P}(C_n; \mu, \lambda)=\Image \, g \sqcup P_0(C_n; \mu, \lambda)$.
\item \label{item:g-Tab}
$\{ T(\bp) \mid  \bp \in \Image \, g \} = 
\{ T\in \widetilde{\Tab}(C_n; \lambda/\mu) \mid  
\text{$T$ contains 
\begin{picture}(24,20)(0,16)
\multiput(0,0)(12,0){3}{\line(0,1){36}}
\multiput(0,0)(0,12){4}{\line(1,0){24}}
\multiput(3,3)(0,12){3}{$\overline{n}$}
\multiput(15,3)(0,12){3}{$n$}
\end{picture}
} 
\vspace{0.2cm}
\}$.
\end{enumerate}
\end{lem}
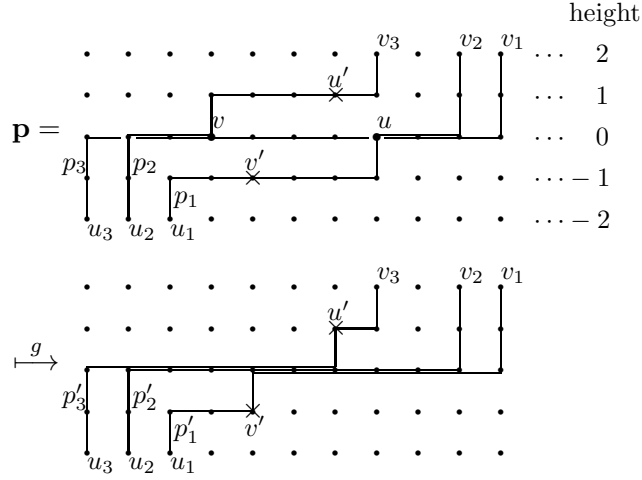
\begin{figure}[t]
{\setlength{\unitlength}{1.1mm}
\begin{tabular}{cc}
$\bp =$
&
\begin{picture}(55,18)(-15,0)
\multiput(-15,-10)(5,0){11}{\circle*{0.5}}
\multiput(-15,-5)(5,0){11}{\circle*{0.5}}
\multiput(-15,0)(5,0){11}{\circle*{0.5}}
\multiput(-15,5)(5,0){11}{\circle*{0.5}}
\multiput(-15,10)(5,0){11}{\circle*{0.5}}
{\small
\put(39,-11){$\cdots -2$}
\put(39,-6){$\cdots -1$}
\put(39,-1){$\cdots \quad 0$}
\put(39,4){$\cdots  \quad 1$}
\put(39,9){$\cdots  \quad 2$}
\put(39,14){\quad \ height}
\put(-5,-12){$u_1$}
\put(-10,-12){$u_2$}
\put(-15,-12){$u_3$}
\put(35,11){$v_1$}
\put(30,11){$v_2$}
\put(20,11){$v_3$}
\put(-4.5,-8){$p_1$}
\put(-9.5,-4){$p_2$}
\put(-18,-4){$p_3$}
}
%
\put(0,0){\circle*{1}}
\put(20,0){\circle*{1}}
{\small
\put(0,1){$v$}
\put(20,1){$u$}
}
\put(3.7,-5.8){$\times$}
\put(13.7,4.2){$\times$}
{\small
\put(4,-4){$v'$}
\put(14,6){$u'$}
}
\put(-15,-10){\line(0,1){9.8}}
\put(-15,-0.2){\line(1,0){4}}
\put(-9,-0.2){\line(1,0){28}}
\put(21,-0.2){\line(1,0){14}}
\put(35,-0.2){\line(0,1){10.2}}
\put(-10,-10){\line(0,1){10.2}}
\put(-10,0.2){\line(1,0){10}}
\put(0,0.2){\line(0,1){4.8}}
\put(0,5){\line(1,0){20}}
\put(20,5){\line(0,1){5}}
\put(-5,-10){\line(0,1){5}}
\put(-5,-5){\line(1,0){25}}
\put(20,-5){\line(0,1){5.2}}
\put(20,0.2){\line(1,0){10}}
\put(30,0.2){\line(0,1){9.8}}
\end{picture}
\\
$\overset{g}{\longmapsto}$ 
&
\begin{picture}(55,27)(-15,0)
\multiput(-15,-10)(5,0){11}{\circle*{0.5}}
\multiput(-15,-5)(5,0){11}{\circle*{0.5}}
\multiput(-15,0)(5,0){11}{\circle*{0.5}}
\multiput(-15,5)(5,0){11}{\circle*{0.5}}
\multiput(-15,10)(5,0){11}{\circle*{0.5}}
{\small
%
\put(-5,-12){$u_1$}
\put(-10,-12){$u_2$}
\put(-15,-12){$u_3$}
\put(35,11){$v_1$}
\put(30,11){$v_2$}
\put(20,11){$v_3$}
\put(-4.5,-8){$p'_1$}
\put(-9.5,-4){$p'_2$}
\put(-18,-4){$p'_3$}
}
\put(3.7,-5.8){$\times$}
\put(13.7,4.2){$\times$}
{\small
\put(4,-8){$v'$}
\put(14,6){$u'$}
}
\put(-15,-10){\line(0,1){10.4}}
\put(-15,0.4){\line(1,0){30}}
\put(15,0.4){\line(0,1){4.6}}
\put(15,5){\line(1,0){5}}
\put(20,5){\line(0,1){5}}
\put(-10,-10){\line(0,1){10}}
\put(-10,0){\line(1,0){40}}
\put(30,0){\line(0,1){10}}
\put(-5,-10){\line(0,1){5}}
\put(-5,-5){\line(1,0){10}}
\put(5,-5){\line(0,1){4.6}}
\put(5,-0.4){\line(1,0){30}}
\put(35,-0.4){\line(0,1){10.4}}
\end{picture}
\end{tabular}
}
\\
\vspace{1cm}
\caption{An example of $\bp \in P_2(C_n; \mu, \lambda)^{\times}$ and the map $g$.}
\label{fig:C-PP}
\end{figure}

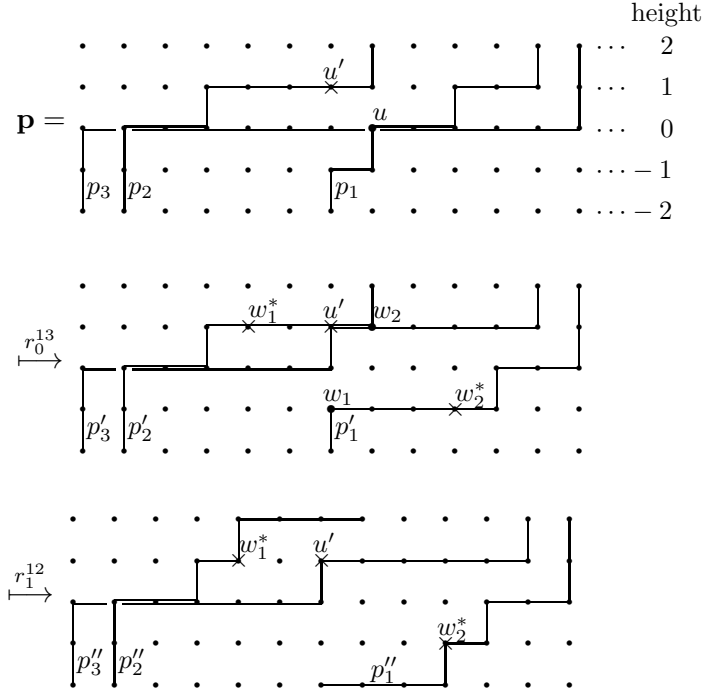
\begin{figure}
{\setlength{\unitlength}{1.1mm}
$\bp =$\hspace{0.1cm}
\begin{picture}(60,17)(-30,0)
\multiput(-30,-10)(5,0){13}{\circle*{0.5}}
\multiput(-30,-5)(5,0){13}{\circle*{0.5}}
\multiput(-30,0)(5,0){13}{\circle*{0.5}}
\multiput(-30,5)(5,0){13}{\circle*{0.5}}
\multiput(-30,10)(5,0){13}{\circle*{0.5}}
{\small
\put(32,-11){$\cdots -2$}
\put(32,-6){$\cdots -1$}
\put(32,-1){$\cdots \ \phantom{ - }0$}
\put(32,4){$\cdots  \ \phantom{ - }1$}
\put(32,9){$\cdots  \ \phantom{ - }2$}
\put(32,13){\quad \ height}
%
%
}
\put(-25,-10){\line(0,1){10.2}}
\put(-25,0.2){\line(1,0){10}}
\put(-15,0.2){\line(0,1){4.8}}
\put(-15,5){\line(1,0){20}}
\put(5,5){\line(0,1){5}}
\put(-30,-10){\line(0,1){9.8}}
\put(-30,-0.2){\line(1,0){4}}
\put(-24,-0.2){\line(1,0){28}}
\put(6,-0.2){\line(1,0){24}}
\put(30,-0.2){\line(0,1){10.2}}
\put(0,-10){\line(0,1){5}}
\put(0,-5){\line(1,0){5}}
\put(5,-5){\line(0,1){5.2}}
\put(5,0.2){\line(1,0){10}}
\put(15,0.2){\line(0,1){4.8}}
\put(15,5){\line(1,0){10}}
\put(25,5){\line(0,1){5}}
{\small
\put(5,0){\circle*{1}}
\put(5,1){$u$}
%
%
%
%
\put(-29.5,-8){$p_3$}
\put(-24.5,-8){$p_2$}
\put(0.5,-8){$p_1$}
\put(-1.2,4.2){$\times$}
\put(-1,6){$u'$}
}
\end{picture}
\\ 
$\overset{r_0^{13}}{\longmapsto}$ \
\begin{picture}(60,28)(-30,0)
\multiput(-30,-10)(5,0){13}{\circle*{0.5}}
\multiput(-30,-5)(5,0){13}{\circle*{0.5}}
\multiput(-30,0)(5,0){13}{\circle*{0.5}}
\multiput(-30,5)(5,0){13}{\circle*{0.5}}
\multiput(-30,10)(5,0){13}{\circle*{0.5}}
\put(-25,-10){\line(0,1){10.2}}
\put(-25,0.2){\line(1,0){10}}
\put(-15,0.2){\line(0,1){5}}
\put(-15,5.2){\line(1,0){20}}
\put(5,5.2){\line(0,1){4.8}}
\put(-30,-10){\line(0,1){9.8}}
\put(-30,-0.2){\line(1,0){4}}
\put(-24,-0.2){\line(1,0){24}}
\put(0,-0.2){\line(0,1){5}}
\put(0,4.8){\line(1,0){25}}
\put(25,4.8){\line(0,1){5.2}}
\put(0,-10){\line(0,1){5}}
\put(0,-5){\line(1,0){20}}
\put(20,-5){\line(0,1){5}}
\put(20,0){\line(1,0){10}}
\put(30,0){\line(0,1){10}}
{\small
%
%
\put(-11.2,4.2){$\times$}
\put(5,5){\circle*{1}}
\put(-10,6.2){$w_1^*$}
\put(5,6){$w_2$}
\put(0,-5){\circle*{1}}
\put(13.8,-5.8){$\times$}
\put(-1,-4){$w_1$}
\put(15,-4){$w_2^*$}
\put(-29.5,-8){$p'_3$}
\put(-24.5,-8){$p'_2$}
\put(0.5,-8){$p'_1$}
\put(-1.2,4.2){$\times$}
\put(-1,6){$u'$}
}
\end{picture}
\\ 
$\overset{r_1^{12}}{\longmapsto}$ \
\begin{picture}(60,28)(-30,0)
\multiput(-30,-10)(5,0){13}{\circle*{0.5}}
\multiput(-30,-5)(5,0){13}{\circle*{0.5}}
\multiput(-30,0)(5,0){13}{\circle*{0.5}}
\multiput(-30,5)(5,0){13}{\circle*{0.5}}
\multiput(-30,10)(5,0){13}{\circle*{0.5}}
\put(-25,-10){\line(0,1){10.2}}
\put(-25,0.2){\line(1,0){10}}
\put(-15,0.2){\line(0,1){4.8}}
\put(-15,5){\line(1,0){5}}
\put(-10,5){\line(0,1){5}}
\put(-10,10){\line(1,0){15}}
\put(-30,-10){\line(0,1){9.8}}
\put(-30,-0.2){\line(1,0){4}}
\put(-24,-0.2){\line(1,0){24}}
\put(0,-0.2){\line(0,1){5.2}}
\put(0,5){\line(1,0){25}}
\put(25,5){\line(0,1){5}}
\put(0,-10){\line(1,0){15}}
\put(15,-10){\line(0,1){5}}
\put(15,-5){\line(1,0){5}}
\put(20,-5){\line(0,1){5}}
\put(20,0){\line(1,0){10}}
\put(30,0){\line(0,1){10}}
{\small
%
%
\put(-11.2,4.2){$\times$}
\put(-10,6){$w_1^*$}
\put(13.8,-5.8){$\times$}
\put(14,-4){$w_2^*$}
\put(-29.5,-8){$p''_3$}
\put(-24.5,-8){$p''_2$}
\put(6,-9){$p''_1$}
\put(-1.2,4.2){$\times$}
\put(-1,6){$u'$}
}
\end{picture}
}
\vspace{1cm}
\caption{An example of $\bp \in P_2(C_n; \mu, \lambda)^{\circ}$ 
and the map $f_2^{13}$ of Case ($f_2^{13}$-b).}
\label{fig:r-paths}
\end{figure}

{\it The maps $f_2^{13}$ and $f_2^{23}$}. 
Let $P^{ij}_1(C_n; \mu, \lambda)$ ($1\le i<j\le 3$) be the set of all
$\bp = (p_1, p_2, p_3) \in P_1(C_n; \mu, \lambda)$ such that $(p_i, p_j)$ is transposed.
Then we have 
$P_1(C_n; \mu, \lambda)= 
P^{12}_1(C_n; \mu, \lambda) \sqcup P^{23}_1(C_n; \mu, \lambda)$. 
Let 
$P_2(C_n; \mu, \lambda)^{\circ} := 
P_2(C_n; \mu, \lambda)\setminus P_2(C_n; \mu, \lambda)^{\times}$. 
We define a map
\begin{align*}
f^{13}_2:P_2(C_n; \mu, \lambda)^{\circ}
\to P^{23}_1(C_n; \mu, \lambda)
\end{align*}  
as follows,  
using the weight-preserving maps $r_y^{ij}$ defined in \eqref{eq:r-map}
(roughly speaking, 
$f^{13}_2$ resolves the transposed pair 
$(p_1,p_3)$ of $\bp = (p_1,p_2, p_3)$ without producing 
any ordinarily intersecting paths):
For any $\bp =(p_1,p_2,p_3)\in 
P_2(C_n; \mu, \lambda)^{\circ}$, 
let $u$ be the leftmost intersecting point of $p_1$ and $p_3$
and $u':= u+(-1,1)$. Set $\bp':=(p'_1,p'_2,p'_3)=r_0^{13}(\bp)$, 
which is well-defined. 
Then 
\vspace{5pt}
\newline
{\bf Case} ($f_2^{13}$-a) \  
If $u'$ is not on any $p_i$, set $f_2^{13}(\bp)=\bp'$, which is in
$P_1^{23}(C_n; \mu, \lambda)$. 
(Otherwise, $u'$ is on $p_2$ and  
$(p'_2,p'_3)$ is ordinarily intersecting.)
\newline
{\bf Case} ($f_2^{13}$-b) \
Otherwise, 
$(p'_1,p'_2)$ satisfies the condition $(\bR_1)$ in Appendix \ref{sec:r}. 
Set $f_2^{13}(\bp)=r^{12}_1(\bp')$, 
which is in $P_1^{23}(C_n; \mu, \lambda)$ (see Figure \ref{fig:r-paths}).
\vspace{5pt}

We remark that if $\bp \in P_2(C_n; \mu, \lambda)^{\times}$, then
$(p'_2,p'_3)$ is ordinarily intersecting, but 
the procedure of ($f_2^{13}$-b) is not well-defined, 
for $(p'_1,p'_2)$ 
does not satisfy $(\bR_1)$. 

We also define 
\begin{align*}
f^{23}_2:P_2(C_n; \mu, \lambda)^{\circ}
\to P^{12}_1(C_n; \mu, \lambda)
\end{align*}
by $\omega \circ f^{13}_2 \circ \omega$, where 
\begin{equation}\label{eq:omega}
\omega : \fP(\bu_{\mu}, \bv_{\lambda}) \to 
\fP(\bu_{\tilde{\mu}}, \bv_{\tilde{\lambda}})
\end{equation}
is a map that rotates $\bp$ by 180 degrees 
around a fixed point $(x,0)$ 
such that $2x -\lambda_1+l(\lambda)-1 \in \bZ_{\ge 0}$  
(so that $\tilde{\lambda}$ and $\tilde{\mu}$
are partitions). 
Then, $f_2^{23}$ resolves the transposed pair $(p_2,p_3)$
of $\bp = (p_1,p_2,p_3)\in P_2(C_n; \mu, \lambda)^{\circ}$.


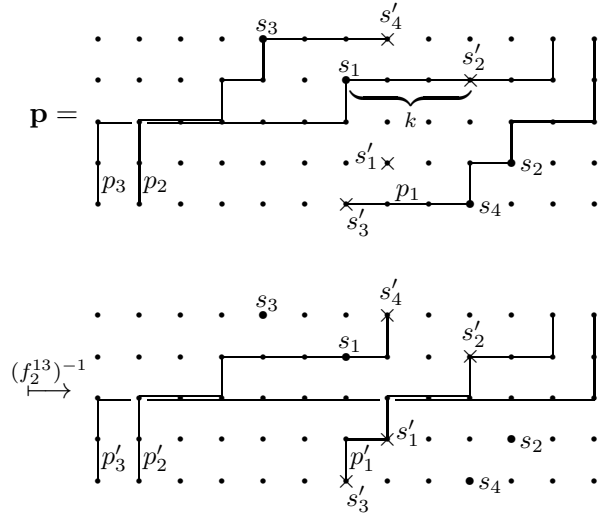
\begin{figure}
{\setlength{\unitlength}{1.1mm}
$\bp =$ \ 
\begin{picture}(60,18)(-30,0)
\multiput(-30,-10)(5,0){13}{\circle*{0.5}}
\multiput(-30,-5)(5,0){13}{\circle*{0.5}}
\multiput(-30,0)(5,0){13}{\circle*{0.5}}
\multiput(-30,5)(5,0){13}{\circle*{0.5}}
\multiput(-30,10)(5,0){13}{\circle*{0.5}}
\put(-25,-10){\line(0,1){10.2}}
\put(-25,0.2){\line(1,0){10}}
\put(-15,0.2){\line(0,1){4.8}}
\put(-15,5){\line(1,0){5}}
\put(-10,5){\line(0,1){5}}
\put(-10,10){\line(1,0){15}}
\put(-30,-10){\line(0,1){9.8}}
\put(-30,-0.2){\line(1,0){4}}
\put(-24,-0.2){\line(1,0){24}}
\put(0,-0.2){\line(0,1){5.2}}
\put(0,5){\line(1,0){25}}
\put(25,5){\line(0,1){5}}
\put(0,-10){\line(1,0){15}}
\put(15,-10){\line(0,1){5}}
\put(15,-5){\line(1,0){5}}
\put(20,-5){\line(0,1){5}}
\put(20,0){\line(1,0){10}}
\put(30,0){\line(0,1){10}}
{\small
%
%
%
%
\put(-29.5,-8){$p_3$}
\put(-24.5,-8){$p_2$}
\put(6,-9){$p_1$}
\put(0,5){\circle*{1}}
\put(-1,6){$s_1$}
\put(-10,10){\circle*{1}}
\put(15,-10){\circle*{1}}
\put(20,-5){\circle*{1}}
\put(-11,11){$s_3$}
\put(16,-11){$s_4$}
\put(21,-6){$s_2$}
\put(3.8,9.2){$\times$}
\put(13.8,4.2){$\times$}
\put(3.8,-5.8){$\times$}
\put(-1.2,-10.8){$\times$}
\put(4,12){$s'_4$}
\put(14,7){$s'_2$}
\put(1,-5){$s'_1$}
\put(0,-13){$s'_3$}
}
\put(0.5,4){$\underbrace{\hspace{1.6cm}}_k$}
\end{picture}
\\
\vspace{1.7cm}
$\overset{(f^{13}_2)^{-1}}{\longmapsto}$
\begin{picture}(60,17)(-30,0)
\multiput(-30,-10)(5,0){13}{\circle*{0.5}}
\multiput(-30,-5)(5,0){13}{\circle*{0.5}}
\multiput(-30,0)(5,0){13}{\circle*{0.5}}
\multiput(-30,5)(5,0){13}{\circle*{0.5}}
\multiput(-30,10)(5,0){13}{\circle*{0.5}}
{\small
%
%
}
\put(-25,-10){\line(0,1){10.2}}
\put(-25,0.2){\line(1,0){10}}
\put(-15,0.2){\line(0,1){4.8}}
\put(-15,5){\line(1,0){20}}
\put(5,5){\line(0,1){5}}
\put(-30,-10){\line(0,1){9.8}}
\put(-30,-0.2){\line(1,0){4}}
\put(-24,-0.2){\line(1,0){28}}
\put(6,-0.2){\line(1,0){24}}
\put(30,-0.2){\line(0,1){10.2}}
\put(0,-10){\line(0,1){5}}
\put(0,-5){\line(1,0){5}}
\put(5,-5){\line(0,1){5.2}}
\put(5,0.2){\line(1,0){10}}
\put(15,0.2){\line(0,1){4.8}}
\put(15,5){\line(1,0){10}}
\put(25,5){\line(0,1){5}}
{\small
%
%
%
%
\put(-29.5,-8){$p'_3$}
\put(-24.5,-8){$p'_2$}
\put(0.5,-8){$p'_1$}
\put(0,5){\circle*{1}}
\put(-1,6){$s_1$}
\put(-10,10){\circle*{1}}
\put(15,-10){\circle*{1}}
\put(20,-5){\circle*{1}}
\put(-11,11){$s_3$}
\put(16,-11){$s_4$}
\put(21,-6){$s_2$}
\put(3.8,9.2){$\times$}
\put(13.8,4.2){$\times$}
\put(3.8,-5.8){$\times$}
\put(-1.2,-10.8){$\times$}
\put(4,12){$s'_4$}
\put(14,7){$s'_2$}
\put(6,-5){$s'_1$}
\put(0,-13){$s'_3$}
}
\end{picture}
}
\vspace{1.5cm}
\caption{An example of $\bp \in P_1^{23}(C_n; \mu, \lambda)$ which satisfy 
($\bFFab$), and the inverse procedure of ($f_2^{13}$-b)}\label{fig:r-paths-inverse}
\end{figure}

Next, 
we give the conditions to describe the images of 
$f_2^{13}$ and $f_2^{23}$.  
For any $\bp \in P_1^{23}(C_n; \mu, \lambda)$
(see Figure \ref{fig:r-paths-inverse}), 
let $s_1$ be the leftmost point of height $1$ on $p_3$, 
$s_2$ be the rightmost point of height $-1$ on $p_1$,
$s_3$ be the leftmost point of height $2$ on $p_2$, and 
$s_4$ be the rightmost point of height $-2$ on $p_1$.
For each $i=1,\dots,4$, set $s'_i := s_i+(y,-2y)$, where $y$ is the height of $s_i$. 
If $s'_2$ is on $p_3$, then we define $k$ as the number of steps of $p_3$
between $s_1$ and $s'_2$. 
Then define conditions ($\bFFaa$) and ($\bFFab$) 
for $\bp\in P_1^{23}(C_n; \mu, \lambda)$ as follows:
\vspace{5pt}
\newline
\begin{tabular}{ll}
($\bFFaa$) & \ $\bp$ 
satisfies all of the following conditions:
\end{tabular}
\begin{itemize}
\item $s'_1$ is on $p_1$.
\item $s'_2$ is on $p_3$ and $k$ is odd.
\end{itemize}
\begin{tabular}{ll}
($\bFFab$) & \ $\bp$ 
satisfies all of the following conditions:
\end{tabular}
\begin{itemize}
\item $s'_1$ is not on $p_1$.
\item $s'_2$ is on $p_3$ and $k$ is odd.
\item $s'_3$ is on $p_1$.
\item $s'_4$ is on $p_2$.
\end{itemize}
\vspace{5pt}
We also define  conditions  
($\bFFba$) and ($\bFFbb$) for $\bp\in P_1^{12}(C_n; \mu, \lambda)$ 
as follows: 
\vspace{5pt}
\newline
\begin{tabular}{ll}
($\bFFba$) & \ $\omega (\bp)$
satisfies the condition ($\bFFaa$). \\
($\bFFbb$) & \ $\omega (\bp)$
satisfies the condition ($\bFFab$).
\end{tabular}
\vspace{5pt}
\newline
For example, $\bp\in P_1^{23}(C_n;\mu, \lambda)$ 
in Figure \ref{fig:r-paths-inverse}
satisfies the condition ($\bFFab$). 
Note that ($\bFFaa$) and ($\bFFab$) 
(resp.\ ($\bFFba$) and ($\bFFbb$)) are exclusive with each other.

We have
\begin{lem}\label{lem:f2-map}
\begin{enumerate}
\item \label{item:condition}
For $\bp \in P_1(C_n; \mu, \lambda)$, 
\begin{enumerate}
\item
$\bp \in \Image \, f^{13}_2$ if and only if either
{\rm ($\bFFaa$)} or {\rm ($\bFFab$)} is satisfied.
\item 
$\bp \in\Image \, f_2^{23}$  if and only if 
either {\rm ($\bFFba$)} or {\rm ($\bFFbb$)} is satisfied.
\end{enumerate}
\item 
$f^{13}_2$ and $f^{23}_2$ are weight-preserving, sign-inverting injections.
\end{enumerate}
\end{lem}

\begin{proof}
We prove it for $f^{13}_2$. 
We can check that 
$\bp$ in the image of ($f_2^{13}$-a) satisfy 
($\bFFaa$). 
Conversely, one can invert the procedure of ($f_2^{13}$-a) 
for any $\bp\in P_1^{23}(C_n; \mu, \lambda)$ that satisfies ($\bFFaa$). 
The same holds when ($f_2^{13}$-a) (resp.\ ($\bFFaa$)) are replaced
with ($f_2^{13}$-b) (resp.\ ($\bFFab$))
(see Figure \ref{fig:r-paths-inverse} for example).
\end{proof}

\begin{figure}
{\setlength{\unitlength}{1.1mm}
\hspace{0.54cm}
\begin{picture}(60,17)(-30,0)
\multiput(-25,-10)(5,0){12}{\circle*{0.5}}
\multiput(-25,-5)(5,0){12}{\circle*{0.5}}
\multiput(-25,0)(5,0){12}{\circle*{0.5}}
\multiput(-25,5)(5,0){12}{\circle*{0.5}}
\multiput(-25,10)(5,0){12}{\circle*{0.5}}
{\small
\put(32,-11){$\cdots -2$}
\put(32,-6){$\cdots -1$}
\put(32,-1){$\cdots \ \phantom{ - }0$}
\put(32,4){$\cdots  \ \phantom{ - }1$}
\put(32,9){$\cdots  \ \phantom{ - }2$}
\put(32,13){\quad \ height}
%
%
}
\put(-25,-10){\line(0,1){10.2}}
\put(-25,0.2){\line(1,0){10}}
\put(-15,0.2){\line(0,1){4.8}}
\put(-15,5){\line(1,0){20}}
\put(5,5){\line(0,1){5}}
\put(-20,-10){\line(0,1){9.8}}
\put(-20,-0.2){\line(1,0){24}}
\put(6,-0.2){\line(1,0){24}}
\put(30,-0.2){\line(0,1){10.2}}
\put(0,-10){\line(0,1){5}}
\put(0,-5){\line(1,0){5}}
\put(5,-5){\line(0,1){5.2}}
\put(5,0.2){\line(1,0){10}}
\put(15,0.2){\line(0,1){4.8}}
\put(15,5){\line(1,0){10}}
\put(25,5){\line(0,1){5}}
{\small
\put(5,0){\circle*{1}}
\put(5,1){$u$}
%
%
%
\put(0,-5){\circle*{1}}
\put(-1,-4){$w_1$}
\put(-11.2,4.2){$\times$}
\put(-11,6){$w^*_1$}
\put(-24.5,-8){$p_3$}
\put(-19.5,-8){$p_2$}
\put(0.5,-8){$p_1$}
\put(-1.2,4.2){$\times$}
\put(-1,6){$u'$}
}
\end{picture}
\\ 
$\overset{r_0^{12}}{\mapsto}$ \
\begin{picture}(60,28)(-30,0)
\multiput(-25,-10)(5,0){12}{\circle*{0.5}}
\multiput(-25,-5)(5,0){12}{\circle*{0.5}}
\multiput(-25,0)(5,0){12}{\circle*{0.5}}
\multiput(-25,5)(5,0){12}{\circle*{0.5}}
\multiput(-25,10)(5,0){12}{\circle*{0.5}}
\put(-25,-10){\line(0,1){10.2}}
\put(-25,0.2){\line(1,0){10}}
\put(-15,0.2){\line(0,1){5}}
\put(-15,5.2){\line(1,0){20}}
\put(5,5.2){\line(0,1){4.8}}
\put(-20,-10){\line(0,1){9.8}}
\put(-20,-0.2){\line(1,0){20}}
\put(0,-0.2){\line(0,1){5}}
\put(0,4.8){\line(1,0){25}}
\put(25,4.8){\line(0,1){5.2}}
\put(0,-10){\line(0,1){5}}
\put(0,-5){\line(1,0){20}}
\put(20,-5){\line(0,1){5}}
\put(20,0){\line(1,0){10}}
\put(30,0){\line(0,1){10}}
{\small
%
%
\put(-11.2,4.2){$\times$}
\put(5,5){\circle*{1}}
\put(-10,6.2){$w_1^*$}
\put(5,6){$w_2$}
\put(0,-5){\circle*{1}}
\put(13.8,-5.8){$\times$}
\put(-1,-4){$w_1$}
\put(15,-4){$w_2^*$}
\put(-24.5,-8){$p'_3$}
\put(-19.5,-8){$p'_2$}
\put(0.5,-8){$p'_1$}
\put(-1.2,4.2){$\times$}
\put(-1,6){$u'$}
}
\end{picture}
\\ 
\hspace{0.27cm}
$\overset{r_1^{13}}{\mapsto}$ 
\begin{picture}(60,28)(-30,0)
\multiput(-25,-10)(5,0){12}{\circle*{0.5}}
\multiput(-25,-5)(5,0){12}{\circle*{0.5}}
\multiput(-25,0)(5,0){12}{\circle*{0.5}}
\multiput(-25,5)(5,0){12}{\circle*{0.5}}
\multiput(-25,10)(5,0){12}{\circle*{0.5}}
\put(-25,-10){\line(0,1){10.2}}
\put(-25,0.2){\line(1,0){10}}
\put(-15,0.2){\line(0,1){4.8}}
\put(-15,5){\line(1,0){5}}
\put(-10,5){\line(0,1){5}}
\put(-10,10){\line(1,0){15}}
\put(-20,-10){\line(0,1){9.8}}
\put(-20,-0.2){\line(1,0){20}}
\put(0,-0.2){\line(0,1){5.2}}
\put(0,5){\line(1,0){25}}
\put(25,5){\line(0,1){5}}
\put(0,-10){\line(1,0){15}}
\put(15,-10){\line(0,1){5}}
\put(15,-5){\line(1,0){5}}
\put(20,-5){\line(0,1){5}}
\put(20,0){\line(1,0){10}}
\put(30,0){\line(0,1){10}}
{\small
%
%
\put(-11.2,4.2){$\times$}
\put(-10,6){$w_1^*$}
\put(13.8,-5.8){$\times$}
\put(14,-4){$w_2^*$}
\put(-24.5,-8){$p''_3$}
\put(-19.5,-8){$p''_2$}
\put(6,-9){$p''_1$}
\put(-1.2,4.2){$\times$}
\put(-1,6){$u'$}
}
\end{picture}
}
\vspace{1cm}
\caption{An example of $\bp \in P_1^{12}(C_n; \mu, \lambda)$ 
and the map $f_1^{12}$ of Case ($f^{12}_1$-b).}\label{fig:r1-path}
\end{figure}
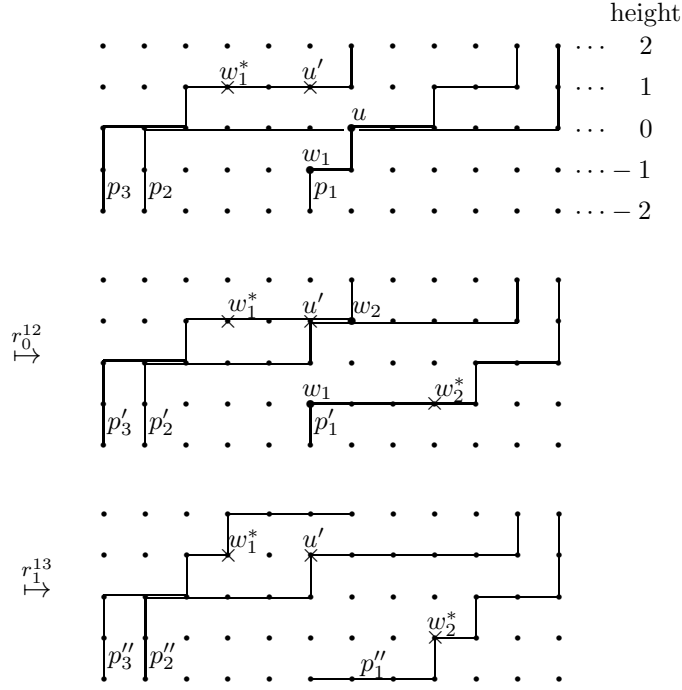

{\it The maps $f_1^{12}$ and $f_1^{23}$}. 
We define a map
$$f^{12}_1:P^{12}_1(C_n; \mu, \lambda) \to P_0(C_n; \mu, \lambda), $$
as follows, using the weight-preserving maps $r_y^{ij}$ defined in \ref{eq:r-map}
(roughly speaking, $f^{12}_1$ 
resolves the transposed pair $(p_1,p_2)$ of $\bp = (p_1,p_2,p_3)$
without producing any ordinarily intersecting paths):
For any $\bp=(p_1, p_2, p_3)\in P^{12}_1(C_n; \mu, \lambda)$
(see Figure \ref{fig:r1-path}), 
let $w_1$  be the leftmost point on $p_1$ at height $-1$ and
$w_1^*=w_1 + (-2,2)$.
Let $u$ be the leftmost intersecting point of $p_1$ and $p_2$
and $u':= u+(-1,1)$.  
Set $\bp' := (p'_1, p'_2, p'_3) = r_0^{12}(\bp)$, which is well-defined. 
Then
\vspace{5pt}
\newline
{\bf Case} ($f_1^{12}$-a) \ 
If $u'$ is not on any $p_i$, 
set $f^{12}_1(\bp)=\bp'$, which is in $P_0(C_n; \mu, \lambda)$. 
(Otherwise, $u'$ is on $p_3$ and  
$(p'_2,p'_3)$ is ordinarily intersecting. )
\newline
{\bf Case} ($f_1^{12}$-b) \ 
If $u'$ is on $p_3$  and  
$w_1^*$ is on $p_3$, then $(p'_1,p'_3)$ satisfies $(\bR_1)$. 
Set $f^{12}_1(\bp)=r^{13}_1(\bp')$, 
which is in $P_0(C_n; \mu, \lambda)$. 
(If $w_1^*$ is not on $p_3$, then $r^{13}_1(\bp')$ is not defined.)
\newline 
{\bf Case} ($f_1^{12}$-c) \ 
Otherwise, $(p_2, p_3)$ satisfies $(\bR_0)$. If we set  
$\bp''=r^{23}_0(\bp')=(p''_1,p''_2, p''_3)$, then  
$(p''_1,p''_2)$ is ordinarily intersecting. 
For $(p''_1,p''_3)$ satisfies $(\bR_1)$, 
we set $f^{12}_1(\bp)=r^{13}_1(\bp'')$, 
which is in $P_0(C_n; \mu, \lambda)$.
\vspace{5pt}

We also define 
$$f^{23}_1:P^{23}_1(C_n; \mu, \lambda) \to P_0(C_n; \mu, \lambda)$$
as $f^{23}_1:=\omega \circ f^{12}_1 \circ \omega$, 
where $\omega$ is the map defined in \eqref{eq:omega}.
Then $f_1^{23}$ resolves the transposed pair 
$(p_2,p_3)$ of $\bp =(p_1,p_2,p_3)$. 

\begin{figure}
{\setlength{\unitlength}{0.8mm}
\begin{tabular}{lcl}
\begin{picture}(61,0)(-46,-1)
\put(-46,10){($\bFaa$)}
\multiput(-15,-10)(5,0){7}{\circle*{0.5}}
\multiput(-15,-5)(5,0){7}{\circle*{0.5}}
\multiput(-15,0)(5,0){7}{\circle*{0.5}}
\multiput(-15,5)(5,0){7}{\circle*{0.5}}
\multiput(-15,10)(5,0){7}{\circle*{0.5}}
\multiput(-15,0)(2,0){3}{\line(1,0){1}}
\put(-10,0){\line(0,1){5}}
\put(-10,5){\line(1,0){15}}
\put(-5,-5){\line(1,0){15}}
\put(10,-5){\line(0,1){5}}
\multiput(10,0)(2,0){3}{\line(1,0){1}}
\put(-9.5,5){$\underbrace{\hspace{1.1cm}}_{k= \text{odd}}$}
\put(3.1,3.8){$\times$}
\put(-6.9,-6.2){$\times$}
\put(-10,5){\circle*{1}}
\put(10,-5){\circle*{1}}
{\small 
\put(-10,6){$s_1$}
\put(-5,-9){$s'_1$}
\put(4,6){$s'_2$}
\put(10,-8){$s_2$}
}
\put(-14,1){$b$}
\put(11,1){$a$}
\end{picture}
& 
$\overset{T}{\longmapsto}$ \ \   
& 
\begin{picture}(25,18)(-43,0)
\put(-4,-4){$\scriptstyle b$}
\put(0,5){$\overbrace{\hspace{1.1cm}}^{k=\text{odd}}$}
\put(-0.5,-5){\line(1,0){15}}
\put(-0.5,0){\line(1,0){15}}
\put(-0.5,5){\line(1,0){15}}
\put(-0.5,-5){\line(0,1){10}}
\put(14.5,-5){\line(0,1){10}}
\put(1,1){$\scriptstyle n$}
\put(6,1){$\scriptstyle n$}
\put(11,1){$\scriptstyle n$}
\put(1,-4){$\scriptstyle \overline{n}$}
\put(6,-4){$\scriptstyle \overline{n}$}
\put(11,-4){$\scriptstyle \overline{n}$}
\put(16,1){$\scriptstyle a$}
\put(-10,-10){$a \succeq \overline{n}$, $b \preceq n$}
\end{picture}
\\
\begin{picture}(61,37)(-46,-1)
\put(-46,20){($\bFabone$)}
\multiput(-30,-10)(5,0){10}{\circle*{0.5}}
\multiput(-30,-5)(5,0){10}{\circle*{0.5}}
\multiput(-30,0)(5,0){10}{\circle*{0.5}}
\multiput(-30,5)(5,0){10}{\circle*{0.5}}
\multiput(-30,10)(5,0){10}{\circle*{0.5}}
\multiput(-30,5)(2,0){3}{\line(1,0){1}}
\put(-25,5){\line(0,1){5}}
\put(-25,10){\line(1,0){15}}
\put(-20,0){\line(1,0){10}}
\put(-10,0){\line(0,1){5}}
\put(-10,5){\line(1,0){15}}
\put(-9.5,5){$\underbrace{\hspace{1.1cm}}_{k=\text{odd}}$}
\put(-15,-10){\line(1,0){15}}
\put(0,-10){\line(0,1){5}}
\put(0,-5){\line(1,0){10}}
\put(10,-5){\line(0,1){5}}
\multiput(10,0)(2,0){3}{\line(1,0){1}}
\put(3.1,3.8){$\times$}
\put(-6.9,-6.2){$\times$}
\put(-16.9,-6.2){$\times$}
\put(-16.9,-11.2){$\times$}
\put(-11.9,8.8){$\times$}
\put(-10,5){\circle*{1}}
\put(10,-5){\circle*{1}}
\put(-25,10){\circle*{1}}
\put(0,-10){\circle*{1}}
{\small
\put(-14,5){$s_1$}
\put(-10,-6){$s'_1$}
\put(6,5){$s'_2$}
\put(10,-8){$s_2$}
\put(-30,11){$s_3$}
\put(-16,-14){$s'_3$}
\put(-21,-6){$s''_3$}
\put(0,-13){$s_4$}
\put(-10,11){$s'_4$}
}
\put(11,1){$a$}
\put(-29,1){$b$}
\end{picture}
& $\overset{T}{\longmapsto}$ \ \ 
& 
\begin{picture}(15,35)(-23,7.5)
\put(17,15){$\overbrace{\hspace{1.3cm}}^{k=k_4+k_5}$}
\put(0,15){$\overbrace{\hspace{1.3cm}}^{2k_3}$}
\put(0,0){\line(1,0){24}}
\put(0,5){\line(1,0){34}}
\put(0,10){\line(1,0){34}}
\put(0,15){\line(1,0){34}}
\put(0,0){\line(0,1){15}}
\put(16.5,5){\line(0,1){5}}
\put(24,0){\line(0,1){5}}
\put(24,10){\line(0,1){5}}
\put(34,5){\line(0,1){10}}
\put(-4,1){$\scriptstyle b$}
\put(1,1){$\scriptstyle \overline{n-1}$}
\put(9,1){$\scriptstyle \overline{n-1}$}
\put(17,1){$\scriptstyle \overline{n-1}$}
\put(1,11){$\scriptstyle n-1$}
\put(9,11){$\scriptstyle n-1$}
\put(17,11){$\scriptstyle n-1$}
\put(2.5,6){$\scriptstyle \overline{n}$}
\put(10.5,6){$\scriptstyle n$}
\put(19.5,6){$\scriptstyle \overline{n}$}
\put(25,6){$\scriptstyle \overline{n}$}
\put(30,6){$\scriptstyle \overline{n}$}
\put(25,11){$\scriptstyle n$}
\put(30,11){$\scriptstyle n$}
\put(35,11){$\scriptstyle a$}
\put(-10,-5){$a \succeq \overline{n}$, $b \preceq \overline{n}$}
\put(-10,-12){$k_4+k_5$: odd, $k_4 \ne 0$}
\end{picture}
\\
& 
%

& 

%
%
\begin{picture}(45,40)(-23,7.5)
\put(-28,7.5){\text{or}}
\put(-8,15){$\overbrace{\hspace{1.9cm}}^{2k_3+1}$}
\put(17,15){$\overbrace{\hspace{1.3cm}}^{k=k_4+k_5}$}
\put(-13,0){\line(1,0){37}}
\put(-13,5){\line(1,0){47}}
\put(-13,10){\line(1,0){47}}
\put(-8,15){\line(1,0){42}}
\put(-13,0){\line(0,1){10}}
\put(-8,0){\line(0,1){5}}
\put(-8,10){\line(0,1){5}}
\put(16.5,5){\line(0,1){5}}
\put(24,0){\line(0,1){5}}
\put(24,10){\line(0,1){5}}
\put(34,5){\line(0,1){10}}
\put(-7,1){$\scriptstyle \overline{n-1}$}
\put(1,1){$\scriptstyle \overline{n-1}$}
\put(9,1){$\scriptstyle \overline{n-1}$}
\put(17,1){$\scriptstyle \overline{n-1}$}
\put(-7,11){$\scriptstyle n-1$}
\put(1,11){$\scriptstyle n-1$}
\put(9,11){$\scriptstyle n-1$}
\put(17,11){$\scriptstyle n-1$}
\put(-11.5,1){$\scriptstyle \overline{n}$}
\put(-11.5,6){$\scriptstyle \overline{n}$}
\put(-4.5,6){$\scriptstyle n$}
\put(2.5,6){$\scriptstyle \overline{n}$}
\put(10.5,6){$\scriptstyle n$}
\put(19.5,6){$\scriptstyle \overline{n}$}
\put(25,6){$\scriptstyle \overline{n}$}
\put(30,6){$\scriptstyle \overline{n}$}
\put(25,11){$\scriptstyle n$}
\put(30,11){$\scriptstyle n$}
\put(35,11){$\scriptstyle a$}
\put(-10,-5){$a \succeq \overline{n}$}
\put(-10,-12){$k_4+k_5$: odd, $k_4 \ne 0$}
\end{picture}
\\
\begin{picture}(61,47)(-48,-1)
\put(-48,20){($\bFabtwo$)}
\multiput(-50,-10)(5,0){14}{\circle*{0.5}}
\multiput(-50,-5)(5,0){14}{\circle*{0.5}}
\multiput(-50,0)(5,0){14}{\circle*{0.5}}
\multiput(-50,5)(5,0){14}{\circle*{0.5}}
\multiput(-50,10)(5,0){14}{\circle*{0.5}}
\multiput(-50,0)(2,0){3}{\line(1,0){1}}
\put(-45,0){\line(0,1){5}}
\put(-45,5){\line(1,0){10}}
\put(-35,5){\line(0,1){5}}
\put(-35,10){\line(1,0){25}}
\put(-39.5,-5){$\overbrace{\hspace{1.5cm}}^{\text{even}}$}
\put(-40,-5){\line(1,0){20}}
\put(-20,-5){\line(0,1){5}}
\put(-20,0){\line(1,0){10}}
\put(-10,0){\line(0,1){5}}
\put(-10,5){\line(1,0){15}}
\put(-9.5,5){$\underbrace{\hspace{1.1cm}}_{k= \text{odd}}$}
\put(-25,-10){\line(1,0){25}}
\put(0,-10){\line(0,1){5}}
\put(0,-5){\line(1,0){10}}
\put(10,-5){\line(0,1){5}}
\multiput(10,0)(2,0){3}{\line(1,0){1}}
\put(3.1,3.8){$\times$}
\put(-6.9,-6.2){$\times$}
\put(-26.9,-6.2){$\times$}
\put(-26.9,-11.2){$\times$}
\put(-11.9,8.8){$\times$}
\put(-41.9,-6.2){$\times$}
\put(-10,5){\circle*{1}}
\put(10,-5){\circle*{1}}
\put(-35,10){\circle*{1}}
\put(-45,5){\circle*{1}}
\put(0,-10){\circle*{1}}
\put(-20,-5){\circle*{1}}
{\small
\put(-14,5){$s_1$}
\put(-10,-6){$s'_1$}
\put(6,5){$s'_2$}
\put(10,-8){$s_2$}
\put(-40,11){$s_3$}
\put(-26,-14){$s'_3$}
\put(-30,-8){$s''_3$}
\put(0,-13){$s_4$}
\put(-10,11){$s'_4$}
\put(-48,5){$t$}
\put(-39,-9){$t'$}
\put(-19,-6){$u$}
}
\put(11,1){$a$}
\put(-49,-4){$b$}
\end{picture}
& 
$\overset{T}{\longmapsto}$ \ \ 
& 
\begin{picture}(39,15)(-23,7.5)
\put(-25,0){$\underbrace{\hspace{2cm}}_{k_1+k_2}$}
\put(1,15){$\overbrace{\hspace{1.2cm}}^{2k_3}$}
\put(17,15){$\overbrace{\hspace{1.3cm}}^{k=k_4+k_5}$}
\put(-26,0){\line(1,0){50}}
\put(-26,5){\line(1,0){60}}
\put(-26,10){\line(1,0){60}}
\put(-16,15){\line(1,0){50}}
\put(-26,0){\line(0,1){10}}
\put(-16,0){\line(0,1){5}}
\put(-16,10){\line(0,1){5}}
\put(0.5,5){\line(0,1){5}}
\put(16.5,5){\line(0,1){5}}
\put(24,0){\line(0,1){5}}
\put(24,10){\line(0,1){5}}
\put(34,5){\line(0,1){10}}
\put(-15,1){$\scriptstyle \overline{n-1}$}
\put(-7,1){$\scriptstyle \overline{n-1}$}
\put(1,1){$\scriptstyle \overline{n-1}$}
\put(9,1){$\scriptstyle \overline{n-1}$}
\put(17,1){$\scriptstyle \overline{n-1}$}
\put(-25,1){$\scriptstyle \overline{n}$}
\put(-20,1){$\scriptstyle \overline{n}$}
\put(-15,11){$\scriptstyle n-1$}
\put(-7,11){$\scriptstyle n-1$}
\put(1,11){$\scriptstyle n-1$}
\put(9,11){$\scriptstyle n-1$}
\put(17,11){$\scriptstyle n-1$}
\put(-25,6){$\scriptstyle n$}
\put(-20,6){$\scriptstyle n$}
\put(-13,6){$\scriptstyle n$}
\put(-4.5,6){$\scriptstyle n$}
\put(3,6){$\scriptstyle \overline{n}$}
\put(10.5,6){$\scriptstyle n$}
\put(19.5,6){$\scriptstyle \overline{n}$}
\put(25,6){$\scriptstyle \overline{n}$}
\put(30,6){$\scriptstyle \overline{n}$}
\put(25,11){$\scriptstyle n$}
\put(30,11){$\scriptstyle n$}
\put(35,11){$\scriptstyle a$}
\put(-30,1){$\scriptstyle b$}
\put(-20,-13){$a \succeq \overline{n}$, $b \preceq n$}
\put(-20,-19){$k_1+k_2$: even, $k_4+k_5$: odd} 
\put(-20,-25){$k_2 \ne 0$, $k_4 \ne 0$}
\end{picture}
\end{tabular}
}
\vspace{2.8cm}
\caption{Examples of $\bp\in P_0(C_n; \mu, \lambda)$ 
which satisfy one of the conditions of $\bp\in \Image \, f_1^{12}$ 
in Lemma \ref{lem:f1-map} \eqref{item:conditions} and 
the corresponding subtableau in $T(\bp)$.
If the $(E)$ step for $a$ (resp.\ $b$) exists, 
then $a$ (resp.\ $b$) satisfies the condition as above.}\label{fig:condition-f-1-12}
\end{figure}
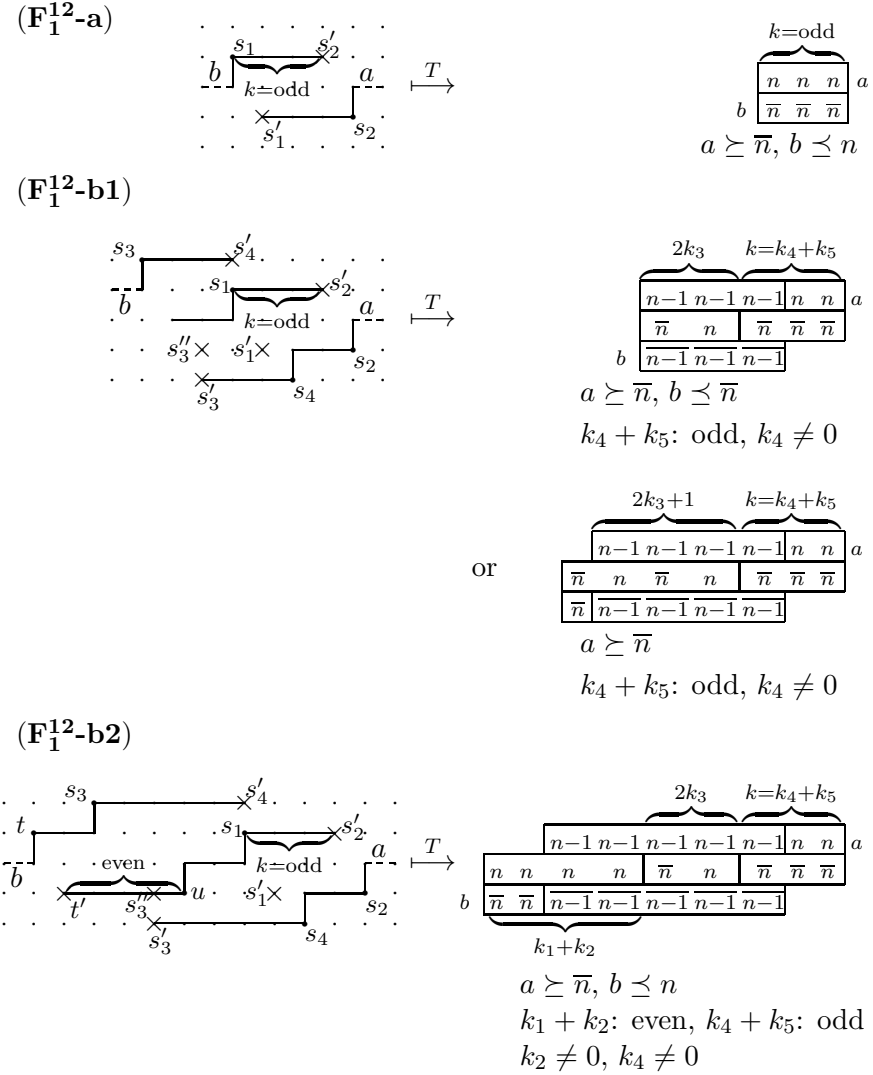

Next, we give the conditions 
to describe the image of $f_1^{12}$ and $f_1^{23}$. 
For any $\bp \in P_0(C_n; \mu, \lambda)$, 
let $s_i$ and $s'_i$ ($i=1,\dots ,4$)
be the points as in 
the conditions ($\bFaa$) and ($\bFab$),
with the roles of $p_2$ and $p_3$ interchanged. 
Namely, (see Figure \ref{fig:condition-f-1-12})
let $s_1$ be the leftmost point of height 1 on $p_2$, 
$s_2$ be the rightmost point of height $-1$ on $p_1$, $s_3$ be 
the leftmost point of height 2 on $p_3$, and $s_4$ be the rightmost point of 
height $-2$ on $p_1$, and set $s'_i:= s_i + (y,-2y)$, where $y$ is the height of $s_i$. 
If $s'_2$ is on $p_2$, then let $k$ be the number of steps of $p_2$ between 
$s_1$ and $s_2$. Then define conditions 
($\bFaa$) and ($\bFab$) for $\bp \in P_0(C_n; \mu, \lambda)$ 
similar to the conditions ($\bFFaa$) and ($\bFFab$), 
with the roles of $p_2$ and $p_3$ 
in the conditions interchanged. Namely, 
\vspace{5pt}
\newline 
\begin{tabular}{ll}
($\bFaa$) & \ $\bp$
satisfies all of the following conditions:
\end{tabular}
\begin{itemize}
\item $s'_1$ is on $p_1$.
\item $s'_2$ is on $p_2$ and $k$ is odd.
\end{itemize}
\begin{tabular}{ll}
($\bFab$) & \ $\bp$
satisfies all of the following conditions:
\end{tabular}
\begin{itemize}
\item $s'_1$ is not on $p_1$.
\item $s'_2$ is on $p_2$ and $k$ is odd.
\item $s'_3$ is on $p_1$.
\item $s'_4$ is on $p_3$.
\end{itemize}
\vspace{5pt}

We also define conditions ($\bFba$) and ($\bFbb$) 
for $\bp \in P_0(C_n; \mu, \lambda)$ as follows: \vspace{5pt}
\newline 
\begin{tabular}{ll}
($\bFba$) & 
\ $\omega(\bp)$ satisfies the condition ($\bFaa$).
\end{tabular}
\newline
\begin{tabular}{ll}
($\bFbb$) & 
\ $\omega(\bp)$ satisfies the condition ($\bFab$).
\end{tabular}
\vspace{5pt}
\newline
Note that ($\bFaa$) and ($\bFab$) 
(resp.\ ($\bFba$) and ($\bFbb$)) are exclusive with each other.

For any $\bp =(p_1,p_2,p_3)\in P_0(C_n; \mu, \lambda)$, 
let $s_3$ 
be the leftmost point of height 2 on $p_3$ 
(as we defined in ($\bFba$) and ($\bFbb$)) 
and $s''_3:=s_{3} + (2,-3)$. 
Let $t$ be  the leftmost point of height $1$ on $p_3$  
and $t':=t +(1,-2)$. 
Let $u$  be the rightmost point of height $-1$ on $p_2$. 
We define conditions 
{\rm ($\bFabone$)}, {\rm ($\bFabtwo$)},  {\rm ($\bFbbone$)}, and {\rm ($\bFbbtwo$)}
as follows (see Figure \ref{fig:condition-f-1-12}):
\vspace{5pt}
\newline
{\rm ($\bFabone$)} \
$\bp$ satisfies {\rm ($\bFab$)} and $s''_3$ is not on $p_2$. 
\newline 
{\rm ($\bFabtwo$)} \
$\bp$ satisfies {\rm ($\bFab$)}, $s''_3$ is on $p_2$, 
$t'$ is on $p_2$, and the number of 
\phantom{{\rm ($\bFabtwo$)} }
the steps from 
$t'$ to $u$ 
is even. 
\newline
{\rm ($\bFbbone$)} \ 
$\omega(\bp)$ satisfies {\rm ($\bFabone$)}.
\newline
{\rm ($\bFbbtwo$)} \ 
$\omega(\bp)$ satisfies {\rm ($\bFabtwo$)}.
\vspace{5pt}

Then we have 
\begin{lem}\label{lem:f1-map}
\begin{enumerate}
\item \label{item:conditions}
Let $\bp \in P_0(C_n; \mu, \lambda)$. Then,  
$\bp \in \Image \, f^{12}_1$ if and only if
one of the conditions {\rm ($\bFaa$)}, {\rm ($\bFabone$)} and {\rm ($\bFabtwo$)}
is satisfied. 
Similarly, $\bp \in \Image \, f^{23}_1$ if and only if
one of the conditions {\rm ($\bFba$)}, {\rm ($\bFbbone$)} and {\rm ($\bFbbtwo$)} 
is satisfied. 
\item 
$f^{12}_1$ and $f^{23}_1$ are weight-preserving sign-inverting injections. 
\end{enumerate}
\end{lem}

The proof of Lemma \ref{lem:f1-map} is 
similar to that of Lemma \ref{lem:f2-map}.

Finally, we give two lemmas which are used in Section \ref{sec:bijective-T3}. 
\begin{lem}\label{lem:ff}
\begin{enumerate}
\item \label{item:ff-and-f1-map}
$\Image \, f^{12}_1 \cap \Image \, f^{23}_1 = \Image \, (f^{23}_1 \circ f^{13}_2)$. 
\item \label{item:ff}
$\Image \, (f^{12}_1\circ f^{23}_2) = \Image \, (f^{23}_1\circ f^{13}_2)$.
\end{enumerate}
\end{lem}

\begin{proof}
Using the conditions in Lemma \ref{lem:f1-map} \eqref{item:conditions}, 
we can check that 
$\bp \in \Image \, f^{12}_1 \cap \Image \, f^{23}_1$ if and only if 
$\bp$ satisfies one of the following: 
\begin{enumerate}
\renewcommand{\labelenumi}{(\alph{enumi})}
\item $\bp$ satisfies ($\bFaa$) and ($\bFba$).
\item $\bp$ satisfies ($\bFaa$) and ($\bFbb$)
($\Leftrightarrow$ ($\bFaa$) and ($\bFbbone$)).
\item $\bp$ satisfies ($\bFab$) and ($\bFba$)
($\Leftrightarrow$ ($\bFabone$) and ($\bFba$)).
\end{enumerate}
On the other hand, 
$f^{23}_1 \circ f^{13}_2$ is given as 
one of the following cases, by the conditions of $\bp \in \Image \, f_2^{13}$: 
\begin{enumerate}
\item 
$f_2^{13}$ as in Case ($f_2^{13}$-a) and $f_1^{23}$ as in Case ($f_1^{23}$-a). 
\item 
$f_2^{13}$ as in Case ($f_2^{13}$-a) and $f_1^{23}$ as in Case ($f_1^{23}$-b). 
\item 
$f_2^{13}$ as in Case ($f_2^{13}$-b) and $f_1^{23}$ as in Case ($f_1^{23}$-a). 
\end{enumerate}
As in the proof of Lemma \ref{lem:f2-map}, all $\bp\in \Image\, f_2^{13}$ 
of Case ($f_2^{13}$-a) (resp.\ Case ($f_2^{13}$-b)) satisfy 
($\bFFaa$) (resp.\ ($\bFFab$)).  
For the conditions ($\bFFaa$)
and ($\bFFab$) of $\bp \in \Image \, f^{13}_2$ 
turn out to be the conditions ($\bFaa$) and 
($\bFab$) respectively  
after $\bp$ is sent by $f^{23}_1$, 
all $\bp \in \Image \, f^{23}_1 \circ f^{13}_2$ of case (1) satisfy (a), 
while that of (2) satisfy (b) and that of (3) satisfy (c). Thus, 
we obtain 
$\Image \, (f^{23}_1 \circ f^{13}_2) \subset 
\Image \, f^{12}_1 \cap \Image \, f^{23}_1$. 
Conversely, $f^{12}_1 \cap \Image \, f^{23}_1 \subset 
\Image \, (f^{23}_1 \circ f^{13}_2)$ is obvious, 
and therefore, \eqref{item:ff-and-f1-map} is proved. 
The condition of $\Image \, (f^{12}_1 \circ f^{23}_2)$ is given similarly, 
and we obtain \eqref{item:ff}.
\end{proof}

\begin{lem}\label{lem:p-T}
For $\bp \in \tilde{P}(C_n; \mu, \lambda)$, 
$\bp \in \Image \, f^{12}_1 \cup \Image \, f^{23}_1$
if and only if 
$T(\bp)\not\in \Tab (C_n, \lambda/\mu)$,
where $\Tab (C_n, \lambda/\mu)$ is the set of tableaux 
defined in Section \ref{sec:bijective-T3}.
\end{lem}

\begin{proof}
If $\bp \in P_0(C_n; \mu, \lambda)$ satisfies one of the conditions of 
$\bp \in \Image \, f^{12}_1$ or that of $\bp \in \Image \, f^{23}_1$
in Lemma \ref{lem:f1-map} \eqref{item:conditions}, then 
$T(\bp)$ does not satisfy either the extra rule $(\bEtwoR )$ or 
the extra rule $(\bEthreeR)$ 
(see Figure \ref{fig:condition-f-1-12}). 

Conversely, let $\bp \in \tilde{P}(C_n; \mu, \lambda)$
and $T(\bp)\not\in \Tab (C_n, \lambda/\mu)$. 
By Lemma \ref{lem:g-map} \eqref{item:g-Tab}, there does not exist any  
$\bp \in \Image \, g$ such that $T(\bp)$ contains one of the subtableaux 
\eqref{eq:arr-2},  \eqref{eq:arr-3-1}, \eqref{eq:arr-3-2} and  
\eqref{eq:arr-3-3}, and therefore,  $T(\bp) \not\in \Tab(C_n, \lambda/\mu)$ 
implies that $\bp \in P_0(C_n; \mu, \lambda)$,  
by Lemma \ref{lem:g-map} \eqref{item:g-map}. 
By assumption, $T(\bp)$ contains a subtableau $T'$ 
described as in \eqref{eq:arr-2},  \eqref{eq:arr-3-1}, \eqref{eq:arr-3-2} or 
\eqref{eq:arr-3-3} which does not satisfy the extra rule
$(\bEtwoR )$ or $(\bEthreeR)$.  
We can check that $\bp$ satisfies one of the conditions in 
Lemma \ref{lem:f1-map} \eqref{item:conditions}
for all such $T'$. Namely (see Figure \ref{fig:condition-f-1-12}),  
\begin{enumerate}
\item 
If $T'$ is the subtableau \eqref{eq:arr-2} prohibited by 
the extra rule $(\bEtwoR )$, 
then $\bp$ satisfies ($\bFaa$) or ($\bFba$).
\item
If $T'$ is the subtableau \eqref{eq:arr-3-1} prohibited by 
the extra rule $(\bEthreeR)$ and  
\begin{enumerate}
\item If $k_4+k_5$ is odd, $k_4\ne 0$ and $k_2=0$, 
then $\bp$ satisfies ($\bFabone$).
\item If $k_4+k_5$ is odd, $k_4\ne 0$ and $k_2\ne 0$, 
then $\bp$ satisfies ($\bFabtwo$). 
\item If $k_4+k_5$ is odd, $k_4=0$ and $k_2\ne 0$, 
then $\bp$ satisfies ($\bFaa$). 
\item If $k_1+k_2$ is odd, $k_2\ne 0$ and $k_4=0$, 
then $\bp$ satisfies ($\bFbbone$). 
\item If $k_1+k_2$ is odd, $k_2\ne 0$ and $k_4\ne 0$, 
then $\bp$ satisfies ($\bFbbtwo$). 
\item If $k_1+k_2$ is odd, $k_2\ne 0$ and $k_4=0$, 
then $\bp$ satisfies ($\bFba$). 
\end{enumerate}
\item If $T'$ is the subtableau \eqref{eq:arr-3-2} 
(resp.\ \eqref{eq:arr-3-3}) 
prohibited by the extra rule $(\bEthreeR)$, then $\bp$ satisfies ($\bFabone$)
(resp.\ ($\bFbbone$)). 
\end{enumerate}
\end{proof}

\end{document}